\documentclass[10pt,reqno]{article}
\usepackage[all]{xy}
\usepackage{amsmath}
\usepackage{amsthm}
\usepackage{amscd}
\usepackage{amssymb}

\input{epsf.tex}

\newcounter{enumitemp}
\newenvironment{enumeratecontinue}{
  \setcounter{enumitemp}{\value{enumi}}
  \begin{enumerate}
  \setcounter{enumi}{\value{enumitemp}}
}
{
  \end{enumerate}
}

\numberwithin{equation}{section}

\newcommand{\centeredepsfbox}[1]{\centerline{\epsfbox{#1}}}
\newcommand{\nb}[1]{#1\nobreakdash-}

\theoremstyle{definition}
\newtheorem*{definition}{Definition}
\newtheorem*{remark}{Remark}
\newtheorem*{example}{Example}

\theoremstyle{plain}
\newtheorem{theorem}{Theorem}[section]
\newtheorem{proposition}[theorem]{Proposition}
\newtheorem{lemma}[theorem]{Lemma}
\newtheorem{corollary}[theorem]{Corollary}
\newtheorem{claim}[theorem]{Claim}

\newtheorem*{corollary*}{Corollary}
\newtheorem*{theorem*}{Theorem}

\newcounter{remarks}
{\paragraph*{Remarks}\smallskip
     \begin{list}{\arabic{remarks}. }{\usecounter{remarks}%
          \setlength{\leftmargin}{0in}%
          \setlength{\rightmargin}{0in}%
          \setlength{\labelsep}{0pt}%
          \setlength{\labelwidth}{0pt}%
          \setlength{\listparindent}{0pt}%
     }
}
{
\end{list}
}

\DeclareMathOperator\Fr{Fr}
\newcommand\inv\inverse

\DeclareMathOperator{\Ends}{Ends}
\DeclareMathOperator\QI{QI}
\DeclareMathOperator\interior{int}

\DeclareMathOperator\SL{SL}
\DeclareMathOperator\Stab{Stab}
\DeclareMathOperator\CStab{CStab}
\DeclareMathOperator\closure{cl}

\DeclareMathOperator\diam{Diam}
\DeclareMathOperator\image{image}
\DeclareMathOperator\midpoint{mid}
\DeclareMathOperator\link{Lk}
\DeclareMathOperator\Vertices{Verts}
\DeclareMathOperator\VE{\V\E}
\DeclareMathOperator\depth{depth}
\DeclareMathOperator\BS{BS}

\DeclareMathOperator\Isom{Isom}
\DeclareMathOperator\supp{supp}

\newcommand\R{{\mathbf R}}
\newcommand\reals{\R}

\newcommand\hyp{{\mathbf H}}

\newcommand\Z{{\mathbf Z}}

\renewcommand\H{{\mathcal H}}
\newcommand\F{{\mathcal F}}
\newcommand\E{{\mathcal E}}
\newcommand\C{{\mathcal C}}

\newcommand\inject{\hookrightarrow}
\newcommand\Sum{\sum}
\newcommand\infinity{\infty}
\newcommand{\bdy}{\partial}
\newcommand{\from}{\colon}
\newcommand\composed{\circ}
\newcommand\suchthat{\bigm|}
\newcommand\inverse{{-1}}
\newcommand\union{\cup}

\newcommand\abs[1]{\left| #1 \right|}
\newcommand\card[1]{\abs{#1}}
\newcommand\subgroup{<}

\newcommand\Id{{\text{Id}}}
\newcommand\intersect{\cap}

\newcommand\restrict{\bigm|}
\newcommand\semidirect{\rtimes}

\newcommand\cross{\times}

\newcommand\Haus{{\mathcal H}}
\renewcommand\O{{\mathcal O}}

\newcommand\Poincare{Poincar\'e}

\newcommand\csubset[1]{\stackrel{\scriptscriptstyle{#1}}{\subset}}

\newcommand\ceq[1]{\stackrel{\scriptscriptstyle{#1}}{=}}

\newcommand\cintersect[1]{\stackrel{\scriptscriptstyle{#1}}\cap}

\newcommand\cstrict[1]{\stackrel{\scriptscriptstyle{#1}}\subsetneq}

\newcommand\<\langle
\renewcommand\>\rangle
\newcommand\disjunion{{\textstyle\coprod}}

\newcommand\Min{{\min}}

\newcommand\U{{\mathcal U}}
\newcommand\V{{\mathcal V}}
\renewcommand\H{{\mathcal H}}
\newcommand\G{{\mathcal G}}
\DeclareMathOperator\PD{PD}
\DeclareMathOperator\rank{rank}

\newcommand\wt{\widetilde}
\renewcommand\P{{\mathcal P}}
\newcommand\tensor\otimes
\newcommand\forest{F}

\newcommand\Xsub[1]{X_{#1}}
\newcommand\norm[1]{\left\| #1 \right\|}

\newcommand\A{{\mathcal A}}
\newcommand\K{{\mathcal K}}
\newcommand\fin{{\text{fin}}}
\newcommand\homeo\approx
\renewcommand\S{{\mathcal S}}

\newcommand\vp{{\vphantom\prime}}

\begin{document}

\setcounter{tocdepth}{3}

\title{Quasi-actions on trees II: \\ Finite depth Bass-Serre trees}
\author{Lee Mosher\thanks{Supported in part by NSF grant
DMS 0103208.}, 
Michah Sageev,
and Kevin Whyte\thanks{Supported in part by NSF grant DMS 0204576}}
\maketitle

\begin{abstract}
This paper addresses questions of quasi-isometric rigidity and
classification for fundamental groups of finite graphs of groups, under
the assumption that the Bass-Serre tree of the graph of groups has finite
depth. The main example of a finite depth graph of groups is one whose
vertex and edge groups are coarse Poincare duality groups. The main
theorem says that, under certain hypotheses, if $\G$ is a finite graph of
coarse Poincare duality groups then any finitely generated group
quasi-isometric to the fundamental group of $\G$ is also the fundamental
group of a finite graph of coarse Poincare duality groups, and any
quasi-isometry between two such groups must coarsely preserves the vertex
and edge spaces of their Bass-Serre trees of spaces. Besides some simple
normalization hypotheses, the main hypothesis is the ``crossing graph
condition'', which is imposed on each vertex group $\G_v$ which is an
$n$-dimensional coarse Poincare duality group for which every incident
edge group has positive codimension: the crossing graph of
$\G_v$ is a graph $\epsilon_v$ that describes the pattern in which the
codimension~1 edge groups incident to $\G_v$ are crossed by other edge
groups incident to $\G_v$, and the crossing graph condition requires that
$\epsilon_v$ be connected or empty.
\end{abstract}

\section{Introduction}

One important part of geometric group theory is the classification of finitely generated groups up to
quasi-isometry.  This paper addresses that question for many groups which split as graphs of groups.  A
typical question in this context is: suppose $G$ splits as a graph of groups and $H$ is quasi-isometric to
$G$, must $H$ split in a similar way? 

The first example of this type is Stallings Ends Theorem, which
implies that if $G$ splits with finite edge groups, then so must $H$. This result was recently refined in
work of Papasoglu and Whyte \cite{PapasogluWhyte:ends} which implies that if $G$ and $H$ are accessible
groups then they have the same quasi-isometry types of one-ended factors. After Stallings theorem, the next
result of this type came in work of Kapovich and Leeb \cite{KapovichLeeb:haken} which shows that if $M$ is
a Haken $3$-manifold then the splitting of $\pi_1(M)$ induced by the geometric decomposition of $M$ is
preserved by quasi-isometries. Farb and Mosher \cite{FarbMosher:BSTwo} showed that if $G$ splits as a graph
of $\Z$'s with solvable fundamental group --- a solvable Baumslag-Solitar group --- then so must $H$ up to
finite groups. Papasoglu \cite{Papasoglu:Zsplittings} showed that if $G$ splits with two-ended edge groups,
then $H$ must as well. Mosher, Sageev and Whyte \cite{MSW:QTOne} showed that if $G$ splits as a
``bushy'' graph of coarse \Poincare\ duality groups of constant dimension, then so must $H$. 

Not all splittings exhibit this type of rigidity. For instance, cocompact lattices in $\SL(2,\R)\times
\SL(2,\R)$ are all quasi-isometric, however, by a theorem of Margulis the irreducible ones do not split
while reducible ones split as amalgams in a myriad of ways (\cite{Margulis:amalgams}, and see also
\cite{MonodShalom:superrigidity}). 

In this paper, we focus on graphs of groups whose vertex groups are ``manifold-like'', for example
fundamental groups of closed aspherical manifolds, and more generally finitely presented \Poincare\ duality
groups. Since we are dealing with quasi-isometry issues, it turns out that the natural class of groups to
work with are coarse \Poincare\ duality groups, introduced by Kapovich and
Kleiner~\cite{KapovichKleiner:duality}.  In a previous paper \cite{MSW:QTOne} we studied a special case in
which all of the vertex and edge groups have the same dimension, equivalently, every edge-to-vertex
injection in the graph of groups has finite index image; such graphs of groups are called
\emph{(geometrically) homogeneous}. In this paper, we address more general cases.

Papasoglu's preprint \cite{Papasoglu:GroupSplittings} independently studies graphs of groups with a
similar focus and overlapping results to ours; we explain below some of the similarities and differences of
our approaches and results.

The main theorems~\ref{TheoremQI}, \ref{TheoremClasses} and~\ref{PDgraphs} are somewhat involved, so before
launching into a full discussion of them, we begin with some applications, chosen from among the
gamut of settings to which the main theorems apply.\footnote{The main theorems, and several of the
applications, were announced in \cite{MSW:QTannc}.}

\subsection{Example applications}
\label{SectionExampleApplications}

First, here is a theorem regarding two fundamental groups of hyperbolic manifolds amalgamated along
infinite cyclic subgroups. For the proof see Theorem~\ref{GenHyp}, which relies essentially on a theorem of
Schwartz \cite{Schwartz:Symmetric}. 

\begin{theorem}\label{Neat} Let $A$ and $B$ be fundamental groups of closed hyperbolic manifolds
of dimension at least $3$, and let $G$ be the amalgamated free product $A *_\Z B$ where $\Z$
includes as a maximal cyclic subgroup in both $A$ and $B$. If $H$ is a finitely generated group
quasi-isometric to $G$ then $H$ splits as a graph of groups each of whose vertex groups is
commensurable to $A$ or $B$. In particular, there are infinitely many quasi-isometry classes of groups of
this form.
\end{theorem}
 
In fact, in the context of Theorem~\ref{Neat} one obtains a complete computation of the
quasi-isometry group; see Proposition~\ref{PropDehn}. 
 
Notice that in the above example, the vertex-to-edge inclusions are codimension-2. As we will see
below, Alexander duality can be applied in a straightforward  manner. In the codimension-1
setting, one typically needs that the edge subgroups in a vertex group cross one another
in some suitable sense. Here is one such example.
 
A collection of geodesics in the hyperbolic plane is said to be {\it filling} if the complementary
regions are all bounded. A collection of cyclic subgroups in a surface group is said to be {\it
filling} if the pattern of geodesics in the universal cover is filling. 

\begin{theorem}\label{Filling} Suppose that $G$ splits as a graph of groups where each vertex
group is a surface group, each edge group is a cyclic group and the edge-to-vertex inclusions
provide a filling collection of curves on the vertex group. Then any torsion-free group
quasi-isometric to $H$ splits as a graph of groups whose edge groups are cyclic and whose vertex
groups are surface groups or cyclic groups. 
\end{theorem}

A stronger version of Theorem~\ref{Filling} is given in Theorem~\ref{TheoremStrongSurface}, with
commensurability information built into the conclusion, by applying a theorem of Kapovich and
Kleiner~\cite{KapovichKleiner:LowDBoundaries}. Using this we will obtain infinitely many
distinct quasi-isometry classes of groups of this type. 

The next theorem gives examples with higher dimensional edge groups, namely fundamental groups of surfaces.
For the proof see Theorem~\ref{FiberedManifolds}, and for even stronger rigidity conclusions in many
examples see Theorem~\ref{StrongFibered}.

\begin{theorem}\label{Fibered} Let $A$ be the fundamental group of a hyperbolic $3$-manifold which
fibers over the circle in two ways, with fibers $F_1$ and $F_2$.  Let $\phi$ be any
isomorphism between (finite index subgroups of) $\pi_1(F_1)$ and $\pi_1(F_2)$ and define
$G$ as the HNN extension $A*_\phi$.  If $H$ is any torsion free group quasi-isometric to $G$, then
$H$ splits as a graph of groups, whose edge groups are surface groups, and whose vertex groups are 
hyperbolic 3-manifold groups commensurable to $A$ and surface groups.
\end{theorem}

The commensurability information in Theorem~\ref{Fibered} is an application of a theorem of Farb and
Mosher~\cite{FarbMosher:sbf}. 

Here is a theorem for abelian groups. 

\begin{theorem} 
\label{TheoremBabyAbelian}
Suppose $A$ is a finitely generated abelian group of rank $n$, $G_1, G_2\subgroup A$ are two rank $n-1$
subgroups that span $A$, and $G$ is an HNN-extension identifying $G_1$ to $G_2$. Then any group $H$
quasi-isometric to $G$ splits as a graph of virtually abelian groups with edge groups of rank $n-1$ and
vertex groups of rank at most $n$.
\end{theorem}

A much more general result is given in Corollary~\ref{abrigid}, applying to a broad spectrum of graphs of
abelian groups.

\subsection{The methods of proof: a special case.} 
\label{SectionMethods}
 
We illustrate some of the ideas involved in proving our main theorems by sketching the proof
in the special  case of theorem \ref{Neat}. Let $M$ and $N$ be closed hyperbolic $n$-manifolds for
some $n \geq 3$.  Let $A$ and $B$ be their fundamental groups.   Choosing maximal cyclic subgroups
of $A$ and $B$ amounts to choosing primitive closed geodesics in $M$ and $N$.   Let $Y$ be the
space built by gluing an annulus to $M$ and $N$, with the boundary circles attached to the chosen
geodesics.  The group $G = A *_\Z B$ is
$\pi_1 (Y)$, and so acts properly discontinuously and cocompactly on $X = \wt{Y}$.  The space
$X$ is built from copies of $\hyp^n$ glued together with strips attached along the geodesics which
cover the chosen curves in $M$ and $N$.   There is a natural $G$-equivariant map, $\pi: X\to T$,
from $X$ to the Bass-Serre tree $T$ of the splitting of $G$ in which each copy of $\hyp^n$ is
mapped to a single vertex, and each strip maps to an edge.  The is an example of what we call a
{\em Bass-Serre complex}, which serves as the standard model space for our considerations. (See
Section~\ref{SectionBassSerre} for the general construction.)  We call the copies of $\hyp^n$ in
$X$ the vertex spaces $X_v$, one for each vertex $v \in T$. The strips are called the edge
spaces $X_e$, one for each edge $e \subset T$.

We start by understanding the self quasi-isometries of $G$. This is typically a crucial step in
proving rigidity results. The salient point is that a quasi-isometry from some mystery group $H$
to $G$ provides a natural map from  $H$ to $\QI(G)$, the quasi-isometry group of $G$. Thus one
seeks geometric patterns which are invariant under $\QI(G)$. As $X$ is quasi-isometric to $G$, we
can work with quasi-isometries of $X$. Since $X$ is uniformly contractible, we can move any
quasi-isometry a bounded distance to a continuous map. Let $f \from X \to X$ be such a continuous
quasi-isometry. Our method of proof consists of two main steps.

\paragraph{Step 1 -- Vertex Rigidity:} There is an $R>0$ such that for any vertex space $X_v
\subset X$ there is a vertex space $X_{v'} \subset X$ with $d_\Haus(f(X_v),X_{v'}) \leq R$, where
$d_\Haus(\cdot,\cdot)$ denotes the Hausdorff metric.

\paragraph{Step 2 -- Tree Rigidity:} There is a automorphism $\tau$ of $T$ and a
quasi-isometry $f' \from X \to X$ at bounded distance from $f$ such that $\pi \circ f' = \tau \circ
\pi$.

\paragraph{Conclusion:} If we establish steps 1 and 2, then we may conclude that if $H$ is
quasi-isometric to $G$ then $H$ acts on $T$ with finite quotient;  for every vertex $v$ of $T$,
$H_v$ is quasi-isometric to $G_v$ and similarly, for every edge $e$, $H_e$ is quasi-isometric to
$G_e$.  In particular, $H$ has a splitting as a graph of groups with virtually cyclic edge groups
and vertex groups quasi-isometric to $\hyp^3$. A theorem of Schwartz \cite{Schwartz:Symmetric}
then tells us that the vertex and edge groups of the splitting for $H$ are
commensurable to the corresponding ones for $G$.

Let us now see how steps 1 and 2 are proven in this setting. 

\paragraph{Vertex rigidity.} Suppose that vertex rigidity fails. Then one finds an edge space
$X_e$ such that $f(X_v)$ travels deeply into both complementary components of $X_e$. This will
implies that $X_e$  (which we recall is quasi-isometric to a line) coarsely separates $f(X_v)$.
But now note that $f(X_v)$ is quasi-isometric to $\hyp^n$, so we get a subset of a line
coarsely separating $\hyp^n$ with $n\geq 3$. This contradicts a coarse version of Alexander
duality.

\paragraph{Tree rigidity.} The vertex rigidity step yields a bijection $f_\# \from \Vertices(T)\to
\Vertices(T)$, so that $f(X_v)$ has finite Hausdorff distance from $X_{f_\#{(v)}}$. We now wish to see that
$f$ carries neighboring vertex spaces to neighboring vertex spaces. The key point here is that in $X$,
adjacent vertex spaces are characterized coarsely by the fact that they have neighborhoods with unbounded
intersection, whereas any intersection of neighborhoods of non-adjacent vertex spaces is bounded. This fact
relies on the assumption that the cyclic edge groups are maximal in their respective vertex
spaces. Thus $f_\#$ carries adjacent vertices to adjacent vertices and hence extends to an
isometry $\tau$ of $T$ as required.

\subsection{The general setting}
\label{SectionSetting}

The proofs of our main theorems for more general graphs of groups follow a similar outline. Let us
discuss some of the issues that arise in the more general setting, and the various hypotheses
needed to deal with these issues.

For starters, our graphs of groups $\G$ will always satisfy the standard hypothesis of
\emph{irreducibility}, meaning that for every edge of $\G$ with distinct endpoints, both of its
edge-to-vertex injections have image of index $\ge 2$. Every finite graph of groups can be made
irreducible by inductively collapsing edges, without changing the fundamental group. 

Let $T$ denote the Bass-Serre tree of $\G$, and let $X$ denote the Bass-Serre complex of $\G$. The group
$\pi_1\G$ acts on $X$ and $T$, and there is a $\pi_1\G$ equivariant projection map $X \to T$. The action on
$X$ is properly discontinuous and cocompact, and the action on $T$ is cocompact with quotient  graph of
groups $\G$. The inverse image under the projection map $X \to T$ of a vertex $v\in T$ is a \emph{vertex
space} denoted $X_v$, and the inverse image of the midpoint of an edge $e\subset T$ is an \emph{edge space}
denoted~$X_e$.

\paragraph{Vertex Rigidity.} The proof of the vertex rigidity step in Section~\ref{SectionMethods}
above relied on the fact that no geodesic can coarsely separate $\hyp^n$ when $n \ge 3$.   In
general we want to know that edge spaces cannot disconnect vertex spaces.  In most applications
the essential focus is on the maximal vertex spaces, those not strictly coarsely contained in any
other vertex space.  We call these the {\em depth zero} vertex spaces. An important hypothesis of
our results will be that every vertex space is coarsely contained in a depth zero vertex space. 

To see that edge spaces cannot coarsely separate the depth zero vertex spaces we generally use a
coarse version of Alexander duality to understand the separation properties of subsets.  The
right context for these arguments is that of coarse Poincare duality spaces introduced in
\cite{KapovichKleiner:duality}. Almost all of our results assume that the depth zero vertex
spaces are coarse $\PD$ spaces (see section \ref{SectionPD} for detailed definitions). To put it
another way, the stabilizer groups of depth zero vertices are required to be ``coarse $\PD$
groups''. This is an essential hypothesis for all of our results.

There are several different vertex rigidity phenomena. For example, given a coarse $\PD(n)$
vertex space, if the incident edge spaces are of dimension at most $n-2$ then we can use a version
of the vertex rigidity step of Section~\ref{SectionMethods}: in a coarse $\PD(n)$ space, no
subspace whose ``coarse codimension'' is $\ge 2$ can coarsely separate the ambient space. We must
also study vertex rigidity phenomena in the presence of incident edge spaces of coarse codimension
zero and one. In order to make sense of these different phenomena, for example in order for the
coarse codimension of an edge space in a vertex space to be well-defined, we will add a hypothesis
that the edge spaces satisfy a coarse version of the finite type property. In other words, the
stabilizers of edges (and also of positive depth vertices) are required to be ``coarse finite type
groups'' (see section \ref{SectionPD}). The coarse finite type hypothesis is a convenient assumption for our
results, although not essential like the coarse $\PD$ assumption; see the comments after the
statement of Proposition~\ref{PropEdgeSeparatesVertex}.

Continuing with a discussion of vertex rigidity phenomena, a codimension zero inclusion corresponds
to an edge group which has finite index in the vertex group.  A \emph{depth zero raft} $R$ is a
maximal subtree of the Bass-Serre tree whose vertices have depth zero, and whose edge stabilizers
include with finite index in incident vertex stabilizers. The union of vertex and edge spaces
corresponding to vertices and edges in $R$ is a portion of the Bass-Serre complex $X$ called a
\emph{raft space} denoted $X_R$. A depth zero raft $R$ is a Bass-Serre tree for a geometrically homogeneous
graphs of groups, whose fundamental group is the stabilizer of the raft, and the raft space $X_R$ is the
corresponding Bass-Serre complex. We can use our previous results \cite{MSW:QTOne} to study depth
zero rafts.  Those results apply to give vertex rigidiy for a geometrically homogeneous graph of
coarse $\PD(n)$ groups, but \emph{only} when its Bass-Serre tree is not a line. That means we must
add a hypothesis that no depth zero raft is a line, equivalently, each depth zero raft is either a
point or a bushy tree (see Proposition~\ref{PropTrichotomy}).

Codimension one edges also require certain additional hypotheses. To explain these assumptions,
consider the group $F_2 \cross \Z \approx \<b,t\> \cross \<a\>$, which splits as an HNN
amalgamation of $\Z\oplus\Z \approx \<a,b\>$ over a cyclic group, with stable letter $t$ --- that
is, $\<a,b,t \suchthat [a,b] = 1,  a^t = a\>$. The Bass-Serre complex $X$ of this splitting does
not satisfy vertex rigidity: the automorphism $b\mapsto t$, $t \mapsto b$, $a\mapsto a$ fails to
respect the vertex spaces of $X$. This failure stems from the fact each edge space $X_e$ coarsely
separates each incident vertex space $X_v$ in $X$ --- that is, for any sufficiently large
neighborhood of $X_e$ in $X$, the set $X_v$ enters deeply into two distinct complementary
components of this neighborhood. 

In order to avoid this problem, we introduce the \emph{crossing graph condition}. Consider a coarse
$\PD(n)$ vertex space $X_v$ of depth zero which is not coarsely equivalent to any other vertex space;
equivalently, the vertex $v$ is a depth zero raft of the Bass-Serre tree. Define the
\emph{crossing graph} of $X_v$ to be a graph $\epsilon_v$ with one vertex for each coarse
$\PD(n-1)$ edge space incident to $X_v$, and with edges as follows. If $X_e$ is a coarse $\PD(n-1)$ edge
space incident to $X_v$, then a subset $A \subset X_v$ is said to \emph{cross $X_e$ in $X_v$} if $A$ goes
deeply into both complementary components of $X_v - X_e$; see Sections~\ref{SubsectionLanguage}
and~\ref{SectionAbstractVertexRigidity} for details on crossing. Consider two vertices $X_e$,
$X_{e'}$ of $\epsilon_v$. We connect these vertices by an edge of $\epsilon_v$ if one of two things
happens: either there is an another edge $e''$ incident to $v$ (with no restriction on the dimension of
$X_{e''}$) such that $X_{e''}$ crosses both $X_e$ and $X_{e'}$ in $X_v$; or $X_e$ and $X_{e'}$ cross each
other in $X_v$. 

\begin{description}
\item[Crossing graph condition:] For each $n$, and for each depth zero vertex $v$ which is a raft and whose
vertex space $X_v$ is coarse $\PD(n)$,
\begin{enumerate}
\item\label{ItemCrossing1}
For any edge $e$ incident to $v$, if the edge space $X_e$ coarsely separates $X_v$ then $X_e$ is coarse
$\PD(n-1)$.
\item The crossing graph $\epsilon_v$ is connected or empty.
\end{enumerate}
\end{description}

Crossing graphs have been used before in geometric group theory. For example, in the proof of
the algebraic torus theorem of Dunwoody and Swenson \cite{DunwoodySwenson:torus}, a crossing graph similar
to ours is constructed and its connectivity properties are used in the proof, although their crossing graph
allows an edge only when two codimension~1 subspaces cross each other and not when they are both crossed by
a third subspace.

Here are two good examples to keep in mind. First, for vertex group $\Z^n$, connectedness of the
crossing graph condition is equivalent to the assumption that the edge groups rationally span the
vertex group.  The second example is a surface group as the vertex group, with edge groups $\Z$. 
Connectedness of the crossing graph means that the curves corresponding to the edge groups fill
the surface.

To summarize, here are the key hypotheses designed to guarantee depth zero vertex rigidity, as given in
our general depth zero vertex rigidity theorem~\ref{TheoremVertexRigidity}:
\begin{itemize}
\item All depth zero vertices are coarse $\PD$.
\item No depth zero raft is a line.
\item The crossing graph condition holds.
\end{itemize}
Another convenient hypothesis,
designed to make various concepts of dimension well defined, is:
\begin{itemize}
\item Each vertex and edge group is coarse finite type.
\end{itemize}

\paragraph{Tree Rigidity.} In the proof of Theorem \ref{Neat} we identified the tree using the coarse
intersections of the vertex spaces.  In general this is not possible. For example, in the free product of
two groups $G = A*B$, the coarse intersection of any pair of vertex spaces is bounded, so there is no
direct way to reconstruct the tree.  Indeed, $G$ has many quasi-isometries which do not 
arise from tree automorphisms.  A similar situation occurs if one tries to extend Theorem
\ref{Neat} to cover amalgamations over non-maximal cyclic subgroups. In order to control the complexity
of the coarse inclusion lattice of edge spaces, we introduce the following:
\begin{description}
\item[Finite depth hypothesis] There is a bound on the length of a chain of strict coarse
inclusions of edge spaces. 
\end{description} 
When this hypothesis is satisfied we say that $\G$, or its Bass-Serre tree $T$, have \emph{finite depth}.
See Section \ref{SectionFiniteDepth} for details.  It turns out that in the setting of graphs of
groups with coarse $\PD$ vertex and edge groups, the finite depth hypothesis is always satisfied.
When the depth zero vertices are required to be coarse $\PD$, but not the remaining edge and
vertex groups, then the finite depth hypothesis may or may not be satisfied.

The finite depth hypothesis can be thought of as a way to tame the coarse inclusion lattice of
edge spaces, and it plays a key role in reconstructing a tree structure, as needed to prove tree
rigidity.

The general statement of tree rigidity is given in Theorem~\ref{TheoremTreeRigidity}.

Here is an intriguing example, which highlights the limitations of our techniques. Let $K$ be a
$\PD(n)$ group with a free subgroup, let $\phi \from A \to B$ be an isomorphism between two free
subgroups of $K$ of rank $\ge 2$, and consider the HNN amalgamation $G = K *_\phi$. Let $T$ be the
Bass-Serre tree of $G$.  Suppose that $B$ has infinite index in $A$, or suppose that one of the
subgroups $A$ or $B$ is conjugate in $K$ to an infinite index subgroup of itself. In this case
there exists a bi-infinite sequence of edges $\ldots,e_{-1},e_0,e_1,\ldots$ in $T$ such that if
$G^{e_i}$ denotes the subgroup of $G$ stabilizing $e_i$ then each $G^{e_i}$ is an infinite
index subgroup of $G^{e_{i+1}}$. In other words, the Bass-Serre complex has infinite depth. Our
techniques do not apply in this situation: if $H$ is quasi-isometric to $G$ then we cannot
conclude using our techniques that $H$ splits as the fundamental group of a graph of groups with
coarse $\PD(n)$ vertex groups and virtually free edge groups. Nonetheless, there are two
partial results in this situation: our vertex rigidity theorem~\ref{TheoremVertexRigidity}
\emph{does} apply provided $n \geq 3$, and it shows that $H$ has subgroups corresponding to the
vertex groups of $G$; and on the other hand Papasoglu's results in
\cite{Papasoglu:GroupSplittings} apply to show that $H$ splits over a virtually free group. We
have not investigated how to combine these results; in particular, we do not know how to get
complete control over a graph of groups for $H$ in the manner of our main quasi-isometric rigidity
theorem~\ref{TheoremQI}.

\subsection{Statements of results} 
\label{SectionStatements}

We now give full statements of our main results on quasi-isometric rigidity and classification, together
with an application to graphs of coarse $\PD$ groups that generalizes all of the examples of
Section~\ref{SectionExampleApplications}. The graphs of groups $\G$ to which our theorems apply are
characterized by the following hypotheses:

\begin{enumerate}
\item[(1)] $\G$ is finite type, irreducible, and of finite depth.
\item[(2)] No depth zero raft of the Bass-Serre tree of $\G$ is a line.
\item[(3)] Every depth zero vertex group of $\G$ is coarse $\PD$.
\item[(4)] For every one vertex, depth zero raft of the Bass-Serre tree, the crossing graph
condition holds.
\item[(5)] Every vertex and edge group is coarse finite type.
\end{enumerate}

\begin{theorem}[Quasi-isometric rigidity theorem]
\label{TheoremQI}
Let $\G$ be a graph of groups satisfying (1)--(5) above. If $H$ is a finitely generated group
quasi-isometric to $\pi_1\G$ then $H$ is the fundamental group of a graph of groups satisfying
(1--5) above.
\end{theorem}

In studying the large scale geometry of groups, besides proving that a particular class of
groups is closed under quasi-isometry (rigidity), one wants to obtain a complete
classification up to quasi-isometry of the groups in that class. In the present situation, our
class of groups is too broadly defined to expect a complete classification. We will settle for
reducing the classification of the whole group to classification of the vertex and edge groups
and the pattern in which those groups are pieced together. There are many examples where one
can go further and get a complete classification of a specific class of groups that is
closed up to quasi-isometry: in the geometrically homogeneous case see
\cite{FarbMosher:BSOne}, \cite{Whyte:bs}, \cite{FarbMosher:ABC}, \cite{FarbMosher:sbf},
\cite{Whyte:FiberedGeometries}; and in the inhomogeneous case see Theorem~\ref{HypManifolds}
above and other examples in Section~\ref{Applications}.

To state the Classification Theorem, consider a graph of groups $\G$ with Bass-Serre tree of
spaces $X\to T$. Let $\VE(T) = V(T) \union E(T)$, where the usual simplicial metric on $T$
induces a metric on $\VE(T)$ by representing each edge $e \in E(T)$ by its midpoint.

\begin{theorem}[Quasi-isometric classification theorem]
\label{TheoremClasses}
Let $\G,\G'$ be graphs of groups satisfying (1--5) above. Let $X \to T$, $X' \to T'$ be
Bass-Serre trees of spaces for $\G,\G'$, respectively. If $f \from X \to X'$ is a quasi-isometry
then $f$ coarsely respects vertex and edge spaces. To be precise, for any $K \ge 1$, $C \ge 0$
there exists $K'\ge 1$, $C' \ge 0$ such that if $f \from X \to X'$ is a $K,C$ quasi-isometry then
there exists a $K',C'$ quasi-isometry $f_\# \from \VE(T) \to \VE(T')$ such that the following hold:
\begin{itemize}
\item If $a \in\VE(T)$ then $d_\Haus(f(X_a),X'_{f_\#(a)}) \le C'$.
\item If $a'\in \VE(T')$ then there exists $a \in \VE(T)$ such that $d_\Haus(f(X_a),X'_{a'})
\le C'$.
\end{itemize}
\end{theorem}

The finite depth hypothesis insures that $T,T'$ each have a nonempty collection of depth zero
rafts, and so in particular this theorem says that $f$ induces a bijection between the depth zero
rafts of $T$ and of $T'$.

The conclusion of Theorem~\ref{TheoremClasses} embraces several further conclusions. For
example, if $a \in \VE(T)$ and $a' = f_\#(a) \in \VE(T)$ then $f$ induces a quasi-isometry
between the subgroup of $\pi_1\G$ stabilizing $a$ and the subgroup of $\pi_1\G'$ stabilizing
$b$. Somewhat deeper implications arise from considering adjacencies of vertex and edge
spaces, for example the conclusion about edge group patterns in a vertex group in
Corollary~\ref{HypManifolds}. We will explore this more in Section~\ref{Applications} under
the moniker of ``pattern rigidity''.

\paragraph{Graphs of coarse \Poincare\ duality groups.} If $G$ is a coarse $\PD(n)$ group and $H \subset
G$ is a coarse $\PD(k)$ subgroup then $k \le n$, with strict inequality if and only if $H$ has
infinite index in $G$. This is well known fact for ordinary \Poincare\ duality groups
\cite{Brown:cohomology}, and the generalization to coarse \Poincare\ duality groups is due to Kapovich and
Kleiner (\cite{KapovichKleiner:duality}, see also Lemma~\ref{LemmaBigInSmall}). 

Suppose now that $\G$ is any finite graph of coarse $\PD$ groups. The fact just noted about dimensions of
subgroups shows that $\G$ has finite depth. We can therefore apply Theorems~\ref{TheoremQI}
and~\ref{TheoremClasses} to $\G$, obtaining the following theorem which applies, for example, to all of the
examples considered in Section~\ref{SectionExampleApplications}:

\begin{theorem}[QI Rigidity and Classification for graphs of coarse $\PD$ groups]
\label{PDgraphs}
Let $\G$ be a finite, irreducible graph of coarse \Poincare\ duality groups satisfying (2) and (4): no
depth zero raft of the Bass-Serre tree $T$ is a line; and the crossing graph condition holds for every one
vertex depth zero raft of $T$. If $H$ is any finitely generated group quasi-isometric to $\pi_{1}\G$ then
$H$ is the fundamental group of a finite type, irreducible graph of coarse \Poincare\ duality groups
$\G'$ satisfying (2) and (4). Furthermore, for any such $\G'$, any quasi-isometry $\pi_1 \G \to
\pi_1\G'$ coarsely respects vertex and edge spaces.
\end{theorem}

As explained at the end of Section~\ref{MainTheorems}, two of the hypotheses --- irreducibility, and no
depth zero raft is a line --- are normalization hypotheses: given any finite graph of
$\PD$ groups, the graph can be normalized in a canonical way so as to make these two hypotheses true,
without changing the fundamental group, and changing the vertex and edge stabilizers in a predictable way.
The most important and restrictive hypothesis is therefore the crossing graph hypothesis. The theorem
applies to any (normalized) finite graph of coarse \Poincare\ duality groups subject to the crossing graph
hypothesis.

There are other situations outside of the coarse $\PD$ world to which Theorems~\ref{TheoremQI}
and~\ref{TheoremClasses} apply, although verifying the finite depth hypothesis can become tricky; see the
end of Section~\ref{SectionHnZ} for one such class of examples.

\paragraph{Relations with work of Papasoglu.} Our results have some overlap with Papasoglu's paper
\cite{Papasoglu:GroupSplittings}, but with somewhat different goals. In particular, many of the special
examples in Section~\ref{SectionExampleApplications} above are handled in Papasoglu's paper as well.

The main difference between the two papers is that the results of \cite{Papasoglu:GroupSplittings} are
focussed not on the overall structure of a graph of groups, but on the question of whether a group splits
with a given type of edge group. As such, the hypotheses of the results in \cite{Papasoglu:GroupSplittings}
are generally weaker than ours, and so have wider applicability. On the other hand, the conclusions are also
weaker. 

To compare the techniques, the use of coarse Alexander duality to understand vertex groups is similar here
and in \cite{Papasoglu:GroupSplittings}. On the other hand, the splittings in
\cite{Papasoglu:GroupSplittings} are obtained in a different manner, by applying work of Dunwoody and
Swenson \cite{DunwoodySwenson:torus}.

\subsection{Structure of the paper}

In Section~\ref{Preliminaries} we introduce the basic notions of coarse geometry, Bass-Serre
theory, and the coarse separation properties which we will need, including a discussion of coarse
\Poincare\ duality groups. All concepts introduced in this introduction are carefully developed in
Section~\ref{Preliminaries}.

Also, in Section~\ref{SectionMethodsGeneral} we introduce the overall scheme of the proofs of the
main theorems. The proofs are broken into three steps. The first step is the Depth Zero Vertex
Rigidity Theorem~\ref{TheoremVertexRigidity}. The second step, which was suppressed in our earlier
discussion, is the Vertex--Edge Rigidity Theorem~\ref{TheoremVERigidity}; this result gives
rigidity for vertex and edge spaces of arbitrary depth, starting from a hypothesis of depth zero vertex
rigidity. The third step is the Tree Rigidity Theorem~\ref{TheoremTreeRigidity}.

In Section~\ref{SectionDepthZero} we prove Depth Zero Vertex Rigidity~\ref{TheoremVertexRigidity}.
First we prove rigidity for the depth zero rafts, and then we apply geometrically homogeneous methods
(e.g.\ \cite{FarbMosher:ABC} or \cite{MSW:QTOne}) to establish vertex rigidity within the rafts.

In Section~\ref{SectionFiniteDepth} we explore the coarse lattice of vertex and edge spaces, with
respect to the operations of coarse set theory, in order to prove Vertex---Edge
Tigidity~\ref{TheoremVERigidity}. This coarse lattice structure organizes the coarse equivalence
classes of vertex and edge spaces into rafts of higher depth. 

Section~\ref{TreeRigidity} proves the Tree Rigidity Theorem~\ref{TheoremTreeRigidity}. With the
results of the previous sections already established, at this point we may conclude that
quasi-isometries on the Bass-Serre complex descend to quasi-isometries on the Bass-Serre tree.
Thus, a mystery group $H$ quasi-isometric to the Bass-Serre complex $X$ quasi-acts on the
Bass-Serre tree $T$. In Section~\ref{TreeRigidity} we show how to use the coarse lattice
structure, together with Dunwoody's theory of tracks, in order to produce an action of $H$ on a
tree $T'$ quasiconjugate to the quasi-action of $H$ on $T$. This gives the desired splitting of
$H$.

In Section~\ref{MainTheorems} we put together the results of Sections~\ref{SectionDepthZero},
\ref{SectionFiniteDepth}, and~\ref{TreeRigidity} to prove the main quasi-isometric rigidity
Theorem~\ref{TheoremQI}. We also prove the main quasi-isometric classification
Theorem~\ref{TheoremClasses}, using only the results of Section~\ref{SectionDepthZero}
and~\ref{TreeRigidity}; tree rigidity is not needed. We also apply these results to a very broad spectrum
of graphs of coarse \Poincare\ duality groups; see Theorem~\ref{PDgraphs}.

Finally in Section~\ref{Applications} we restrict our discussion to some interesting special
cases, including all of the example applications from above. Together with some pattern rigidity results
(including deep work of R.\ Schwartz), we obtain some stronger quasi-isometric rigidity results.

\paragraph{Bottlenecks.} The bottleneck for further extension of our techniques is depth zero
vertex rigidity. The results of sections~\ref{SectionFiniteDepth} and~\ref{TreeRigidity} show that
the conclusions of the main QI-rigidity and classification theorems are formal consequences of
depth zero vertex rigidity (when working with finite type, finite depth graphs of groups). The
specialized hypotheses of the main QI-rigidity and classification theorems are designed
specifically to guarantee depth zero vertex rigidity. 

Certainly there are other situations, not covered by our hypotheses, in which depth zero vertex
rigidity holds. For example, consider a graph of groups with $\hyp^2$ vertex groups,
and whose depth one edge groups are free of rank $\ge 2$. These would be nice examples to
contemplate, working towards QI-rigidity and classification for all graphs of groups in which each
depth zero raft is a single vertex with an $\hyp^2$ vertex group.

\vfill\break

\tableofcontents

\section{Preliminaries}
\label{Preliminaries}

\subsection{Coarse language}
\label{SubsectionLanguage}

In this section we introduce the language of coarse geometry: Hausdorff distance of subsets,
uniformly proper and quasi-isometric embeddings of metric spaces, quasi-actions, etc. Many of
these concepts come equipped with explicit constants, for example we speak of a $K,C$
quasi-isometry with a multiplicative constant $K\ge 1$ and an additive constant $C \ge 0$.
The notation of coarse language can be used with or without explicit reference to the
constants involved; the context should make clear what the constants are and what roles they
play. But we will be careful to distinguish between multiplicative constants and additive
constants; this contrasts with a common practice of referring to a $K$ quasi-isometry where
$K$ stands for both the multiplicative and the additive constant. On the other hand, we will
often take the liberty to conflate several additive constants into one.

\paragraph{Coarse set theory.} Let $X$ be a metric space. Given $A \subset X$ and $r \ge 0$,
denote $N_r(A) = \{x \in X \suchthat \exists a \in A \quad\text{such that}\quad d(a,x) \le r\}$.
Given subsets $A, B \subset X$, let $A \csubset{[r]} B$ denote $A \subset N_r(B)$. Let
$A\csubset{c} B$ denote the existence of $r \ge 0$ such that $A \csubset{[r]} B$; this is called
\emph{coarse containment} of $A$ in $B$. Let $A \ceq{[r]} B$ denote the conjunction of $A
\csubset{[r]} B$ and $B\csubset{[r]} A$, equivalently $d_\Haus(A,B) \le r$ where
$d_\Haus(\cdot,\cdot)$ denotes Hausdorff distance. Let $A\ceq{c} B$ denote the existence of $r$
such that $A\ceq{[r]} B$; this is called \emph{coarse equivalence} of $A$ and $B$. If $A
\csubset{c} B$ but $A \not\ceq{c} B$ then we say that $A$ is \emph{strictly coarsely contained}
in $B$, denoted $A \cstrict{c} B$. 

Note the two versions of coarse set theory notation, for example $\ceq{c}$ versus $\ceq{[r]}$. The
unbracketed symbol ``$c$'' simply stands for ``coarse''; whereas in a bracketed symbol ``$[r]$''
the $r$ stands for a distance.

Given metric spaces $X,Y$ and functions $f,g \from X \to Y$ let $f \ceq{[r]} g$ denote $f(x) \ceq{[r]}
g(x)$ for all $x \in X$, equivalently the sup norm distance between $f$ and $g$ is at most $r$.

Given a metric space $X$ and subsets $A,B$, let $A \cintersect{[r]} B = N_r(A) \intersect
N_r(B)$. We say that a subset $C$ is a \emph{coarse intersection} of $A$ and $B$, written~$C =
A \cintersect{c} B$, if $C \ceq{c} A \cintersect{[r]} B$ for all sufficiently large $r$. A
coarse intersection of $A$ and $B$ may not exist, but if one does exist then it is well-defined
up to coarse equivalence. 

For an example in the real line where coarse intersection does not exist, let $A = \{2^n
\suchthat n=1,2,3,\ldots\}$. Choose a decomposition $A=A_1 \union A_2 \union \cdots$ where each
$A_i$ is infinite. Let $B_i = \{a+2i \suchthat a \in A_i\}$, and let $B=B_1 \union B_2
\union \cdots$. Then $A \cintersect{[i]} B \ceq{c} A_1 \union\cdots\union A_i$, but
$A_1\union\cdots\union A_i \not\ceq{c} A_1\union\cdots\union A_j$ for $i < j$, and so
$A\cintersect{c} B$ does not exist.

More generally, we say that $C$ is a coarse intersection of a finite collection of subsets
$A_1,\ldots,A_n$ if for all sufficiently large $r$ we have 
$$\bigl( (A_1 \cintersect{[r]} A_2) \cintersect{[r]}
\cdots \bigr)\cintersect{[r]} A_n \ceq{c} C
$$ 
Although the binary operator $\cintersect{[r]}$ is not associative, it is ``coarsely
associative'', meaning for each $r \ge 0$ there exists $r' \ge 0$ such that $(A \cintersect{[r]}
B) \cintersect{[r]} C \csubset{[r']} A \cintersect{[r']} (B \cintersect{[r']} C)$, and similarly
for the other direction. This implies that the coarse intersection operator $\cintersect{c}$
is associative as long as it is defined, and so we can write sentences like $A_1\cintersect{c}
A_2\cintersect{c}\cdots \cintersect{c} A_n \ceq{c} C$.  For example,
Lemma~\ref{LemmaCoarseIntersectionSubgroup} will show that if $X$ is a finitely generated group with
the word metric and $A_1,\ldots,A_n$ are subgroups of $X$ then $A_1\cintersect{c}\cdots\cintersect{c}
A_n \ceq{c} A_1 \intersect
\cdots \intersect A_n$.

\paragraph{Deepness.} Consider a metric space $X$.  Given a subset $A \subset X$ and $r>0$, we say that a point $x \in
A$ is \emph{$r$ deep in $A$} if $N_r(x) \subset A$, and we say that $A$ is \emph{$r$ deep} if it contains an $r$ deep
point. We also say that $A$ is \emph{deep} if it is $r$ deep for each $r \ge 0$. A subset $B \subset X$ is said to have
\emph{deep intersection} with $A$ if $B$ contains an $r$ deep point of $A$ for each $r \ge 0$.

Given two subsets $A,B \subset X$, we say that $A$ \emph{crosses $B$ in $X$} if for some $r \ge 0$ the set $A$ has deep
intersection with at least two components of $N_r(B)$; equivalently, there exists $R \ge 0$ such that for all $r \ge
R$, $A$ has deep intersection with at least two components of $N_r(B)$. Note that this property is quasi-isometrically
invariant: if $f\from X\to Y$ is a quasi-isometry of metric spaces and $A,B \subset X$, then $A$ crosses $B$ in $X$ if
and only if $f(A)$ crosses $f(B)$ in $Y$.

\paragraph{Quasi-isometries.} Consider metric spaces $X,Y$ and a map $f \from X \to Y$. The
map $f$ is \emph{$K,C$ coarse Lipschitz}, with $K \ge 1$ and $C \ge 0$, if $d_Y(fx,fy) \le K
d_X(x,y) + C$ for all $x,y \in X$. Also, $f$ is a \emph{uniformly proper embedding} if it is
coarse Lipschitz and there exists a proper, increasing function $\rho \from [0,\infinity) \to
[0,\infinity)$, called a \emph{gauge function}, such that
$$\rho(d_X(x,y)) \le d_Y(fx,fy)
$$
If $f$ is $K,C$ coarse Lipschitz and uniformly proper with gauge function $\rho(d) = \frac{1}{K} d
- C$ then we say that $f$ is a $K,C$ \emph{quasi-isometric embedding}.

For example, if $G$ is a finitely generated group and $H$ is a finitely generated subgroup,
then the inclusion map $H \subset G$ is a uniformly proper embedding from $H$ with its word
metric to $G$ with its word metric, but it need not be a quasi-isometric embedding.

A map $f \from X \to Y$ is \emph{$C$ coarsely surjective}, with $C \ge 0$, if $f(X) \ceq{[C]} Y$. 
Two maps $f \from X \to Y$, $\bar f \from Y \to X$ are \emph{$C$ coarse inverses} if $\bar
f \composed f \ceq{[C]} \Id_X$ and $f \composed \bar f \ceq{[C]} \Id_Y$.

A $K,C$ quasi-isometric embedding $f \from X \to Y$ which is $C$ coarsely surjective is called
a \emph{$K,C$ quasi-isometry}. A uniformly proper map $f \from X \to Y$ which is coarsely surjective
is called a \emph{uniform equivalence} or a \emph{uniformly proper equivalence}.

For every quasi-isometry $f \from X \to Y$ there exists a quasi-isometry $\bar f \from Y \to X$ such
that $f, \bar f$ are coarse inverses, and the quasi-isometry constants for $\bar f$ and the coarse
inverse constant depend only on the quasi-isometry constants of $f$. For every uniform equivalence $f
\from X \to Y$ there exists a uniform equivalence $\bar f \from Y \to X$ such that $f, \bar f$ are
coarse inverses; the gauge function and coarse Lipschitz constants for $\bar f$, and the coarse
inverse constant, depend only on gauge function and coarse Lipschitz constants for $f$.

Let $X$ be a metric space. A \emph{geodesic} in $X$ connecting $x$ to $y$ is a rectifiable
path from $x$ to $y$ whose length equals $d(x,y)$. A \emph{$K,C$-quasigeodesic} in $X$
connecting $x$ to $y$ is a $K,C$-quasigeodesic embedding $p \from [a,b] \to X$ such that
$p(a)=x$, $p(b)=y$. We say that $X$ is a \emph{geodesic metric space} if any two points of $X$
are connected by geodesic. We say that $X$ is a \emph{proper metric space} if closed metric
balls in $X$ are compact.

\begin{proposition}[Coarse convergence principle] 
\label{PropCC}
For each $K,C$ there
exists $A$ such that if $f_i \from X \to Y$ is a sequence of $K,C$
quasi-isometric embeddings, and if there exists a point $x \in X$ such
that $f_i(x)$ is a bounded subset of $Y$, then there exists a $K,A$
quasi-isometry $f_\infinity \from X \to Y$ and a subsequence
$f_{i_n}$ satisfying the following property:
\begin{itemize}
\item For each $x$ there exists $N$ such that if $n \ge N$ then
$$d(f_{i_n}(x),f_\infinity(x)) \le A$$
\end{itemize}
We say that $f_\infinity$ is an \emph{$A$-coarse limit} of $f_{i_n}$.
\qed\end{proposition}

\paragraph{Quasi-actions.}
Let $G$ be a group and $X$ a metric space. A \emph{$K,C$ quasi-action} of $G$ on $X$ is a map
$G \cross X \to X$, denoted $(g,x) \mapsto A_g(x)=g\cdot x$, with the following properties:
\begin{itemize}
\item For each $g \in G$ the map
$A_g \from X \to X$ is a $K,C$ quasi-isometry of $X$.
\item For each $g,h
\in G$ we have
$$g \cdot (h \cdot x) \ceq{[C]} (gh) \cdot x, \forall x \in X,
\quad\text{equivalently,}\quad A_g \composed A_h \ceq{[C]} A_{gh}
$$
\end{itemize}
A quasi-action is \emph{cobounded} if there exists a constant $r$ such that for each $x \in X$
we have $G \cdot x\ceq{[r]} X$. A quasi-action is \emph{proper} if for each $r$ there exists
$M$ such that for all $x,y\in X$, the cardinality of the set $\{g\in G \suchthat \bigl(g \cdot
N(x,r)\bigr) \intersect N(y,r) \ne\emptyset\}$ is at most $M$. Note that if $G \cross X \to X$
is an isometric action on a proper metric space, then ``cobounded'' is equivalent to
``cocompact'' and ``proper'' is equivalent to ``properly discontinuous''.

Given a group $G$ and quasi-actions of $G$ on metric spaces $X,Y$, a \emph{quasiconjugacy} is a
quasi-isometry $f \from X \to Y$ such that for some $r \ge 0$ we have $f(g \cdot x) \ceq{[r]}
g \cdot fx$ for all $g \in G$, $x  \in X$. Properness and coboundedness are invariants of
quasiconjugacy.

A fundamental principle of geometric group theory says that if $G$ is a finitely generated
group equipped with the word metric, and if $X$ is a proper geodesic metric space on which $G$
acts properly discontinuously and cocompactly by isometries, then $G$ is quasi-isometric to
$X$, in fact for any base point $x_0 \in X$ the map $g \mapsto g \cdot x_0$ is a
quasi-isometry $G \mapsto X$.

More generally, if $G$ quasi-acts properly and coboundedly on two proper, geodesic metric
spaces $X,X'$, then there is an induced quasi-isometry $f \from X\to X'$ which is $G$-almost
equivariant, and $f$ is uniquely determined up to finite distance. To define such an $f$,
choose quasi-isometries $h \from G \to X$, $h' \from G \to X'$ as above, choose $\bar h \from
X \to G$ to be a coarse inverse for $h$, and let $f = h' \composed \bar h$.

A converse to this result is the \emph{quasi-action principle} which says that if $G$
is a finitely generated group with the word metric and $X$ is a metric space quasi-isometric
to $G$ then there is a cobounded, proper quasi-action of $G$ on $X$: take the left action of
$G$ on itself and quasiconjugate by a quasi-isometry $G \to X$ and its coarse inverse. The
constants for this quasi-action depend only on the quasi-isometry constants for $G \to X$.

\subsection{Coarse properties of subgroups}
\label{SectionCoarseSubgroups}
In a finitely generated group $\Gamma$ with a symmetric generating set, we
use $\norm{\gamma}$ to denote the word norm of $\gamma$, which is the
length of the shortest word in the generating set which evaluates to be
$\gamma$. The (left invariant) word metric is then $d(g,h) = \norm{g^\inv h}$.

\begin{lemma}[Coarse intersection of subgroups is well defined]
\label{LemmaCoarseIntersectionSubgroup}
Let $\Gamma$ be a finitely generated group. Let $G_1,\ldots,G_n$ be a 
finite collection of subgroups. Then the coarse intersection of 
$G_1,\ldots,G_n$ is well-defined, and is represented by $G_1 
\intersect \cdots \intersect G_n$.
\end{lemma}

Note that the subgroups need not be finitely generated; the lemma holds
for \emph{any} subgroups $G_1,\ldots,G_n$.

\begin{proof}
We'll give the proof just for $n=2$. Let the subgroups be $G,H$. Let the
generating set be $A$.

Pick an $r$ and consider the set $N_r(G) \intersect N_r(H)$ in the Cayley
graph of $\Gamma$; we must show that this set is coarsely contained in $G
\intersect H$. Pick a vertex $x$ of the Cayley graph lying in $N_r(G)
\intersect N_r(H)$. It follows that there is an edge path $a_1 *
\cdots * a_k$ in the Cayley graph from an $H$ vertex $h_x$ to a $G$ vertex
$g_x$, with $k \le 2r$, such that this edge path passes over the vertex
$x$.  Let $\gamma_x =  a_1 \cdot\ldots\cdot a_k = h_x^\inv g_x$, so
$\norm{\gamma_x} \le 2r$.

Consider the set of all elements $\gamma \in \Gamma$ such that
$\norm{\gamma}\le 2r$ and such that the equation $\gamma = h^\inv g$ has
a solution $(g,h)\in G \cross H$; choose a solution $(g_\gamma,h_\gamma)$
so that $h_\gamma$ has minimal norm $N_\gamma$. Let $N$ be the maximum of
$N_\gamma$. 

With $x$ as before, setting $\gamma=\gamma_x$, we have
$$\gamma = h_\gamma^\inv g_\gamma = h_x^\inv g_x
$$
so
$$h_x h_\gamma^\inv = g_x g_\gamma^\inv
$$
is an element of $G \intersect H$, and since
$$h_x = (g_x g_\gamma^\inv) h_\gamma
$$
it follows that $h_x$ is within distance $N$ of an element of $G
\intersect  H$. The point $x$ is therefore within distance $N+r$ of an
element of $G
\intersect H$.
\end{proof}

\begin{lemma}[Coarse equivalence of subgroups] 
\label{LemmaSubgroupCoarseEquivalence}
Let $\Gamma$ be a finitely generated group with subgroups
$B \subgroup A \subgroup \Gamma$, and suppose that $B \ceq{c} A$. Then
$B$ has finite index in $A$.
\end{lemma}

\begin{proof} There exist $\gamma_1,\ldots,\gamma_k \in \Gamma$ such that
$$B \subset A \subset B \gamma_1 \union \cdots \union B \gamma_k
$$
We may assume that $A \intersect B \gamma_i \ne\emptyset$ for each
$i=1,\ldots,k$, for otherwise we can just drop the term $B \gamma_i$.
It follows that for each $\gamma_i$ we have $a = b \gamma_i$ for some $a
\in A, b
\in  B$, and so $\gamma_i = b^\inv a \in A$. Thus, 
$$B \subset A \subset B \gamma_1 \union \cdots \union B \gamma_k 
\subset A
$$
and so 
$$A = B \gamma_1 \union \cdots \union B \gamma_k 
$$
\end{proof}

Recall that two subgroups $A,B$ of a group $\Gamma$ are \emph{commensurable} if $A \intersect
B$ has finite index in both $A$ and $B$. This is an equivalence relation on subgroups of
$\Gamma$. Define a partial order on subgroups of $\Gamma$ where $A \prec B$ if $A \intersect
B$ has finite index in $A$, equivalently, $A$ is commensurable to a subgroup of $B$. The
equivalence relation generated by this partial order is the same as commensurability. The
partially ordered set obtained from the relation~$\prec$ by passing to commensurability classes
is called the \emph{commensurability lattice} of $G$.

\begin{corollary}[Commensurability lattice $=$ coarse inclusion lattice]
\label{CorollaryCommensurable}
If $\Gamma$ is a finitely generated group with subgroups $A,B$, then the following hold:
\begin{enumerate}
\item $A \csubset{c} B$ if and only if $A \prec B$.
\item $A \ceq{c} B$ if and only if $A,B$ are commensurable in $\Gamma$.
\end{enumerate}
\end{corollary}

\begin{proof} Clearly (2) follows from (1). The ``if'' direction of (1) is obvious.

To prove the ``only if'' direction of (1), assume that $A \csubset{c} B$, and so $A
\csubset{c} A \cintersect{c} B$. Applying Lemma~\ref{LemmaCoarseIntersectionSubgroup} we
obtain $A \csubset{c} A \intersect B \subset A$, and so $A \ceq{c} A \intersect B$.
Lemma~\ref{LemmaSubgroupCoarseEquivalence} now implies that $A \prec B$.
\end{proof}

\subsection{Coboundedness principle} 
\label{SectionCobounded}


Given a proper, cobounded quasi-action of a group on a proper, geodesic metric space, we often
want to study the ``coarse stabilizer subgroups'' of certain subsets of the metric space. In
order for this to make sense, and in order to obtain stabilizers with good properties, we must
be careful about the properties of the subsets. We follow a method of Kapovich and Leeb
\cite{KapovichLeeb:haken}.

Let $X$ be a proper, geodesic metric space. A \emph{pattern} in $X$ is a collection $\A$ of
subsets $A \subset X$, each of which is a proper, geodesic metric space in its own right,
mapped into $X$ by a $K,C$ coarse lipschitz, $\rho$-uniformly proper map, with coarse
lipschitz constants $K,C$ and uniform properness gauge $\rho$ independent of $A \in \A$.

Given a pattern $\A$ and a cobounded, proper quasi-action of a group $G$ on
$X$, consider the following properties:
\begin{description}
\item[(a) Local Finiteness Condition:] Any set of finite diameter in $X$ intersects
at most finitely many elements of $\A$.
\item[(b) Coarse Discreteness Condition:] Any two distinct elements of $\A$ have infinite
Hausdorff distance from each other.
\item[(c) Coarse Action Condition:] There exists $K \ge 1$, $C \ge 0$ such that for any $g \in
G$, $A \in \A$ there exists $A' \in \A$ such that:
\begin{itemize}
\item $g \cdot A \ceq{[C]} A'$, 
\item the map $A \to g \cdot A \to A'$, composing $x \mapsto g \cdot x$ with a
closest point map $g \cdot A \to A'$, is a $K,C$ quasi-isometry.
\end{itemize}
\end{description}
The coarse discreteness condition (b) combined with the coarse action condition (c) imply
that, for each $g \in G$ and $A \in \A$, the element $A'$ in (c) is unique. Conditions (b) and
(c) therefore imply:
\begin{description}
\item[(d) Action Condition:] There exists an action of $G$ on the set $\A$, denoted
$g(A) \in \A$ for each $g \in G$ and $A \in \A$, and there exist $K \ge 1$, $C \ge 0$, such
that for each $g \in G$, $A \in \A$ we have:
\begin{itemize}
\item $g \cdot A \ceq{[C]} g(A)$ for $g \in G$, $A \in \A$
\item The map $A \to g \cdot A \to g(A)$, composing $x \mapsto g \cdot x$ with a
closest point map $g \cdot A \to g(A)$, is a $K,C$ quasi-isometry.
\end{itemize}
\end{description}
The conjunction of conditions (b) and (c) is strictly stronger than condition (d); consider,
for instance, the pattern of vertical lattice lines in $\reals^2$ and the action of $\Z^2$.
When conditions (b) and (c) hold, the action in (d) is uniquely characterized by the condition
that $g(A) \ceq{c} g \cdot A$, for all $g \in G$, $A \in \A$.

Assuming the action condition (d) holds, for each $A \in \A$ we have a subgroup $\Stab(A) =
\{g\in G \suchthat g(A) = A \}$ of the group $G$.

There is another interesting subgroup which might be called the ``coarse stabilizer'' of $A$,
defined by $\CStab(A) = \{g\in G \suchthat g \cdot A \ceq{c} A\}$. When the coarse
discreteness condition (b) and the coarse action condition (c) hold then $\Stab(A) =
\CStab(A)$. If these conditions fail, however, then $\CStab(A)$ may be strictly larger than
$\Stab(A)$. The subgroup $\Stab(A)$ is generally more useful. For example, if $\G$ is a
homogeneous graph of groups with Bass-Serre tree of spaces $X \to T$, e.g.\ the standard
circle of $\Z$'s presentation for the Baumslag-Solitar group $\BS(p,q) = \{a,t
\suchthat t a^p t^\inv = a^q\}$, then for each vertex or edge space $A$ of $X$ the group
$\CStab(A)$ is just the whole group $\pi_1\G$, which gives no information.

For each $g \in \Stab(A)$ there is an induced quasi-isometry $\alpha_g \from A\to A$, defined
as the composition of the map $x \mapsto g \cdot x$ with a closest point map from $g\cdot A$
to~$A$. Note that the quasi-action constants of the maps $\alpha_g$ are uniformly bounded for
$g \in \Stab(A)$, and in fact the map $g \mapsto \alpha_g$ defines a quasi-action of
$\Stab(A)$ on the set $A$. This quasi-action is automatically proper, because the ambient
quasi-action of $\Stab(A)$ on $X$ is proper, the embedding $A \to X$ is uniformly proper, and
the closest point maps $g \cdot A \to A$ move points a uniformly bounded amount. Although each
element of the larger subgroup $\CStab(A)$ also induces a quasi-isometry of $A$, the
quasi-isometry constants may not be uniform, and so $\CStab(A)$ may not quasi-act on $A$.

Coboundedness of the quasi-action of $\Stab(A)$ on $A$ is provided by the following:

\begin{proposition}[The Coboundedness Principle]
\label{PropCoboundednessPrinciple}
Given a metric space $X$, a pattern $\A$ in $X$ satisfying the local finiteness condition (a),
and a cobounded, proper quasi-action of a group $G$ on $X$ satisfying the action
condition (d), the following properties hold:
\begin{enumerate}
\item There are only finitely many orbits of the action of $G$ on $\A$.
\item For each $A \in \A$ the quasi-action of $\Stab(A)$ on $A$ given by $g \mapsto
\alpha_g$ is a cobounded. 
\end{enumerate}
\end{proposition} 

Condition 2 of this proposition is essentially the same as \cite{KapovichLeeb:haken}, Lemma
5.2, whose proof requires the additional assumption that there exists a ball $B \subset X$
such that $G \cdot B = X$.

\begin{proof} To prove (1), since the quasi-action of $G$ on $X$ is cobounded, we can find
within every $G$-orbit of $\A$ a representative which intersects some fixed ball in $X$. 
Since this ball is finite, there can only be finitely many such representatives and hence
only finitely many orbits.

To prove (2), since the embedding $A \to X$ is uniformly proper, it suffices to check
coboundedness with respect to the restriction of the ambient metric on $X$ to the subset $A$.
Fix $K,C$ so that each map $x\to g\cdot x$ is a $K,C$ quasi-isometry of $X$, so
that $d(g\cdot (g' \cdot x), gg'\cdot x) \le C$, and so that the $G$ quasi-action is
$C$-cobounded on $X$. 

Fix a base point $x_0 \in A$. Fix a constant $D$, to be determined later, and let $\A_0$ be
the finite set consisting of all elements of $\A$ contained in the $G$-orbit of $A$
which intersect the $D$ ball around $x_0$. For each $B \in \A_0$ choose $g_B\in G$ so that
$g_B(A) = B$. Let $r$ be the maximum of $d(x_0,g_B \cdot x_0)$, for $B \in\A_0$.

Consider a point $x \in A$. We must find an element of $\Stab(A)$ moving
$x_0$ to within a fixed distance of $x$. There exists $g \in G$ such that
$d(g \cdot x_0,x) \le C$. Note that 
\begin{align*}
d(x_0,g^\inv \cdot x) &\le d(x_0, g^\inv \cdot (g \cdot x_0)) +
d(g^\inv \cdot (g \cdot x_0), g^\inv x) \\
  &\le C + (KC+C) = KC+2C
\end{align*}
Since $x \in A$ it follows that $g^\inv \cdot x$ is within distance $C$
of a point of $g^\inv(A)$. Now we set $D=KC+3C$, and so $g^\inv(A) \in
\A_0$. Setting $B = g^\inv \cdot A$ we thus have $g^\inv(A) = g_B(A)$ and
so $gg_B \in
\Stab(A)$. We also have $d(x_0,g_B \cdot x_0) \le r$. It follows that
\begin{align*}
d(g^\inv \cdot x, g_B \cdot x_0) &\le d(g^\inv \cdot x, x_0) + d(x_0, g_B
\cdot x_0) \\
  & \le KC+2C+r
\end{align*}
and so
\begin{align*}
d(x,(gg_B) \cdot x_0) &\le d(x, g \cdot (g^\inv \cdot x)) +
   d(g \cdot (g^\inv \cdot x), g \cdot (g_B \cdot x_0)) \\
& \qquad\qquad\qquad\qquad\qquad\qquad + d(g \cdot (g_B \cdot x_0), (g
g_B) \cdot x_0)
\\
  &\le C + [K(KC+2C+r)+C] + C \\
  &= K^2C + 2KC + Kr + 3C
\end{align*}
\end{proof}

\subsection{Bass-Serre trees and Bass-Serre complexes}
\label{SectionBassSerre}

Associated to a graph of groups is its \emph{fundamental group}, an algebraic definition of
which is given in \cite{Serre:trees}. We shall content ourselves with the topological
definition, following \cite{ScottWall}, which leads quickly to the concepts of Bass-Serre
trees and complexes. Then we introduce a related idea, a coarse Bass-Serre complex.

\paragraph{Graphs of groups.} Let $\G$ be a graph, that is, a 1-dimensional CW-complex. The two ends of an edge $e$ of
$\G$ form a set $\Ends(e)$, and each $\eta \in\Ends(e)$ is attached to a vertex $v(\eta)$. 

A \emph{graph of groups} is a graph $\G$ equipped with a \emph{vertex group} $\G_v$ for each
vertex $v$, an \emph{edge group} $\G_e$ for each edge $e$, and an injective
\emph{edge-to-vertex homomorphism} $\phi_\eta \from\G_e \to \G_{v(\eta)}$ for each edge
$e$ and each $\eta \in \Ends(e)$. If $\G$ is finite, its vertex groups are finitely presented,
and its edge groups are finitely generated, then we say that $\G$ has \emph{finite type}.

A \emph{graph of spaces} $\K$ over $\G$ is constructed as follows. For each vertex $v$ choose
a pointed, connected CW-complex $\K_v$, the \emph{vertex space}, and an
isomorphism $\pi_1 \K_v\to \G_v$; if $\G$ is of finite type we can and shall choose $\K_v$ to
be compact. For each edge $e$ choose a pointed, connected CW-complex
$\K_e$ for each edge $e$, the \emph{edge space}, and an epimorphism $\pi_1 \K_e \to \G_e$
(this map is \emph{not} required to be an isomorphism); if $\G$ is of finite type we can
and shall choose $\K_e$ to be compact. For each edge $e$ and each $\eta \in
\Ends(e)$ choose a pointed CW-map $f_\eta \from \K_e \to \K_{v(\eta)}$, the
\emph{edge-to-vertex map}, such that the
following diagram commutes:
$$
\xymatrix{
\pi_1(\K_e) \ar[r]^{(f_\eta)_*} \ar[d] & \pi_1(\K_{v(\eta)}) \ar[d]\\
\G_e \ar[r]^{\phi_\eta} & \G_{v(\eta)}
}
$$
For each edge $e$ let $\bar e$ denote the end compactification of $\interior(e)$, homeomorphic to a
closed interval, with an endpoint also denoted $\eta$ for each $\eta\in\Ends(e)$. The identity map on
$\interior(e)$ extends to a continuous map $\bar e \to e$ (this notation gets around the possibility
if an edge with both ends attached to the same vertex). Now construct a CW-complex from the disjoint
union
$$\coprod_{v \in V(\G)} \K_v \quad\disjunion\quad \coprod_{e \in
E(\G)} \K_e \cross \bar e
$$
by making the following identifications: for each endpoint $\eta$ of $\bar e$, glue $\K_e
\cross \eta$ to $\K_v$ by identifying $(x,\eta)$ with $f_\eta(x)$ for each $x \in \K_e$. We
obtain a CW-complex $X$ and a projection map $p \from \K \to \G$ induced by the collection of
constant maps $\K_v \to v$ and projection maps $\K_e \cross \bar e\to\bar e \to e$. Let $\K_t
= p^\inv(t)$ for each point $t \in \G$, so $\K_t = \K_v$ if $t=v$ is a vertex, and $\K_t$ is a
copy of $\K_e$ if $t \in \interior(e)$ for an edge $e$. If $\G$ is of finite type then $\K$ is
finite.

An application of Van Kampen's theorem shows that $\pi_1 \K$ is well-defined up to isomorphism
regardless of the choices in the construction of $\K$, and we define the fundamental group of
$\G$ to be $\pi_1\G \approx \pi_1\K$. Van Kampen's theorem also connects this topological
definition of $\pi_1\G$ with the algebraic definition in \cite{Serre:trees}. Note that if $\G$
has finite type then $\pi_1\G$ is finitely presented. 

Given a finite type graph of groups $\G$ and a graph of spaces $\K \to \G$, we construct the
\emph{Bass-Serre tree of spaces} $X \to T$ as follows. Let $\pi \from X \to \K$ be a universal
covering map, with $\pi_1\G$ acting on $X$ by deck transformations, properly and cocompactly.
The decomposition of $\K$ into subsets $\{\K_t\}$ induces a decomposition of $X$ whose
elements are path connected lifts of decomposition elements $\K_t$; this decomposition of $X$
is preserved by the action of $\pi_1\G$. The quotient space of this decomposition of $X$ is
denoted $T$. It is a tree on which $\pi_1\G$ acts without edge inversions, and the quotient
map $X \to T$ is a $\pi_1\G$ equivariant map called a \emph{Bass-Serre tree of spaces} for
$\G$. Indeed, the $\pi_1\G$-tree $T$ is equivariantly isomorphic to the Bass-Serre tree of
$\G$ as defined in \cite{Serre:trees}. For each $t\in T$ the set $p^\inv(t) \subset X$ is
denoted $X_t$. If $v$ is a vertex then $X_v$ is called a \emph{vertex space}. If $t$ is the
midpoint of an edge $e$ of $T$ then we set $X_e = X_t$ and call this an \emph{edge space}; also,
the two components of $X-X_e$ are labelled $U^\pm_e$.

For each $t \in T$ its stabilizer $\Stab(t)$ in $\pi_1\G$ is also the stabilizer of $X_t$, and
the projection from $X_t$ to $\K_{\pi(t)}$ is a covering map with deck transformation group
$\Stab(t)$; if moreover $t$ is a vertex then $X_t \mapsto \K_{\pi(t)}$ is a universal covering
map. If $t$ is the midpoint of an edge $e$ then $X_e = X_t \mapsto \K_{\pi(t)} =
\K_{\pi(e)}$ need not be a universal covering map, but it is as long as we set things up so that
the epimorphism $\pi_1\K_{\pi(e)} \to \G_{\pi(e)}$ is an isomorphism.

\subparagraph{Remark on terminology:} Henceforth the terminology of ``vertex spaces'' and ``edge spaces'' always
applies to subsets of the Bass-Serre complex $X$, not to its finite quotient complex $\K$.

\paragraph{Depth zero rafts.} We define \emph{depth zero rafts} of $T$ as follows. A vertex or edge $a$ of
$T$ is said to have \emph{depth zero} if $X_a$ is maximal with respect to coarse inclusion of
vertex or edge spaces of $X$, that is, for any other vertex or edge $b$ of $T$, if $X_a
\csubset{c} X_b$ then $X_a \ceq{c} X_b$. Two depth zero vertices or edges $a,b$ are \emph{in
the same raft} if $X_a \ceq{c} X_b$; this defines an equivalence relation on depth zero
vertices and edges, and the equivalence classes are called \emph{depth zero rafts}. Each depth
zero raft is a closed subcomplex of $T$, for if $e$ is a depth zero edge and $v \in \bdy e$
then $X_e \csubset{c} X_v$ and so $X_e \ceq{c} X_v$ implying that $X_v$ also has depth zero.
In fact, each depth zero raft is a subtree of $T$, for suppose that $v,w$ are depth zero
vertices, consider the arc $\overline{vw}$ in $T$ with endpoints $v,w$, and let $e$ be an edge
in $\overline{vw}$. We clearly have $X_v \cintersect{c} X_w \csubset{c} X_e$, but since $X_v
\ceq{c} X_w$ it follows that $X_v \cintersect{c} X_w \ceq{c} X_v \ceq{c} X_w$, and hence $X_v
\csubset{c} X_e$. Since $v$ is of depth zero it follows that $X_v \ceq{c} X_e$ and so $e$ is
also of depth zero and is contained in the same depth zero raft as~$v$.

Given a depth zero raft $R \subset T$, the group $\Stab(R)$ acts cocompactly on $R$. The quotient graph of
groups $R / \Stab(R)$ has fundamental group identified with $\Stab(R)$. The inverse image of $R$ under the
projection $X \to T$ is called the \emph{raft space} corresponding to $R$, denoted $X_R$. Note that 
$$X_R \ceq{c} \bigcup_{v \in R} X_v \ceq{c} \left( \bigcup_{v \in R} X_v \right) \union \left( \bigcup_{e
\in R} X_e \right)
$$
The raft space $X_R$ is identified with the Bass-Serre complex for the quotient graph of groups $R /
\Stab(R)$.

See Section~\ref{SectionFiniteDepth} for a deeper study of rafts, in particular rafts deeper than
depth zero.

We have seen in Proposition~\ref{CorollaryCommensurable} that two subgroups $A,B$ of a
finitely generated group $\Gamma$ are commensurable in $\Gamma$ if and only if $A,B$ are
coarsely equivalent with respect to the word metric on $\Gamma$. From this we obtain an
equivalent definition of a depth zero raft in $T$, namely, a subtree $R \subset T$ such that if
$e$ is an edge of $R$ incident to a vertex $v$ of $R$ then $\Stab(e)$ has finite index in
$\Stab(v)$, whereas if $e$ is an edge not in $R$ that is incident to a vertex $v$ of
$R$ then $\Stab(e)$ has infinite index in $\Stab(v)$. This is precisely the definition of a
depth zero raft given in the introduction.

We can also identify depth zero rafts in terms of the original graph of groups $\G$. Define a depth zero
raft in $\G$ as a maximal subgraph $\G' \subset \G$ with the property that if an edge $e$ of $\G'$ is
incident to a vertex $v$ of $\G'$ then the injection $\G_e \inject \G_v$ has finite index image, whereas if
an edge $e$ not in $\G'$ is incident to a vertex $v$ of $\G'$ then the injection $\G_e \inject \G_v$ has
infinite index image. The depth zero rafts in the Bass-Serre tree of $\G$ are just the connected components
of the pre-images of the depth zero rafts in $\G$ itself.

\subsection{Irreducible graphs of groups}
\label{SectionIrreducible}

Consider a finite type graph of groups $\G$.

We say that $\G$ is \emph{reducible} along an edge $e$ if $e$ has ends at
distinct vertices $v,w$ and the injection $\phi_{ev} \from \G_e \to \G_v$ is surjective. In
this case we obtain a new graph of groups $\G'$ by \emph{collapsing $\G$ along $e$},
identifying $v,w,e$ to a single vertex $u$ with $\G'_u = \G_w$; edge to vertex injections at
vertices other than $u$ are unchanged; at $u$, each incident edge $e'$ at $w$ other than $e$
defines an injection $\G'_{e'} \to \G'_u$ agreeing with $\G_{e'} \to \G_w$; and each edge to
vertex injection $\phi_{e'v} \from\G_{e'} \to \G_v$ is replaced by the injection $\G'_{e'} \to
\G'_u$ given by $\phi_{ew}\composed\phi_{ev}^\inv \composed \phi_{e'v}$. Reduction of $\G$
along an edge decreases the number of edges and so must come to a halt after finitely many
iterations, a process called \emph{complete reduction}. 

The end product of a complete reduction is a graph of groups that is \emph{irreducible} meaning
that it is not reducible along any edge.

\subparagraph{Example.} The complete reduction process is not canonical: a graph of groups $\G$
can have nonisomorphic complete reductions $\G_1,\G_2$, indeed it may happen that no vertex group
of $\G_1$ is isomorphic to a vertex group of $\G_2$. Here is an example. Consider the groups
$$H_n = \<x,y,z \suchthat [x,z]=[y,z]=1, [x,y]=z^n\>
$$
The group $H_1$ is just the integer Heisenberg group. For each integer $p>0$ there is an index
$p^2$ self-embedding $f \from H_1 \to H_1$ defined by $x \to x^p$, $y \to y^p$, $z \to
z^{p^2}$. In between $f(H_1)$ and $H_1$ is an index $p$ subgroup of $H_1$ and an index $p$
supergroup of $f(H_1)$ isomorphic to $H_p$, namely the image of the embedding $H_p \to
H_1$ given by $x \to x$, $y \to y^p$, $z \to z^p$. Construct a circle of groups $\G$ with two
edges $e,f$ and two vertices $v,w$, let $\G_e=\G_v=H_1$ with $\G_e\inject\G_v$ the identity,
let $\G_f=\G_w=H_p$ with $\G_f\inject\G_w$ the identity, let $\G_f\to\G_v$ be the given
index~$p$ inclusion $H_p \to H_1$, and let $\G_e\to\G_w$ be the index~$p$ injection $H_1
\xrightarrow{f} f(H_1) \inject H_p$. Collapsing $f$ gives an ascending HNN extension of $H_1$,
wherease collapsing $e$ gives an ascending HNN extension of $H_p$, and clearly $H_1
\not\approx H_p$.

\begin{example}\label{universal}
Here is an even worse example: by allowing edge-to-vertex injections of infinite index, we can
arrange that there are two different complete reductions whose vertex and edge groups are not
even of the same quasi-isometry type. Higman constructed a finitely presented group
$H$ in which any finitely presented group can be embedded \cite{Higman:subgroups}.  Let $K=H
\cross \Z$ and $L=H * \Z$.  Since both are finitely presented, they embed as subgroups of $H$
and hence of each other.  On the other hand, $K$ is one ended while $L$ is not, so they are not
quasi-isometric.

Choose infinite index embeddings $\alpha \from L \to K$, $\beta \from K \to L$. Construct a
circle of groups $\G$ with two edges $e,f$ and two vertices $v,w$, let $\G_e = \G_v = L$ and let
$\G_e\inject \G_v$ be the identity, let $\G_f = \G_w = K$ and let $\G_f \inject \G_w$ be the
identity, let $\G_e \to\G_w$ be $\alpha$ and let $\G_f \to \G_v$ be $\beta$. Collapsing $e$
results in a strictly ascending HNN extension of $K$ by the map $\beta\alpha$, whereas
collapsing $f$ results in a strictly ascending HNN extension of $L$ by the map $\alpha\beta$.

Note that there is an index $k$ subgroup of $\pi_1\G$ which has a circular graph of groups
decomposition that is a degree $k$ covering of $\G$. By collapsing all but one of the $2k$ edges
in the covering graph, we can express this subgroup as a strictly ascending HNN extension of $K$
by the map $(\beta\alpha)^k$, and also as a strictly ascending HNN extension of $L$ by the map
$(\alpha\beta)^k$.
\end{example}

Notice one crucial difference between these two examples: in the first example, which is of
finite depth, the vertex groups of the reduction \emph{are} well defined up to quasi-isometry.
Our Classification Theorem~\ref{TheoremClasses} shows that, in certain contexts, the end
product of reduction \emph{is} canonical in a geometric sense, in that the vertex and edge
groups obtained at the end of the reduction are well-defined. For a precise statement of this
result see Proposition~\ref{PropCollapse}.

\paragraph{Depth zero rafts for irreducible graphs of groups.} Recall that a bounded valence tree
$T$ is \emph{bushy} if each vertex $v \in T$ is a uniformly bounded distance from a vertex $w$
such that that $T-w$ has at least three unbounded components. When $T$ has no valence~1 vertices,
this is equivalent to saying that each vertex $v \in T$ is a uniformly bounded distance from a vertex
of valence~$\ge 3$.

The following result is obtained by specializing the trichotomy of Bass and Kulkarni \cite{BassKulkarni}
to depth zero rafts of an irreducible graph of groups.

\begin{proposition}
\label{PropTrichotomy}
Suppose that $\G$ is a finite type, irreducible graph of groups, with Bass-Serre tree of spaces $X
\to T$. If $R \subset T$ is a depth zero raft then $R$ has no valence~1 vertices, and $R$ is either a point,
a line, or a bushy tree.
\end{proposition}

\begin{proof} The fact that $R$ has no valence~1 vertices is an immediate consequence of
irreducibility of $\G$. The Coboundedness Principle~\ref{PropCoboundednessPrinciple} implies that the
stabilizer group $\Stab(R)$ acts properly discontinuously and cocompactly on the raft space $X_R$, and so
the graph of groups $R / \Stab(R)$ is a geometrically homogeneous graph of groups presentation for the group
$\Stab(R)$. It follows from \cite{BassKulkarni} that $R$ is either a point, a line, or a bushy tree.
\end{proof}

\subsection{Coarse $\PD(n)$ spaces and groups}
\label{SectionPD}

Coarse $\PD(n)$ spaces were introduced in \cite{KapovichKleiner:duality} as the appropriate setting for
arguments of coarse algebraic topology, such as the coarse separation theorem of \cite{FarbSchwartz}, and the
more general coarse Jordan separation theorem and coarse Alexander duality theorem of
\cite{KapovichKleiner:duality}. The material of this section is taken from \cite{KapovichKleiner:duality}, as
extended to CW-complexes in \cite{MSW:QTOne}.

\paragraph{Bounded geometry CW-complexes.} Let $Y$ be a finite dimensional CW-complex. Using induction on
dimension we define what it means for $Y$ to have \emph{bounded geometry}. When $Y$ has dimension~1 this means
the valences of vertices have a finite upper bound. Suppose by induction that $Y$ has dimension~$n$ and
$Y^{(n-1)}$ has bounded geometry. Then $Y$ has bounded geometry if there exists $A > 0$ such that the
following hold:
\begin{itemize}
\item Each point of $Y^{(n-1)}$ touches at most $A$ closed cells of
dimension~$n$.
\item For each $n$-cell $e$ of $Y$ with attaching map $\alpha_e \from S^n
\to Y$, the set $\image(\alpha_e)$ is a subcomplex of $Y$ containing
at most $A$ cells.
\item Up to postcomposition by CW isomorphisms, there are at most $A$ different attaching maps $S^n
\xrightarrow{\alpha_e} \image(\alpha_e)$, as $e$ varies over all $n$-cells of $Y$.
\end{itemize}
Note that if $Y$ is connected then $Y$ supports a geodesic metric in which all \nb{1}cells have
length~1 and all cells of all dimensions have uniformly bounded diameter; indeed, one can arrange that
there are only finitely many isometry types of cells. This metric is well defined up to
quasi-isometry. 

Suppose that $Y$ is a bounded geometry CW-complex. Given a subset $A \subset Y$, define subcomplexes
$N_r A$ by induction as follows: $N_0 A$ is the smallest subcomplex containing $A$; and
$N_{r+1} A$ is the union of $N_r A$ with all closed cells intersecting $N_r A$. In this context,
the subscript $r$ is implicitly assumed to be a non-negative integer. We say that $Y$ is
\emph{uniformly contractible} if for any $r$ there exists $r'$ such that for each subcomplex $L$, if
$\diam(L) < r$ then $L$ is contractible inside $N_{r'} L$. We say that $Y$ is \emph{uniformly
acyclic} if for any $r$ there exists $r'$ such that for each subcomplex $L$, if
$\diam(L) < r$ then the inclusion map $L \inject N_{r'} L$ is trivial on reduced homology.

We remark that the notation $N_r A$ as a subcomplex of $Y$ clashes somewhat with the notation $N_r A$
as a metric neighborhood of $A$ in $Y$. However, the two versions of the set $N_r A$ are at finite
Hausdorff distance, where the distance is bounded by a constant that depends only on $r$ and on $Y$.

\paragraph{Coarse finite type groups.} Recall that a group $G$ is of finite type if it is the fundamental
group of some finite, aspherical CW complex $Y$. Equivalently, $G$ has a free, proper, cocompact, cellular
action on some bounded geometry, uniformly contractible CW complex $X$.

A finitely generated group $G$ is of \emph{coarse finite type} if it is quasi-isometric to some bounded
geometry, uniformly contractible CW-complex $X$, equivalently, $G$ has a proper, cobounded quasi-action on
$X$. By definition, this is a quasi-isometry invariant. Examples of coarse finite type groups include any
finite type group $G$, or any group obtained from $G$ by passage to a finite index subgroup, or to a finite
extension, or to any iteration of these two operations.

Let $H_*(Y)$ denote CW homology and let $H_c^{*}(Y)$ denote compactly supported CW cohomology, with $\Z$
coefficients.

By slightly generalizing the proofs of Gersten \cite{Gersten:dimension} or of Block and Weinberger
\cite{BlockWeinberger}, we obtain the following:

\begin{theorem}
\label{TheoremChainQI}
If $f \from X \to X'$ is a quasi-isometry between bounded geometry, uniformly acyclic CW complexes, then
$f$ induces an isomorphism $f^* \from H^*_c(X) \approx H^*_c(X')$.
\end{theorem}

From this follows the theorem of Gersten and Block--Weinberger, that $H^*(G;\Z G)$ is a quasi-isometry
invariant among groups of finite type. We will need some of the details of the proof of this theorem.

\begin{proof}
By the method of acyclic models, using uniform acylicity, we may construct chain maps $f_\# \from C_*(Y)
\to C_*(Z)$ and $\bar f_\#\from C_*(Z) \to C_*(Y)$ whose supports are uniformly close to $f$, which means for
example that there exists a constant $D$ such that for each simplex $\sigma$ of $Y$, $\supp(f_\#(\sigma))
\csubset{[D]} f(\sigma)$. Moreover, the chain maps $f_\# \composed \bar f_\#$, $\bar f_\# \composed f_\#$ are
chain homotopic to the identity by chain homotopies of bounded support. The cochain maps $f^\# \from C^*(Z) \to
C^*(Y)$ and $\bar f^\#\from C^*(Y) \to C^*(Z)$ which are dual to $f_\#$ and $\bar f_\#$ restrict to cochain
maps $f^\# \from C^*_c(Z)\to C^*_c(Y)$ and $\bar f^\# \from C^*_c(Y)\to C^*_c(Z)$ whose supports are uniformly
close to $f$, $\bar f$, respectively. Moreover, the cochain maps $f^\# \composed\bar f^\#$ and $\bar f^\#
\composed f^\#$ are cochain homotopic to the identity, by cochain homotopies of bounded support. It now
follows that $H^*_c(Y) \approx H^*_c(Z)$. 

If $f'_\# \from C_*(Y) \to C_*(Z)$ is another chain map whose support is uniformly close to $f$, then $f_\#$
and $f'_\#$ are chain homotopic by a chain homotopy $C_*(Y) \to C_{*+1}(Z)$ whose support is uniformly close
to $f$, and similarly for the cochain maps on compactly supported chains. This shows that $f^*$ is well
defined.
\end{proof}

In the above proof, all of the bounds depend only on the uniform acyclicity data of
$Y$ and $Z$, and on the quasi-isometry constants of $f$ and $\bar f$. We say that $f_\#$ and $\bar f_\#$ are
\emph{chain quasi-isometries} induced by $f$ and $\bar f$, and that $f_\#$ and $\bar f_\#$ are \emph{chain
coarse inverses} of each other.

Given a coarse finite type group $G$, the \emph{coarse dimension} of $G$, also called simply the
\emph{dimension}, is the maximum value of $n$ such that $H^n_c(X) \ne 0$, where $X$ is any bounded geometry,
uniformly contractible space quasi-isometric to $G$. By Theorem~\ref{TheoremChainQI}, coarse dimension is
well-defined independent of the choice of $X$, and is a quasi-isometry invariant.

\subparagraph{Questions.} Theorem \ref{TheoremChainQI}, and its connection to the Gersten--Block--Weinberger theorem,
brings up some questions, the hardest of which asks: does every coarse finite type group \emph{act}, properly
discontinuously and cocompactly, on some uniformly contractible bounded geometry CW complex? As a special case, is
every coarse finite type, torsion free group of finite type? The next question, which would follow from a positive
answer to the previous questions, asks: is $H^*(G;\Z G)$ a quasi-isometry invariant, among groups of coarse finite type?
Finally, inspired by Stalling's theorem \cite{Stallings:Dimension1}, is every coarse finite type group of coarse
dimension~1 virtually free?

\paragraph{Coarse $\PD(n)$ spaces.} A bounded geometry, uniformly acylic CW complex $Y$ is said to be a
\emph{coarse $\PD(n)$ space} if there exist chain maps
$$C_*(Y) \xrightarrow{P} C_c^{n-*}(Y) \xrightarrow{\overline P} C_*(Y)
$$
such that the maps $\overline P \composed P$, $P \composed \overline P$
are each chain homotopic to the identity via respective chain homotopies 
$$C_*(Y) \xrightarrow{\Phi} C_{*+1}(Y), \quad C^*_c(Y)
\xrightarrow{\overline\Phi} C_c^{*-1}(Y)
$$
and the supports of these maps are uniformly close to the identity map, that is, there exists a constant $D
\ge 0$ such that for each cell $\sigma$ of $Y$, the supports of each of $P(\sigma), \overline P(\sigma),
\Phi(\sigma), \overline\Phi(\sigma)$ lie in $N_D\sigma$.

\begin{proposition}
\label{PropPDQI}
The coarse $\PD(n)$ property is a quasi-isometry invariant on the class of bounded geometry, uniformly
acyclic CW-complexes.
\end{proposition}

\begin{proof}
Consider a quasi-isometry $f \from Y \to Z$ with coarse inverse $\bar f \from Z \to Y$ where $Y,Z$ are
bounded geometry uniformly acyclic. Let $f_\# \from C_*(X) \to C_*(Y)$ and $\bar f_\# \from C_*(Y) \to C_*(X)$
be chain quasi-isometries induced by $f$, let $f^\# \from  C^*(Y) \to C^*(X)$ and $\bar f^\# \from C^*(X) \to
C^*(Y)$ denote the duals, and let $f^\#_c \from C^*_c(Y) \to C^*_c(X)$ and $\bar f^\#_c \from C^*_c(X) \to
C^*_c(Y)$ denote the restrictions to the compactly supported subcomplexes. There are chain homotopies to the
identity with bounded support of the various maps such as $\bar f_\# \composed f_\#$.

Suppose that $Y$ is coarse $\PD(n)$, and let $P_Y$, $\bar P_Y$, $\Phi_Y$,
$\bar\Phi_Y$ denote the chain maps and chain homotopies in the definition of coarse $\PD(n)$. We may now
define the operators 
$$C_*(Z) \xrightarrow{P_Z} C_c^{n-*}(Z) \xrightarrow{\overline P_Z} C_*(Z)
$$
$$C_*(Z) \xrightarrow{\Phi_Z} C_{*+1}(Z), \quad C^*_c(Z)
\xrightarrow{\overline\Phi_Z} C_c^{*-1}(Z)
$$
by the formulas
\begin{align*}
P_Z      &= \bar f^\#_c \composed P_Y \composed \bar f_\# \\
\bar P_Z &= f_\# \composed \bar P_Y \composed f^\#_c \\
\Phi_Z   &= f_\# \composed \Phi_Y \composed \bar f_\# \\
\bar \Phi_Z &= \bar f^\#_c \composed \bar \Phi_Y \composed f^\#_c
\end{align*}
By a diagram chase, using appropriate chain homotopies of bounded support as mentioned above, these operators
exhibit the coarse $\PD(n)$ property for $Z$.
\end{proof}

A coarse $\PD(n)$ space is defined to be \emph{good} if it is uniformly contractible and its topological dimension is
equal to the formal dimension $n$. Examples of good coarse $\PD(n)$ spaces include bounded geometry, uniformly
contractible spaces which are topological $n$-manifolds. Usually we will need our coarse $\PD(n)$ spaces to be good:
uniform contractibility is used to construct maps as described below, and topological dimension~$n$ is used to apply
certain results of~\cite{MSW:QTOne}.

\paragraph{Coarse $\PD(n)$ groups.} A finitely generated group $G$ is said to be a \emph{coarse
$\PD(n)$ group} if there exists a good coarse $\PD(n)$ space $X$ such that $G$ is quasi-isometric to $X$,
equivalently, $G$ quasi-acts properly and coboundedly on $X$. By definition this is a quasi-isometry invariant. By
convention, finite groups are defined to be coarse $\PD(0)$. 

Examples of coarse $\PD(n)$ groups include any finitely presented $\PD(n)$ group $G$ (see \cite{MSW:QTOne}), or anything
obtained from such a group $G$ by passing to a finite index subgroup, taking a finite extension, or iterating these
operations.

\subparagraph{Question:} Is it true that a coarse $\PD(n)$ group must act, properly and coboundedly, on some good
coarse $\PD(n)$ space?

It is not true, using our definition, that any $\PD(n)$ group is coarse $\PD(n)$: we cannot drop finite
presentability. This is unfortunate because Davis' examples of nonfinitely presentable
\Poincare\ duality groups are not coarse $\PD(n)$ by our definition. This argues for a broader definition in
which we drop the requirement that $X$ be uniformly contractible and require only uniformly acylic; see the
following remarks.

\subparagraph{Remarks on uniform contractibility.} 
Throughout this paper the uniform contractibility property is used chiefly to construct maps. For example,
if $X,Y$ are uniformly contractible and $f \from X \to Y$ is a quasi-isometry, then $f$ is a bounded
distance from a cellular quasi-isometry. This is often called the ``connect the dots'' principle, because of
its proof. First move $f \restrict X^0$ a bounded distance so that its image is in $Y^0$. Then, for each
edge $e$ of $X$, define $f(e)$ by connecting the two dots $f(\bdy e)$, using uniform connectedness to ensure
that the path $f(e)$ has bounded length. Then, for each \nb{2}cell $\sigma$ of $X$, define $f(\sigma)$ by
extending the closed edge path $f(\bdy \sigma)$ to a cellular map of a disc. And so on.

It would be nice to drop uniform contractibility from the definition of coarse $\PD(n)$ groups --- this would
allow us to use the interesting examples of nonfinitely presentable \Poincare\ duality groups due to Davis
\cite{Davis:CoxeterCohomology}. In order to do this, we would have to work with general coarse $\PD(n)$ spaces
rather than just those which are uniformly contractible. The connect-the-dots argument would break down, and
we could no longer depend on the existence of continuous, cellular quasi-isometries. 

Bruce Kleiner has suggested that one ought to get used to working with chain quasi-isometries rather than
quasi-isometries. That is, one ought to replace the category of bounded geometry, uniformly contractible
complexes and quasi-isometries with the category of bounded geometry, uniformly acylic complexes and chain
quasi-isometries. This endeavour would take us too far beyond our present purposes, and so we shall not pursue this
topic further here.

\paragraph{Coarse Bass-Serre complex.} The Bass-Serre tree $T$, regarded as a $\pi_1\G$ tree, is independent of the
choice of a graph of spaces $\K$ over $\G$. The tree of spaces $X$ does indeed depend on the choice of
$\K$. However, the coarse geometric information regarding the quasi-isometry types of the vertex and edge spaces of $X$
and the manner in which they are glued up by coarse uniformly proper embeddings to produce $X$ is independent of $\K$.
To state this more precisely, if $\K'$ is another graph of spaces over $\G$, and if $X' \to T$ is the associated
Bass-Serre tree of spaces, then there is a $\pi_1\G$-equivariant quasi-isometry $X \to X'$ which commutes with the
projections $X \to T$, $X' \to T$. 

The coarse geometric information encoded in a tree of spaces invites us to broaden the concept of a Bass-Serre tree of
spaces $X\to T$ associated to a graph of groups $\G$, in such a way that $\pi_1\G$ may only quasi-act on $X$. This will
be useful when we want to employ model spaces for the vertex and edge groups on which those groups may not strictly act,
but only quasi-act. 

The following lemma captures the key idea that we will need:

\begin{lemma}
\label{LemmaCoarseFiniteType}
If $\G$ is a finite graph of coarse finite type groups, then $\pi_1\G$ has coarse finite type.
\end{lemma}

\begin{proof}
Coarse finite type easily implies finitely presented, so we can construct a finite graph of spaces over $G$, whose
universal cover gives a Bass-Serre tree of spaces $X \to T$ as described above. Given an edge $e$ and an incident
vertex $v$ in $T$, let $F_{ev} \from X_e \to X_v$ denote the ``attaching map'', the lift of the appropriate attaching
map $f_\eta$ of the graph of spaces $\K$, so that $d_X(x,F_{ev}(x))$ is uniformly bounded for all $x \in X_e$ and all
$e,v$. Also, the maps $F_{ev}$ are uniformly proper, with a common properness gauge independent of $e,v$. 

For each vertex or edge $a$ of $T$, since $\Stab(a)$ has coarse finite type, there is a bounded geometry, uniformly
contractible complex $Y_a$ and a quasi-isometry $h_a\from X_a\to Y_a$. Since $T$ has only finitely many orbits of
vertices and edges, we can choose the complexes $Y_a$ to be in finitely many isometry classes, and we can choose the
quasi-isometries $h_a$ to have uniform quasi-isometry constants, and to have coarse inverses $\bar h_a \from Y_a \to
X_a$ with uniform coarse inverse constants, and hence uniform quasi-isometry constants.

For each edge $e$ and incident vertex $v$ in $T$, consider the uniformly proper map $\phi_{ev} \from Y_e \to Y_v$ given
by $h_v \composed F_{ev} \composed \bar h_e$. This collection of maps has a common properness gauge independent of
$e,v$. The spaces $Y_v$ are of uniformly bounded geometry, and they have common uniform contractibility data. By the
connect-the-dots principle, it follows that each of the maps $\phi_{ev}$ can be moved a bounded distance to obtain a
continuous, cellular map $\Phi_{ev}$, so that the maps $\Phi_{ev}$ have a common properness gauge, and so that the
restrictions of $\Phi_{ev}$ to the cells of $Y_e$ have only finitely many different topological types up to pre and
postcomposition. Therefore, by using the maps $\Phi_{ev}$ to glue together the spaces
$(\bigcup_{v \in T} Y_v) \union (\bigcup_{e \subset T} (Y_e \cross e))$, we obtain a tree of spaces $Y$ with the
structure of a bounded geometry, uniformly contractible cell complex, and the maps $h_a$ piece together to give a
quasi-isometry $h\from X \to Y$. It follows that $\pi_1\G$, which is quasi-isometric to $X$, is also quasi-isometric to
$Y$.
\end{proof}

The complex $Y$ constructed in this proof is called a \emph{coarse Bass-Serre complex} for the graph of
groups $\G$. The action of $\pi_1\G$ on $T$ lifts to a proper, cobounded quasi-action of $\pi_1\G$ on $Y$, obtained by
quasiconjugating the action of $\pi_1\G$ on $X$ via the quasi-isometry $h \from X \to Y$.

Henceforth, when $\G$ is a finite graph of coarse finite type groups, and when some of the vertex and edge groups of
$\G$ are coarse $\PD$, then we will assume that the corresponding vertex and edge spaces of the coarse Bass-Serre
complex are coarse $\PD$ as well. Moreover, we will usually subsume the adjective ``coarse'': the distinction between
the (true) Bass-Serre complex and the (coarse) Bass-Serre complex should be evident from the context. This will be
particularly relevant in the proof of Depth Zero Vertex Rigidity given in Section~\ref{CoarsePD}.

\subsection{The methods of proof: the general case.} 
\label{SectionMethodsGeneral}
In this section we formalize the discussions of
vertex rigidity and tree rigidity introduced in Sections~\ref{SectionMethods}
and~\ref{SectionSetting}. We do this by stating three theorems which are the focus of three sections
of the paper: 
\begin{description}
\item[The Depth Zero Vertex Rigidity Theorem~\ref{TheoremVertexRigidity},] covered in
Section~\ref{SectionDepthZero}; 
\item[The Vertex--Edge Rigidity Theorem~\ref{TheoremVERigidity},] covered in
Section~\ref{SectionFiniteDepth}; 
\item[The Tree Rigidity Theorem~\ref{TheoremTreeRigidity},] covered in
Section~\ref{TreeRigidity}. 
\end{description} 
As mentioned earlier, these three theorems will be put together in Section~\ref{MainTheorems} to
prove our main QI-rigidity and QI-classification theorems~\ref{TheoremQI} and~\ref{TheoremClasses}.

Consider $\G$ a finite type graph of groups with Bass-Serre tree of spaces $\pi\from X \to T$.  Recall the
notation $\VE(T)$ for the union of the vertices and the edge midpoints of $T$, regarded as a metric space by
restricting the distance function of $T$. Let $\V_0(T) \subset \VE(T)$ be the set of depth zero vertices of $T$.

\subparagraph{Definition.} The graph of groups $\G$ is \emph{admissible} if it is finite type,
irreducible, and every vertex space is coarsely contained in a depth zero vertex space. 

\bigskip

All the graphs of groups we consider will be admissible. Admissibility often follows from other
assumptions, for instance if every vertex space is coarse Poincare duality, or more generally if the
graph of groups has finite depth. 

Recall the hypotheses of Section~\ref{SectionStatements}. In particular, admissibility of $\G$ is a slight weakening of
the first hypothesis:
\begin{enumerate}
\item[(1)] $\G$ is finite type, irreducible, and of finite depth. 
\end{enumerate}
The other hypotheses are: 
\begin{enumerate}
\item[(2)] No depth zero raft of $T$ is a line.
\item[(3)] The stabilizer of each depth zero vertex of $T$ is coarse $\PD$.
\item[(4)] For each depth zero vertex $v$ of $T$, if $v$ is a one vertex raft then the crossing graph condition holds
at $v$.
\item[(5)] Every edge and vertex group of $\G$ is of coarse finite type.
\end{enumerate}

The following theorem formalizes the discussion of depth zero vertex rigidity in
Section~\ref{SectionSetting}. Note that the theorem does not require the full power of the finite depth property, only
the weaker property of admissibility.

\begin{theorem}[Depth Zero Vertex Rigidity]
\label{TheoremVertexRigidity}
Let $\G,\G'$ be admissible graphs of groups, $X \to T$, $X' \to T'$ their respective Bass-Serre trees of spaces.
Suppose that $\G,\G'$ satisfy hypotheses (2--5). Then any quasi-isometry $f \from X \to X'$ coarsely respects the depth
zero vertex spaces, which means:
\begin{itemize}
\item for any $K \ge 1$, $C \ge 0$ there exists $K' \ge 1$, $C' \ge 0$ such that if $f
\from X \to X'$ is a $K,C$ quasi-isometry, then there exists a $K',C'$ quasi-isometry $f_\#\from\V_0(T)\to\V_0(T')$ such
that:
\begin{itemize}
\item if $v \in\V_0(T)$ then $d_\Haus(f(X_v),X'_{f_\#(v)}) \le C'$, 
\item if $v' \in \V_0(T')$ then there exists $v\in \V_0(T)$ such that $d_\Haus(f(X_v),X'_{v'}) \le C'$.
\end{itemize}
\end{itemize}
\end{theorem}

Next we want to extend the conclusions of Depth Zero Vertex Rigidity to all vertex and edge spaces (this distinction
was blurred somewhat in Sections~\ref{SectionMethods} and~\ref{SectionSetting} where we introduced the methods of
proof). This is where we first use the finite depth property in its full power.

\begin{theorem}[Vertex--Edge Rigidity]
\label{TheoremVERigidity}
Given finite type, finite depth, irreducible graphs of groups $\G,\G'$, with associated
objects as above, if $f \from X \to X'$ is a quasi-isometry that coarsely respects depth zero
vertex spaces, then $f$ coarsely respects vertex and edge spaces, which means:
\begin{itemize}
\item for any $K \ge 1$, $C \ge 0$ there exists $K' \ge 1$, $C' \ge 0$ such that if $f \from X \to X'$ is a $K,C$
quasi-isometry, then there exists a $K',C'$ quasi-isometry $f_\# \from \VE(T) \to \VE(T')$ such that:
\begin{itemize}
\item if $a \in \VE(T)$ then $d_\Haus(f(X_a),X'_{f_\#(v)}) \le C'$.
\item if $a' \in \VE(T)$ then there exists $a \in \VE(T)$ such that $d_\Haus(f(X_a),X'_{a'}) \le C'$.
\end{itemize}
\end{itemize}
\end{theorem}

As we will see in Section~\ref{SectionFiniteDepth}, the Quasi-isometric Classification Theorem~\ref{TheoremClasses}
follows quickly by combining Depth Zero Vertex Rigidity with Vertex--Edge Rigidity. 

Theorem~\ref{TheoremVERigidity} says, intuitively, that the depth zero vertex spaces in $X$ determine the vertex and
edge spaces of any depth, up to Hausdorff equivalence, assuming finite type, finite depth, and irreducibility. When
translated into group theoretic language, it says that the depth zero vertex stabilizers in $\pi_1\G$ determine the
vertex and edge stabilizers of any depth, up to commensurability, again assuming finite type, finite depth, and
irreducibility of $\G$. In Proposition~\ref{PropCollapse} we apply this to obtain a purely group theoretic result,
namely that the process of ``reduction'' of a reducible graph of groups $\G$ produces a well-defined collection of
vertex and edge groups up to commensurability, assuming only that $\G$ is of finite type and finite depth.

\bigskip

Next we turn to the additional arguments needed to prove the Quasi-Isometric Rigidity Theorem~\ref{TheoremQI}. Recall
the setting of that theorem: a graph of groups $\G$ satisfying hypotheses (1--5) above, with Bass-Serre tree of spaces
$X \to T$, and a group $H$ quasi-isometric to $\pi_1\G$. The quasi-isometry $H \to \pi_1\G$ induces a proper, cobounded
quasi-action of $H$ on $X$. By applying the Depth Zero Vertex Rigidity Theorem combined with the Vertex--Edge Rigidity
Theorem, we conclude that the quasi-action of $H$ on $X$ coarsely respects all vertex and edge spaces of $X$. This
conclusion, together with hypotheses (1) and (2) for $\G$, is the starting point for the next theorem. 

\begin{theorem}[Tree rigidity theorem]
\label{TheoremTreeRigidity}
Let $\G$ be a finite type, finite depth, irreducible graph of groups, with Bass-Serre tree of spaces $X \to T$, and
suppose that no depth zero raft of $T$ is a line. Let $H$ be a finitely generated group quasi-acting properly and
coboundedly on $X$, coarsely respecting vertex and edge spaces of $X$. Then there exists a finite type, finite
depth, irreducible graph of groups $\G'$ with Bass-Serre tree of spaces $X' \to G'$, and there exists an isomorphism
$\pi_1\G' \approx H$, such that any induced quasi-isometry $X \to X'$ coarsely respects vertex and edge spaces.
\end{theorem}
\vfill\break

\newcommand\ds\displaystyle

\section{Depth Zero Vertex Rigidity}
\label{SectionDepthZero}

The main goal of this section is to prove the Depth Zero Vertex Rigidity
Theorem~\ref{TheoremVertexRigidity}.  Let $\G$ and $\G'$ be finite type, irreducible graphs of
groups, $X\mapsto T$ and $X' \mapsto T'$ associated trees of spaces, and $f \from X \to X'$ a
quasi-isometry. Recall the notation $\V_0(T)$ for the depth zero vertices of $T$ with metric obtained by
restriction of the metric on $T$, and similarly for $T'$. The desired conclusion is:

\begin{description}
\item[Depth Zero Vertex Rigidity] There exists a $K',C'$ quasi-isometry $f_\# \from\V_0(T) \to \V_0(T')$
such that for each $a\in \V_0(T)$ we have $f(X^{\vphantom{\prime}}_{\vphantom{a'}a}) \ceq{[C']}
X'_{f_\#(a)}$, and for every $a' \in \V_0(T')$ there exists $a \in \V_0(T)$ such that $X'_{f_\#(a)}
\ceq{[C']} X'_{a'}$. The constants $K' \ge 1$, $C' \ge 0$ depending only on $\G,\G'$ and on the
quasi-isometry constants $K \ge 1$, $C\ge 0$ of $f$.
\end{description}
In other words, $f$ induces a bijection between the Hausdorff equivalence classes of depth
zero vertex spaces in $X$ and the Hausdorff equivalence classes of depth zero vertex spaces
in $X'$. In particular, each depth zero vertex space of $\G$ is quasi-isometric to
some depth zero vertex space of $\G'$ and vice versa.  

The statement of the Depth Zero Vertex Rigidity Theorem~\ref{TheoremVertexRigidity} contains several hypotheses
on $\G$ and $\G'$. To understand these hypotheses it is helpful to consider some examples where they fail. 
The simplest problem is that a group can have several splittings as a graph of groups.  Most trivially this
can happen simply by starting with a graph of groups decomposition and collapsing a subgraph to a single
vertex.   Thus we generally must make sufficient assumptions about our graph of groups to assure that they do
not further split, at least within the category of graphs of groups under discussion.  This is one of the
benefits of assuming the vertex groups to be coarse Poincare duality groups --- they have limited
splittings.   

For another kind of example, if $G$ is a Poincare duality group then so is $G \times \Z$. 
We can write $G \times \Z$ as an HNN extension of $G$.  The Bass-Serre tree of this
splitting is a line, with all the vertex and edge spaces depth zero. Thus $G \times \Z$ has
two descriptions as a graph of PD groups: the HNN splitting just described; and the graph
with one vertex labelled $G \times \Z$ and no edge. Further, this problem can be embedded
inside bigger graphs of groups: consider $(G \times \Z) * H $ for some $H$. We get around
these sorts of examples by assuming that there are no line-like rafts, and making
assumptions on vertex spaces which are incident to codimension one edge spaces.

Another difficulty is exhibited by the the second example in Section~\ref{universal}. In this example, by
passing to a finite index subgroup one can collapse edges and eliminate Hausdorff equivalence classes of
vertex spaces, and as a result there seem to be no canonically defined vertex spaces.  To get around this
sort of problem we only discuss the depth zero vertex spaces; in the second example of
Section~\ref{universal} there are none.  Restricting to depth zero vertex spaces is not sufficient on its
own, as this example can be embedded as a subgraph of another graph of groups, and the ambiguity of
splitting in the subgraph prevents vertex rigidity from holding for the other, possibly depth zero, vertex
spaces. To get around this we assume admissibility, which says that every vertex space is coarsely
contained in some depth zero vertex space.  This guarantees that there are enough depth zero vertex
spaces to coarsely control the shape of the tree of any splitting.

\paragraph{An outline of the proof of depth zero vertex rigidity.} Depth zero vertex rigidity for
quasi-isometries between $X$ and $X'$ is essentially equivalent to the statement that under a quasi-isometry
$X \to X'$, the images of depth zero vertex spaces ``cannot cross'' edge spaces.

The first step is to focus on the property of ``depth zero raft rigidity'', which says that every
quasi-isometry $X \to X'$ coarsely respects depth zero rafts. In the Abstract Raft Rigidity
Theorem~\ref{AbstractRaftRigidity} we give a sufficient condition for depth zero raft rigidity called ``mutual
thickness'', which says essentially that vertex spaces in one of $X$ or $X'$ cannot cross positive depth edge
spaces in the other. To handle depth zero edge spaces we invoke some of the hypotheses of our main main
theorems, namely, that all depth zero vertex groups are coarse $\PD$ and no depth zero raft is a line. These
hypotheses, combined with mutual thickness, give a sufficient condition for depth zero vertex rigidity, stated
in Corollary~\ref{GDZVR}.

The second step is to use the hypotheses of the depth zero vertex rigidity theorem to establish the property
of mutual thickness. This is where we apply the Coarse Alexander Duality Theorem
\cite{KapovichKleiner:duality}.

\subsection{A sufficient condition for depth zero vertex rigidity.}
\label{SectionAbstractVertexRigidity}

In this section we consider trees of spaces $X \mapsto T$, $X' \mapsto T'$ associated to admissible graphs
of groups $\G,\G'$. In order for depth zero vertex rigidity to hold, it is necessary that the depth zero
rafts also exhibit a rigidity property:

\begin{definition}  We say \emph{depth zero raft rigidity} holds for quasi-isometries
between $X$ and $X'$ if, for all $K,C$ there is an $R$ so that for every $K,C$
quasi-isometry $f\from X \to X'$, there is a bijection $f_\#$ between the depth zero rafts of
$T$ and of $T'$ so that for any such raft~$\Sigma$ of~$T$ we have $f(X_\Sigma)
\ceq{[R]} X'_{f_\#(\Sigma)}$.
\end{definition}

While this property is not sufficient for depth zero vertex rigidity, the only difference
between the two properties is simply what happens within a depth zero raft. A depth zero raft
is the Bass-Serre tree of a homogeneous graph of groups in the sense of \cite{MSW:QTOne}. Within the raft,
vertex rigidity holds trivially if the raft is bounded, and it also holds if the raft is bushy with coarse
$\PD$ vertex spaces \cite{FarbMosher:ABC}. In other words, if depth zero vertices are coarse
$\PD$, and if no rafts are lines, both of which are guaranteed by the hypotheses of the depth zero vertex
rigidity theorem, then depth zero vertex rigidity holds within the depth zero rafts. In this situation we will
show, in the proof of Corollary~\ref{GDZVR}, that depth zero raft rigidity does indeed imply
depth zero vertex rigidity.

The Abstract Raft Rigidity Theorem~\ref{AbstractRaftRigidity} gives sufficient conditions for depth zero raft
rigidity, which we now discuss.

\paragraph{Crossing and thickness.} The basic idea in the proof of vertex rigidity is that
a subset of a tree of spaces that is not coarsely contained in a vertex space must cross some
edge space. Often one can rule this out by examining how an edge space disconnects the
tree of spaces and subsets thereof. We formalize this as follows.

\begin{definition}  
Let $A$ and $B$ be subsets of $X$.  We say that $A$ \emph{crosses} $B$ if, for every sufficiently large $r>0$,
$A$ has deep intersection with at least two components of $X - N_r B$.
\end{definition}

Note in particular that, for sufficiently large $r$, $N_r B$ must disconnect $X$, and the
set $N_r B \intersect A$ disconnects $A$ into two deep subsets of $A$. This definition is
formulated so that ``$A$ crosses $B$'' is a quasi-isometrically invariant property.

Our interest is mainly in the case where $B$ is an edge space $X_e$. In that case, $T
- \midpoint(e)$ has two components $T^+$ and $T^-$, and there is a corresponding decomposition
of $X - X_e =  \U^+_e \cup \U^-_e$. For now the assignment of the $+$ and $-$ signs is arbitrary, but
these signs will be given geometric significance in the proof of Theorem~\ref{AbstractRaftRigidity}. Any
subset $A$ which goes arbitrarily deep into both of $\U^\pm_e$ certainly crosses $X_e$. 

In sifting among several possible variations of the concept of ``crossing'' for the best definition, it is
important to keep in mind that the property ``$A$ intersects each of $\U^\pm_e$ arbitrarily deeply'' is
strictly stronger than ``$A$ crosses $X_e$'', and is not quasi-isometrically invariant. For example, consider
the graph of groups $\Gamma$ with one vertex group $\Z^2$ and one edge group $\Z$ attached on both ends to the
same $\Z$ subgroup of $\Z^2$. In the corresponding tree of spaces $X \to T$, for each vertex $v$ and
incident edge $e$ of $T$, the vertex space $X_v$ crosses $X_e$, but $X_v$ intersects exactly
one of the sets $\U^\pm_e$. Notice also that in this example, $\pi_1\Gamma$ is quasi-isometric
to the free product of a free group with an infinite cyclic group, and there are
quasi-isometries of $X$ which do \emph{not} coarsely respect the vertex spaces.

On the other hand, the separation hypotheses in our theorems are designed to say, essentially,
that for any positive depth edge $e$ and each $r>0$, each depth zero vertex space is coarsely
contained in a single component of $X - N_r(X_e)$. In other words, the depth zero vertex
spaces do not cross edge spaces in the sense of the above definition.

\bigskip

Consider now $X \mapsto T$ and $X' \mapsto T'$, trees of spaces associated to two
admissible graphs of groups. 

\begin{definition} Given $Y \subset X$ and an edge $e$ of $T'$, we say \emph{$Y$ cannot cross
$e$} if, for any quasi-isometry $f\from X \to X'$, $f(Y)$ does not cross $X'_e$. We say that $Y$ is
\emph{$e$-thick} if $Y$ cannot cross $e$, and if there is no quasi-isometry $f\from X \to X'$ such that $f(V)
\csubset{c} X'_e$. 
\end{definition}

Informally, $Y$ is $e$-thick if the number of complementary components of $X'_e$ that the image of $Y$ deeply
intersects is precisely one.

Note that in the definition of ``$e$-thickness'', one obtains a weaker and quasi-isometrically unnatural
property by replacing the statement ``$f(Y)$ does not cross $X'_e$'' with the statement ``$f(Y)$ does not
deeply intersect both components of $X - X'_e$''; see the discussion above regarding variations on the
definition of ``crossing''.

\begin{definition} Let $X$ and $X'$ be as above. We say that $X$ and $X'$ are \emph{mutually thick} if
every depth zero vertex space of $X$ is $e'$-thick for every positive depth edge $e'$ of $X'$ and
likewise, every depth zero vertex space of $X'$ is $e$-thick for every positive depth edge $e$ of $X$. 
\end{definition}

\paragraph{Abstract raft rigidity.} The abstract raft rigidity theorem says that mutual thickness is a
sufficient condition for depth zero raft rigidity. As shown in Corollary~\ref{GDZVR}, depth zero vertex
rigidity follows from depth zero raft rigidity combined with vertex rigidity within the rafts.

\begin{theorem}[Abstract Raft Rigidity Theorem]
\label{AbstractRaftRigidity} 
Let $X$ and $X'$ be as above. If $X$ and $X'$ are mutually thick, then depth zero raft rigidity holds for
quasi-isometries $f \from X \to X'$.
\end{theorem}

We begin with some preliminary lemmas.
 
Mutual thickness guarantees that for any depth zero vertex space $X_v$ of $X$, positive depth edge $e$ of $X'$ and $f\from X
\to X'$ a quasi-isometry, there is an $R$ so that $f(X_v)$ is contained in the $R$-neighborhood of one side of $X'_e$.  Depth
zero raft rigidity requires that $R$ be uniform over all choices of $v$ and $e$, and all $(K,C)$-quasi-isometries. Since there
are only finitely many orbits of vertices and edges, this uniformity will follow from:

\begin{lemma}  
\label{LemmaUniformThickness}
Given a depth zero vertex space $X_v$ of $X$ and an edge $e$ of $T'$, if $X_v$ is $e$-thick then for
every $K,C$ there is an $R$ so that for all $(K,C)$-quasi-isometries $f \from X \to X'$, $f(X_v)$ is
contained in the $R$-neighborhood of some component $\U^\pm_e$ of $X' - X'_e$. 
\end{lemma}

\begin{proof} Let $\U^\pm_{e,r} = \U^\pm_e - N^\vp_r(X'_e)$. Note that the complement in $X'$ of the $r$-neighborhood
of~$\U^\pm_e$ is~$\U^\mp_{e,r}$. Any subset of $X'$ that intersects both of the sets $\U^\pm_{e,r}$ is said
to \emph{weakly cross} the set $N^\vp_r(X'_e)$. 

By assumption, for any quasi-isometry $f \from X \to X'$ there exists an $R$ so that
$f(X_v)$ does not weakly cross $N^\vp_R(X'_e)$. We need to see that this $R$ can be chosen to depend only on
the quasi-isometry constants of $f$. If not then, fixing $K,C$, there is a sequence $f_i \from X \to
X'$ of $(K,C)$ quasi-isometries such that as $i \to +\infinity$, the infimum of all $R$ that work for $f_i$
tends to $+\infinity$. To be precise, there exists a sequence $R_i \to +\infinity$ such that if $r>R_i$
then $f_i(X_v)$ does not weakly cross $N^\vp_{r}(X'_e)$, but if $r<R_i$ then $f_i(X_v)$ weakly crosses
$N^\vp_r(X'_e)$. By changing $C$, perturbing each $f_i$ a uniform amount, and perturbing each $R_i$ a uniform
amount, we can assume that each $f_i$ is continuous. For the rest of the proof we fix
such a sequence $f_i$. 

By applying the Coarse Convergence Principle~\ref{PropCC} to the sequence $f_i$,
suitably normalized, we will obtain a limiting $K,C$ quasi-isometry $f_\infinity\from X\to X'$ which is
used to contradict $e$ thickness of $X_v$. In Case~1 this is done by showing
that $f_\infinity(X_v)\csubset{[r]} X'_e$ for some $r$, and in Case~2 by showing that $f_\infinity(X_v)$
crosses $X'_e$.

Given $x \in X_v$ and $s  \ge 0$, let $B(x,s)$ denote the closed ball of radius $s$ about $x$ in the
path metric on $X_v$. 

\subparagraph{Case 1:} Suppose that there exists $r>0$ such that for every $s>0$, there exists $i \ge 0$
and $x_i \in X_v$ such that $f_i(B(x_i,s)) \subset N_r(X'_e)$. By pre-composing with elements of
$\Stab{v}$ and post-composing with elements of $\Stab{e}$, we may assume that the sequence $\{x_i\}$ is
bounded in $X$ and the sequence $\{f_i(x_i)\}$ is bounded in $X'$. Passing to a subsequence and applying
the Coarse Convergence Principle~\ref{PropCC}, we have a limit point $x_\infinity = \displaystyle\lim_i
x_i$ and an $A$-coarse limit quasi-isometry $f_\infinity$. Moreover, for each $s>0$ the map
$f_\infinity$ takes $B(x_\infinity,s)$ into $\overline N^{\vphantom\prime}_{r+A}(X'_e)$. It follows that
$f_\infty(X_v)\subset \overline N^{\vphantom\prime}_{r+A}(X'_e)$, violating $e$-thickness of~$X_v$.

\subparagraph{Case 2:} Suppose that for every $r>0$ there exists $s=s(r)>0$ such that for every $i \ge
0$ and $x \in X_v$, letting $B=B(x,s(r))$ we have $f_i(B) \not\subset N_r(X'_e)$. In other words, the
image of every $s(r)$-ball $B$ under every $f_i$ intersects at least one of the two sets $\U^\pm_{e,r}$.
Clearly we can choose the function $s(r)$ to be monotonically increasing.

Note that $s(r) \to +\infinity$ as $r \to +\infinity$. Otherwise, $s(r)$ is bounded above by some constant
$S$. Letting $A = KS+C$, the image of every $S$ ball under every $f_i$ has diameter $\le A$ in $X'$, since $f_i
\restrict X_v \from X_v \to X'$ is $K,C$ coarse Lipschitz. On the other hand, if $i$ is sufficiently large
then $f_i(X_v)$ weakly crosses $X'_e$, and so there exists $x \in X_v$ such that $f_i(x) \in X'_e$; this
follows by continuity of $f_i$.  We therefore have $f_i(B(x,S)) \subset N^\vp_{A}(X'_e)$, and so
$f_i(B(x,s(r))) \subset N^\vp_{R}(X'_e)$ for $r=A$, violating the hypothesis of Case 2.

\begin{claim} 
\label{ClaimBallCross} 
\begin{enumerate}
\item \label{ItemBallCross}
For each $r$ there exists $i$ and there exists a ball $B \subset X_v$ of radius $s(r)$ so that
$f_i(B)$ weakly crosses $N^\vp_{r}(X'_e)$.
\item \label{ItemExhaustingCross}
There exists an increasing sequence $r_n \to \infinity$ so that, letting $s_n=s(r_n)$, the following
holds: given a ball $B' \subset X_v$ of radius $s_{n+1}$, if $f_i(B')$ weakly crosses $N^\vp_{r_{n+1}}(X'_e)$,
then there exists a ball $B \subset X_v$ of radius $s_n$ so that $f_i(B)$ weakly crosses $N^\vp_{r_n}(X'_e)$
and so that $B \subset B'$. 
\end{enumerate}
\end{claim}

\begin{proof} To prove item~\ref{ItemBallCross}, for each $r>0$, we know that there exists $i$ such
that $f_i(X_v)$ intersects both of the sets $\U^\pm_{e,r}$. Choose balls $B_0,B_1 \subset X_v$ of
radius $s(r)$ so that $f_i(B_0) \intersect \U^-_{e,r} \ne \emptyset$ and $f_i(B_1) \intersect
\U^+_{e,r} \ne \emptyset$. Choose a path $p \from [0,1] \to X_v$ from the center of $B_0$ to the
center of $B_1$. Let $B_t = B(p(t),s_n)$. As $t$ varies, the closed balls $B_t$ and their images
$f_i(B_t)$ vary continuously in the Hausdorff topology. For each $t$, $f_i(B_t)$ intersects at least
one of the closed sets $\U^\pm_{e,r}$, and so there must exist a time $t \in [0,1]$ at which
$f_i(B_t)$ intersects both of $\U^\pm_{e,r}$.

The proof of item~\ref{ItemExhaustingCross} uses the same ideas as item~\ref{ItemBallCross}, except
that we arrange to choose the balls $B_t$ to all be contained in $B'$. Suppose by induction that
$r_1,\ldots,r_n$ have been chosen. Choose $r_{n+1}$ so large that the following holds: for any $K,C$
quasi-isometry $f\from X \to X'$, and any ball $B \subset X_v$ of radius $s_n$, the diameter of $f(B)$
is less than $r_{n+1} - r_n$. Suppose that $B'\subset X_v$ has radius $s_{n+1}$ and $f_i(B')$ crosses
$N^\vp_{r_{n+1}}(X'_e)$. Choose $y^\pm \in B'$ so that $f_i(y^\pm) \subset \U^\pm_{e,r_{n+1}}$. Choose
$X_v$ geodesics $\overline{y^\pm x}$ and let $z^\pm \in \overline{y^\pm x}$ be the point at distance
$s_n$ from $y^\pm$. The set $f_i\left(B(z^\pm,s_n)\right)$ has nontrivial intersection with
$\U^\pm_{e,r_{n+1}}$ and so is entirely contained in $\U^\pm_{e,r_n}$, by choice of $r_{n+1}$. Let $B_0 =
B(z^-,s_n)$ and $B_1 = B(z^+,s_n)$, and now proceed as in the proof of item~\ref{ItemBallCross}, using
the path $p(t) = \overline{z^- x} * \overline{x z^+}$ to define $B_t = B(p(t),s_n)$. Clearly $B_t
\subset B'$. By the proof of item~\ref{ItemBallCross}, for some~$t$ the set $f_i(B_t)$ crosses
$N^\vp_{r_n}(X'_e)$.
\end{proof}

To finish Case~2, for each $n$ apply Claim~\ref{ClaimBallCross} to choose $i_n$ and balls $B_{n,n}
\supset B_{n,n-1} \supset \cdots \supset B_{n,1}$ in $X_v$ of respective radii $s_n,s_{n-1},\ldots,s_1$ 
such that $f_{i_n}(B_{n,j})$ weakly crosses $N^\vp_{r_j}(X'_e)$. By precomposing with elements of $\Stab(X_v)$
we may assume that the center of $B_{n,1}$ lies in a predetermined compact set independent of $n$. It
follows that for each $j$ the centers of the balls $B_{n,j}$ lie in a predetermined compact set
independent of $n$ and depending on $j$. Since the group $\Stab(X'_e)$ acts coboundedly on
$N^\vp_{r_1}(X'_e)$, by postcomposing with elements of $\Stab(X'_e)$ we may assume that
$f_{i_n}(B_{n,1})$ lies in a predetermined compact set independent of $n$. It follows that for each $j$
the balls $f_{i_n}(B_{n,j})$ lie in a predetermined compact set independent of $n$ and depending on
$j$. Passing to a subsequence, the balls $B_{n,1}$ converge in the Hausdorff topology to a ball $B_1$
and the functions $f_{i_n}$ coarsely converge to a $K,C$ quasi-isometry which takes
$B_{n,1}$ to a subset that weakly crosses $N^\vp_{r_1}(X'_e)$. Continuing inductively to choose subsequences,
and then diagonalizing, we obtain the following result: for each $j$ the balls $B_{n,j}$ converge in the
Hausdorff topology to a ball $B_j$, and the functions $f_{i_n}$ $A$-coarsely converge to a $K,C$
quasi-isometry $f_\infinity$, so that $f_\infinity(B_j)$ weakly crosses $N^\vp_{r_j}(X'_e)$ for each $j$. 
Since $r_j \to +\infinity$, it follows that $f(X_v)$ crosses $X'_e$, contradicting $e$-thickness of $X_v$. This
finishes Case 2 and completes the proof of Lemma~\ref{LemmaUniformThickness}.
\end{proof}

Next we state a local version of depth zero raft rigidity:

\begin{lemma}[Local raft rigidity] 
\label{LemmaLocalRaftRigidity}
For each quasi-isometry $f \from X \to X'$, the following holds: 

Consider a depth zero raft $\Sigma \subset T$, and suppose that for each vertex $v \in \Sigma$, and each positive depth edge
$e$ of $T$, $X_v$ is $e$-thick. Then there exists a depth zero raft $\Sigma' \subset T'$ such that $f(X_\Sigma) \csubset{[r]}
X'_{\Sigma'}$. The constant $r$ depends only on the quasi-isometry constants of $f \from X \to X'$.
\end{lemma}

The meaning of ``local'' in this statement is that it applies to one depth zero raft at a time, regardless
of what is known about other depth zero rafts. This generality will be useful in later sections.

Before proving local raft rigidity, we apply it to:

\begin{proof}[Proof of Theorem~\ref{AbstractRaftRigidity}] Fix a quasi-isometry $f \from X \to X'$. Note
that the hypotheses of Theorem~\ref{AbstractRaftRigidity} imply that the hypotheses of Lemma~\ref{LemmaLocalRaftRigidity} hold
for each depth zero raft $\Sigma \subset T$, and so the conclusions hold for each $\Sigma$, and so we may define
the map $f_\#$ by $f_\#(\Sigma) = \Sigma'$. A similar argument applies to a coarse inverse $\bar f \from X' \to X$, giving a map
$\bar f_\#$ with the property that $f(X'_{\Sigma'}) \csubset{c} X_{\bar f_\#(\Sigma')}$, with uniform coarse inclusion
constant. We therefore have $X_\Sigma \csubset{c} X_{\bar f_\# \composed \bar f(\Sigma)}$ which implies that $\Sigma = \bar
f_\# \composed \bar f(\Sigma)$, because no raft space in $X$ is coarsely contained in any other. A similar argument in the
other direction proves that $f_\#$ and $\bar f_\#$ are inverse bijections, completing the proof of
Theorem~\ref{AbstractRaftRigidity}.
\end{proof}

\begin{proof}[Proof of local raft rigidity] For each edge $e$ of $T'$, as usual let the two components of
$X'-X'_e$ be denoted $\U^\pm_e$. Since $X_v$ is $e$-thick for positive depth edges of $T'$, for any
quasi-isometry $f\from X \to X'$ the set $f(X_v)$ is coarsely contained in precisely one complementary
component $\U^-_e$ or $\U^+_e$ of each positive depth edge space $X'_e$. Since there are only finitely
$\pi_1\G'$ orbits of edges in $T'$, Lemma~\ref{LemmaUniformThickness} shows that there is an $R$ depending
only on the quasi-isometry constants of $f$ so that $f(X_v) \csubset{[R]} \U^-_e$ or
$\U^+_e$.

Orient each positive depth edge $e$ of $T'$ to point toward the component of $X' - X'_e$ containing $f(X_v)$;
in other words, $\U^+_e$ is the component of $X-X'_e$ whose $R$-neighborhood contains $f(X_v)$. Note that if
$e,e'$ share a vertex $w$, and if $e$ points away from $w$, then $e'$ points toward $w$. This shows that $T'$
has a unique sink towards which every edge of $T'$ points: the sink is either a depth zero
raft of $T'$, or an end of $T'$, or a positive depth vertex of~$T'$. 

First we rule out the case where the sink is a positive depth vertex $w \in T'$. Let $e$ be a positive
depth edge of $T'$ incident to $w$ so that $X'_w \ceq{c} X'_e$. Each edge $e'$ incident to $w$ is oriented
towards $w$, and we have $f(X_v) \subset N^\vp_R(\U^+_{e'})$. But the intersection of the sets
$N^\vp_R(\U^+_{e'})$, over all edges $e'$ incident to $w$, is coarsely equivalent to~$X'_w$. It follows
that $f(X_v) \csubset{c} X'_w \ceq{c} X'_e$, contradicting $e$ thickness of $X_v$.

Next we rule out the case where the sink is an end $\eta \in \Ends(T')$. Choose a vertex $w$ of $T'$ such
that $f(X_v)$ actually intersects $X'_w$. Let $\rho$ be the ray in $T'$ pointing from $w$ to $\eta$.
Choose an edge $e$ on $\rho$ so that $X'_w$ has distance greater than $2R$ from $X'_e$. Note that $X'_w
\subset \U_e^-$, and so $f(X_v)$ intersects $\U_e^- - N^\vp_{2R}(X'_e)$. But by the definition  of the
orientation, we have $f(X_v) \subset N^\vp_R(\U^+_e)$, which is a contradiction because the sets
$\U_e^- - N^\vp_{2R}(X'_e)$ and $N^\vp_R(\U^+_e)$ are disjoint. 

We have proved that for each depth zero vertex $v$ of $T$, the sink of $T'$ associated to $v$ is a depth
zero raft $\Sigma\subset T'$. Letting $s$ be the minimum distance between an edge space and any incident
vertex space in $T'$, and letting $r=R+s$, it follows that $f(X_v) \csubset{[r]} X'_\Sigma$. 

Let $u,v$ be depth zero vertices in the same raft $\Sigma$ of $X$. Let $\Sigma',\Sigma''$ be depth zero rafts of $T'$ such that
$f(X_u) \csubset{[R+s]} X'_{\Sigma'}$, $f(X_v)\csubset{[R+s]} X'_{\Sigma''}$. Since $X_{u} \ceq{c} X_{v}$, it follows that
$f(X_{u}) \ceq{c} f(X_{v})$ are contained in the coarse intersection $X'_{\Sigma'} \cintersect{c} X'_{\Sigma''}$. If $\Sigma'
\ne \Sigma''$ then this coarse intersection is coarsely contained in $X'_e$ for some positive depth edge $e$ of $T'$. But this
contradicts $e$ thickness of $u,v$, which implies that neither $f(X_u)$ nor $f(X_v)$ is coarsely contained in $X'_e$. Thus we
have that for each raft $\Sigma$ of $T$ there exists a \emph{unique} raft $\Sigma'$ of $T'$ such that $f(X_v)\csubset{[R+s]}
X'_{\Sigma'}$ for each vertex $v \in \Sigma$. There is a uniform upper bound $d$ to the distance from each point of
$f(X_\Sigma)$ to $f(X_v)$ for some vertex $v \in \Sigma$, and so we have $f(X_\Sigma)\csubset{[R+s+d]} X'_{\Sigma'}$. 
\end{proof}

\begin{corollary}\label{GDZVR}  
Let $X$ and $X'$ be mutually thick trees of spaces of admissible graphs of groups. If vertex rigidity holds
for depth zero rafts then depth zero vertex rigidity holds.
\end{corollary}

\begin{proof} Applying Theorem~\ref{AbstractRaftRigidity} to get depth zero raft rigidity, and combining
with vertex rigidity for depth zero rafts, we obtain an $R_2$ depending only on $K,C$ so that for
any $(K,C)$ quasi-isometry $f \from X \to X'$, the image of each depth zero vertex space $X_v$ is
contained in the $R_2$ neighborhood of some depth zero vertex space $X'_{f_\#(v)}$. We need to see that
the map $f_\# \from\V_0(T) \to\V_0(T')$ is a quasi-isometry. Since $\pi_1\G$ and $\pi_1\G'$ act with
bounded quotient on $X$ and $X'$, there is $r$ so that every vertex space of $X$ or $X'$ is contained in
the $r$-neighborhood of a depth zero vertex space. In particular, the depth zero vertices are
$r$-dense in $T$ and $T'$.

Let $u$ and $v$ be depth zero vertices of $T$, and let $u=w_0, \cdots , w_D = v$ be the vertices
on a geodesic between $u,v$ in $T$, with $D=d_T(u,v)$.  We can find a sequence of lifts to~$X$,
$a_0,a_1,\cdots,a_{D-1}$ and $b_1,\cdots,b_D$, such that $a_i$ and $b_i$ are in the vertex space
$X_{w_i}$ and $d(a_{i-1},b_{i})\leq 1$ for all appropriate values of $i$.  For each
$i$, choose a depth zero vertex space $X_{y_i}$ whose $r$-neighborhood contains $X_{w_i}$.  By moving
$a_i$ and $b_i$ at most $r$, we get sequences $\hat a_i$ and $\hat b_i$, with $\hat a_i$ and
$\hat b_i$ in $X_{y_i}$ and $d(\hat a_{i-1},\hat b_i)\leq 2r+1$ for all appropriate~$i$. 

Consider the images of these points.  The points $f(\hat a_i)$ and $f(\hat b_i)$ are in
the $R_2$ neighborhood of the depth zero vertex spaces $X'_{f_\# (y_i)}$, and
$$d(f(\hat a_{i-1}),f(\hat b_i))\leq K(2r+1)+C
$$
The sequence of points $f_\#(y_i)$ therefore gives a path in $T'$ from $f_\#(u)$ to
$f_\#(v)$ of length at most $D(K(2r+1)+C+2R_2)$. Thus $f_\#$ is coarse lipschitz.  The same
argument applies to $f^{-1}$ to show that $f^{-1}_\#$ is coarse lipschitz.  As they are also
coarse inverses, this proves they are quasi-isometries.
\end{proof}

\subsection{Proof of the Depth Zero Vertex Rigidity Theorem}
\label{CoarsePD}

The results of Section~\ref{SectionAbstractVertexRigidity} show how to derive depth zero vertex rigidity from
the assumption of mutual thickness. In this section we show that when the hypotheses of the depth zero vertex
rigidity theorem hold, the Coarse Alexander Duality theorem of \cite{KapovichKleiner:duality} can be used to
prove mutual thickness.

We shall prove two propositions, \ref{PropCodim2Thick} and~\ref{PropThickCrossing}, each of which is
illustrated by applying it to a special case of the depth zero vertex rigidity theorem. These two propositions
illustrate two different aspects of mutual thickness: the case when an edge space has codimension $\ge 2$
with respect to a vertex space; and the case of codimension~1. The proof of each proposition uses coarse
Alexander duality. We will then combine these two propositions with an inductive argument to prove depth zero
vertex rigidity in full generality.

Throughout this section, when $\G$ is a graph of groups whose vertex and edge groups have coarse finite type some of which are
coarse $\PD$, the terminology ``Bass-Serre complex'' will refer to the bounded geometry, uniformly contractible, coarse Bass-Serre
complex on which $\pi_1\G$ quasi-acts properly and coboundedly, with all edge and vertex spaces being of bounded geometry and
uniformly contractible, and with appropriate edge and vertex spaces being coarse $\PD$ spaces,
as constructed in Section~\ref{SectionBassSerre}. We will usually abuse terminology by speaking about the ``action'' of $\pi_1\G$
on its (coarse) Bass-Serre complex, when only a quasi-action is meant.

\paragraph{The codimension $\ge 2$ case.} Here is one special case of depth zero vertex rigidity, which
can be applied to the class of examples discussed in Theorem~\ref{Neat}:

\begin{corollary}\label{codim2rigidthm}  If $\G$ and $\G'$ are graphs of groups such that
every depth zero vertex group is a coarse $\PD(n)$ group and every positive depth edge group
has dimension at most $n-2$ then depth zero raft rigidity holds, and if neither $\G$ nor
$\G'$ has line-like rafts then depth zero vertex rigidity holds.
\end{corollary}

Corollary~\ref{codim2rigidthm} is an immediate consequence of Corollary~\ref{GDZVR} and the following general
statement:

\begin{proposition}\label{PropCodim2Thick} Let $\G$ and $\G'$ be admissible graphs of groups whose depth zero
vertex groups are all coarse $\PD$. Let $X \to T$, $X' \to T'$ be the Bass-Serre trees of spaces for
$\G,\G'$ respectively. If $v$ is a depth zero vertex of $T$ so that $X_v$ is coarse $\PD(n)$, and if
$e$ is an edge of $T'$ so that $X'_e$ has coarse dimension $\le n-2$, then $X_v$ is $e$-thick.
\end{proposition}

The proof will depend on a lemma about subspaces of codimension $\ge 2$:

\begin{lemma}
\label{LemmaCodim2Fatness} If $Y$ is a coarse $\PD(n)$ space and $Z$ is a uniformly contractible, bounded
geometry space of coarse dimension $\le n-2$, then no subset $B \subset Z$ is uniformly
equivalent to a coarsely separating subset of $Y$, nor to $Y$ itself. 

To be more quantitative: if $A \subset Y$ and if $f \from A \to B$ is a uniformly proper equivalence then there
exists an increasing function $R(r)$ depending only on $Y$, $Z$, and the uniform properness data of $f$, such
that for each $r\ge 0$ there exists a component $U$ of $Y-N_r A$ with the following property: $U$ contains each
$R(r)$ deep point of $Y-N_r A$; and $U$ contains $S$-deep points for any $S \ge R(r)$.
\end{lemma}

\subparagraph{Remark:} Lemma~\ref{LemmaCodim2Fatness} and its codimension~1 analog
Lemma~\ref{LemmaCoarseJordan} each have a qualitative statement followed by a more quantitative statement.
These will be applied in the proof of Proposition~\ref{PropThickCrossing}, and in particular the more
quantitative statements are applied at one crucial point in Case~2 near the end of the proof of
Proposition~\ref{PropThickCrossing}. 

\bigskip

Before proving Lemma~\ref{LemmaCodim2Fatness}, we apply it to:

\begin{proof}[Proof of Proposition \ref{PropCodim2Thick}] Let $X \to X'$, $v \in T$, $e \subset T'$ be as in
the statement of the proposition. Suppose that $f \from X \to X'$ is a quasi-isometry. If $f(X_v) \csubset{c}
X'_e$, then some subset of $X'_e$ is uniformly properly equivalent to $X_v$ itself, contradicting
Lemma~\ref{LemmaCodim2Fatness}. If $f(X_v)$ crosses $X'_e$, then the image under a coarse inverse
$\bar f \from X' \to X$ of some subset of $X'_e$ coarsely separatex $X_v$, but this also violates
Lemma~\ref{LemmaCodim2Fatness}. This proves that $X_v$ is $e$-thick.
\end{proof}

\begin{proof}[Proof of Lemma \ref{LemmaCodim2Fatness}] The flavor of the proof is best appreciated in the
classical setting of Alexander duality, namely when $Y$ is a topological $n$-manifold. In this setting, since
$Y$ is contractible by hypothesis, Alexander duality holds in $Y$. Later, in the general setting, we will apply
the coarse Alexander duality theorem of \cite{KapovichKleiner:duality}.

Assume, then, that $Y$ is a topological $n$-manifold. Fix $r \ge 0$. There exists $s \ge 0$ and a CW
map $N_r A \to N_s B$ whose restriction to $A$ is a bounded distance from $f$; we still denote this map
$f\from N_ rA \to N_s B$. The constant $s$ depends only on $r$, the uniform contractibility data of $Z$,
and the uniform properness data of the original map~$f$. There exists $r'>r$ and a CW map $\bar f \from
N_s B \to N_{r'} A$ such that $\bar f \composed f$ is at distance $\le r'$ from the inclusion map
$N_r A \inject N_{r'} A$. The difference $r'-r$ depends only on $s$, the uniform contractibility data of
$X$, and the uniform properness data of $f$. Increasing~$r'$ by an amount that depends only on the given
$r'$ and the uniform contractibility data of~$X$, the map $\bar f \composed f \from N_r A \to N_{r'} A$
is homotopic to the injection $N_r A \inject N_{r'} A$. 

We shall show that $Y-N_r A$ has at most one component that contains all of the $R$-deep points, where
$R=r'-r$. The proof is a diagram chase using Alexander Duality.

The restriction homomorphism $i^*\from H^{n-1}_c(N_{r'} A)\to H^{n-1}_c(N_r A)$ factors through
$H^{n-1}_c(N_s B)$, as is shown in the triangle at the left of the following diagram:
$$
\xymatrix{
  0=H^{n-1}_c(N_s B) \ar[dr]_{f^*} 
& H^{n-1}_c(N_{r'} A) 
   \ar[d]_(.4){\circlearrowleft \quad i^*}
   \ar[l]_(.4){\bar f^*} 
& \wt H_0(Y-N_{r'} A) \ar[d]^{i_*} \ar@2{~}[l]_{AD} 
\\
&  H^{n-1}_c( N_r A)   
&  \wt H_0(Y -  N_r A) \ar@2{~}[l]^{AD} \ar@{}[ul]|{\circlearrowleft} 
}$$ 
The maps marked $AD$ are Alexander duality isomorphisms, and the right square is
commutative by naturality of Alexander duality. Since $Z$ has dimension $\leq n-2$, so does
its subcomplex $N_s B$, and so $H^{n-1}_c(N_s B)=0$ implying that $i^*$ is the zero map.   It follows
that $i_*\from \tilde H_0(Y - N_{r'} A ) \to \tilde H_0(Y-N_r A)$ is the zero map. In other words,
there is exactly one component of $Y-N_r A$ containing the entirety of $Y-N_{r'} A$. This shows
that $Y-N_r A$ has exactly one component containing all points that are $R=r'-r$ deep in $Y-N_r A$.

For the general proof, we drop the assumption that $Y$ is an $n$-manifold, and we replace ordinary
Alexander duality by:

\begin{theorem}[\cite{KapovichKleiner:duality} Coarse Alexander duality] 
\label{TheoremCoarseAlexDuality}
For every coarse $\PD(n)$ space~$Y$ there exists a constant $D \ge 0$ such that for every subcomplex $A
\subset Y$, every $R \ge D$, and every $* \in \Z$, there exists a homomorphism 
$$H^*_c(N_{R} A) \xrightarrow{AD_R} 
H_{n-*}\left(Y,\overline{Y-N_{R-D} A}\right) \approx \wt H_{n-*-1}\left(\overline{Y-N_{R-D} A}\right)
$$ 
where the isomorphism on the right exists whenever it makes sense (that is, when $* \ge 1$), and the
following hold for every $R \ge 2D$:
\begin{itemize}
\item[(1)] In the following diagram, the left square always commutes, and the right square
commutes whenever it makes sense:
\end{itemize}
$$
\xymatrix@C+25pt{ 
H^*_c(N_{R+D} A) \ar[r]^-{AD_{R+D}} \ar[d]^{i^*} & 
     H_{n-*}\left(Y , \overline{Y-N_R A} \right)  \ar[r]^{\approx} \ar[d]^{i_*} &
     \wt H_{n-*-1}\left( \overline{Y-N_R A} \right)  \ar[d]^{i_*}
\\ 
H^*_c(N_{R-D} A)  \ar[r]^-{AD_{R-D}} &
     H_{n-*} \left(Y , \overline{Y-N_{R-2D} A} \right) \ar[r]^{\approx} &
     \wt H_{n-*-1} \left( \overline{Y - N_{R-2D} A} \right) }
$$ 
\begin{itemize}
\item[(2)] We have
\begin{align*}
\ker(AD_{R+D}) &\subset \ker(i^*) \\
\image(i_*)   &\subset \image(AD_{R-D})
\end{align*}
\end{itemize}
\qed\end{theorem}

We now continue with the proof of Lemma~\ref{LemmaCodim2Fatness}. Assuming that $A$ is uniformly equivalent to
a subcomplex $B$ of the $n-2$ dimensional complex $Z$, it follows by the same argument as before that for each
$R \ge 2D$ there exists $r' \ge R-D$, with the difference $r'-(R-D)$ depending only on $R$, $D$, $Y$, $Z$, and
the uniform properness gauge between $A$ and $B$, such that the restriction homomorphism $i^* \from
H^{n-1}_c(N_{r'} A) \to H^{n-1}_c(N_{R-D} A)$ is the zero map. Enlarging $r'$ by at most $2D$ we get an odd
integer $k$ so that restriction homomorphism $i^* \from H^{n-1}_c(N_{R+kD} A) \to H^{n-1}_c(N_{R-D} A)$ is the
zero map. A diagram chase using coarse Alexander duality shows that the homomorphism $i_* \from \wt
H_0\left(\overline{Y-N_{R+(k+1)D} A}\right) \to \wt H_0\left(\overline{Y-N_{R-2D} A}\right)$ is the zero map.
It follows as before that $Y - N_{R-2D} A$ has at most one component containing all points that are $(k+3)D$
deep in $Y - N_{R-2D} A$. 
\end{proof}

\paragraph{The codimension~$1$ case.} In the codimension~$\ge 2$ case, vertex rigidity holds because no subset
of an edge space can coarsely separate a vertex space, a property which is much stronger than thickness. We now
turn to a case where the former property does not hold --- codimension one subgroups of $\PD$ groups. 

Here is another special case of depth zero vertex rigidity, which can be applied to prove depth zero vertex
rigidity for the classes of examples considered in Theorems~\ref{Filling}, \ref{Fibered},
and~\ref{TheoremBabyAbelian}. 

\begin{corollary} 
\label{CorollaryCodim1}
If $\G$, $\G'$ are graphs of groups such that every depth zero vertex group is a coarse $\PD(n)$ group,
every edge group is a coarse $\PD(n-1)$ group, and the crossing graph of each vertex is connected, then
depth zero vertex rigidity holds for $\G$ and $\G'$.
\qed\end{corollary}

This is an immediate application of Corollary~\ref{GDZVR} combined with the following:

\begin{proposition}\label{PropThickCrossing} Let $\G$ and $\G'$ be admissible graphs of groups whose depth
zero vertex groups are coarse $\PD$, and whose remaining vertex and edge groups are coarse finite type. Let $X
\to T$, $X'\to T'$ be the Bass-Serre trees of spaces for $\G,\G'$ respectively. Suppose that each depth
zero vertex space is coarse $\PD$, and each depth zero vertex which is a raft satisfies the crossing graph
condition. Let $v$ be a depth zero vertex of $T$ which is a raft, suppose that $X_v$ is coarse $\PD(n)$, and
suppose that the crossing graph $\epsilon_v$ is nonempty. Let $e$ be an edge of $\G'$ for which $\Stab(e)$ is
coarse $\PD(n-1)$. Then $X_v$ is $e$-thick.
\end{proposition}

For the proof we need two coarse Alexander duality arguments. Kapovich and Kleiner showed that if $X \subset Y$
where $Y$ is coarse $\PD(n)$ and $X$ is coarse $\PD(k)$, and if $X$ is not coarsely equivalent to $Y$, then $k<n$
\cite{KapovichKleiner:duality}. In group language this says that if $H$ is an infinite index subgroup of $G$
where $G$ is coarse $\PD(n)$ and $H$ is coarse $\PD(k)$, then $k<n$. The following lemma gives a generalization
which will be used several times henceforth:

\begin{lemma} 
\label{LemmaBigInSmall} 
If $Y$ is a coarse $\PD(m)$ space, $X$ is a uniformly contractible, bounded geometry
space of coarse dimension $k$, and $f \from X \to Y$ is a uniformly proper embedding, then $k \le
m$, and if $k = m$ then $X$ is coarse $\PD(m)$ and $f(X) \ceq{c} Y$.
\end{lemma}

\begin{proof} If $k>m$ then, applying coarse Alexander duality~\ref{TheoremCoarseAlexDuality}, for some
$r\ge 0$ we would have $H_{m-k}(Y;N_r(f(X)) \ne 0$, an absurdity.

If $k=m$ then it suffices to prove that $f(X) \ceq{c} Y$, because that implies that $f$ is a quasi-isometry and
so $X$ is coarse $\PD(m)$. Assuming that $A = f(X) \not\ceq{c} Y$, it follows that $H_0(Y,\overline{Y - N_R A})
= 0$ for all $R$. Applying coarse Alexander duality, the restriction maps $i^* \from H^*_c(N_{(2k+2)D} A) \to
H^*_c(N_{2kD} A)$ are zero for $* \ge m$ and all $k$. For each $k$ choose a map $\bar f_{2kD} \from N_{2kD} A
\to X$ that is a coarse inverse to $X \xrightarrow{f} A \inject N_{2k} A$. For each $k$ there exists $k' > k$ so
that the map $N_{2kD} A \xrightarrow{\bar f_{2kD}} X$ is boundedly homotopic to the map $N_{2kD} A \inject
N_{2k'D} A\xrightarrow{\bar f_{2k'D}} X$. Also, the map $X \xrightarrow{f} A \inject N_{2k A} \xrightarrow{\bar
f_{2kD}} X$ is boundedly homotopic to the identity, because $X$ is uniformly contractible. We thus have a
commutative diagram
$$
\xymatrix{
H^*_c(X)  \ar[dr]^(.7){\bar f_{2kD}^*} \ar[drr]^{\bar f_{2k'D}^*} \\
H^*_c(A) \ar[u]^{f^*} & H^*_c(N_{2kD} A) \ar[l] & H^*_c(N_{2k'D} A) \ar[l]
}$$
and any path in this diagram from $H^*_c X$ back to itself is the identity map. One such path factors
through the restriction map $H^*_c(N_{2k'D} A) \to H^*_c(N_{2kD} A)$. For $* \ge m$ the latter is the zero map,
and so $H^*_c X = 0$. But this contradicts that $X$ has coarse dimension~$m$, and so $f(X) \ceq{c} Y$. 
\end{proof}

Next we need to understand separation properties of codimension one subsets:

\begin{lemma}
\label{LemmaCoarseJordan}
If $Y$ is a uniformly contractible coarse $\PD(n)$ space, $Z$ is a uniformly contractible coarse
$\PD(n-1)$ space, and $f\from B \to Y$ is a uniformly proper embedding of a subset $B \subset Z$ to $Y$,
then $f(B)$ coarsely separates $Y$ if and only if $B\ceq{c} Z$. Moreover, if $B \ceq{c} Z$ then there are
precisely two deep components of $Y - N_r f(B)$ for each sufficiently large~$r$.

For a more quantitative version of the ``only if'' direction: for each $r$ there exists $R$ such that if
$B \not\ceq{c} Z$ then all $R$-deep points of $Y - N_r f(B)$ are contained in a single component of $Y -
N_r f(B)$; the constant $R$ depends only on $Y$, $Z$, $r$, and the properness gauge of $f$.
\end{lemma}

This is similar \cite{KapovichKleiner:duality}, Corollary 7.8, except we do not assume that $f \from B
\to Y$ is the restriction of a uniformly proper map $Z \to Y$.

\begin{proof} If $B \ceq{c} Z$ then  $f$ extends to a uniformly proper map $f\from Z\to Y$, and $f(Z)$
coarsely separates $Y$ if and only if $f(B)$ coarsely separates $Y$. The result now follows from
\cite{KapovichKleiner:duality} Corollary~7.8(1), including the desired conclusion about two deep
components.

If $B\not\ceq{c} Z$ then, setting $A = f(B)$, the argument proceeds exactly as in the proof of Lemma
\ref{LemmaCodim2Fatness} above, as long as we show that for each $R>0$ there exists $R' \ge R$ such
that the restriction homomorphism $H^{n-1}_c(N_{R'} A) \to H^{n-1}_c(N_R A)$ is zero, and such that $R' -
R$ depends only on $Y$, $Z$, the uniform properness gauge of $f$, and the number $R$. We may inductively
define sequences $0=S_0 < S_1 < \ldots$ and $0=R_0 < R_1 < \ldots$ and uniformly proper maps $f_i \from
N_{S_i} B \to N_{R_i} A$ and $\bar f_i \from N_{R_i} A \to N_{S_{i+1}} B$ so that $\bar f_i \composed
f_i$ and $f_{i+1} \composed \bar f_i$ are each homotopic to the corresponding inclusion. For example,
given $R_i$ and $f_i$, first one chooses $S_{i+1}$ and $\bar f_i$ so that $\bar f_i\composed f_i$ is a
bounded distance from inclusion, then one increases $S_{i+1}$ so that $\bar f_i\composed f_i$ is
homotopic to inclusion, each time using uniform contractibility of $Y$; the difference $S_{i+1}-S_i$
depends only on $Y$, $Z$, the uniform properness gauge of $f$, and the the previously defined stuff in
this inductive definition. It follows that the restriction homomorphism $H^{n-1}_c(N_{R_{i+k}} A) \to
H^{n-1}_c(N_{R_{i-1}} A)$ factors through the restriction homomorphism $H^{n-1}_c(N_{S_{i+k}} B) \to
H^{n-1}_c(N_{S_i} B)$. It therefore suffices to prove that for each $S>0$ there exists $S'>S$ such that
the restriction homomorphism $H^{n-1}_c(N_{S'} B) \to H^{n-1}_c(N_S B)$ is zero, and so that the
difference $S'-S$ depends only on $Z$. By applying Coarse Alexander Duality and a diagram chase, it
suffices to prove that $H_0(Z,\overline{Z - N_S B})$ is zero for all $S$. But this is immediate because
$Z$ is connected and $N_S B \subsetneq Z$.
\end{proof}

\begin{proof}[Proof of Proposition \ref{PropThickCrossing}] Since $X_v$ is $\PD(n)$ and the edge space
$X'_e$ of $e$ is $\PD(n-1)$, $X_v$ does not uniformly embed in $X'_e$, by Lemma~\ref{LemmaBigInSmall}.

Thus we need to see that $X_v$ cannot cross $e$. Assume to the contrary that $f \from X \to X'$ is a
quasi-isometry such that $f(X_v)$ crosses $X'_e$. It follows that the image under a coarse inverse $\bar f
\from X' \to X$ of some subset of $X'_e$ coarsely separates $X_v$. Lemma~\ref{LemmaCoarseJordan} now applies to
show that $\bar f(X'_e)$ must be coarsely contained in $X_v$, and so the map $\bar f$ may be moved a bounded
amount so that $\bar f (X'_e) \subset X_v$. The desired contradiction is obtained by quoting
Lemma~\ref{LemmaIntrinsicCrossing} below.
\end{proof}

\begin{lemma} 
\label{LemmaIntrinsicCrossing} 
Let $v$ be a depth zero vertex of $T$ which is a one vertex raft, suppose that $X_v$ is coarse $\PD(n)$, and suppose
that the crossing graph $\epsilon_v$ is nonempty. If $Z\subset X_v$ is any uniformly properly embedded, coarse
$\PD(n-1)$ space, then $X_v$ does not cross $Z$ in~$X$.
\end{lemma}

\begin{proof} \textbf{Step 1:} We show that for any edge $e'$ of $T$ incident to $v$, the edge space
$X_{e'}$ does not cross $Z$ in $X$. In proving this, we may assume that $X_{e'}$ is not coarsely
contained in $Z$. 

Let $w$ be the vertex of $e'$ opposite $v$. By admissibility, there is a depth zero vertex space
$X_{w'}$ so that $X_{e'} \csubset{c} X_w \csubset{c} X_{w'}$. If $X_{e'}\not\ceq{c} X_w$ then $e'$
separates $v$ from $w'$. But even if $X_{e'} \ceq{c} X_w$, we can still choose $w'$ so that it is
separated from $v$ by $e'$. To do this, note that $\Stab(e')$ has finite index in $\Stab(w)$. This index
must be greater than one by irreducibility, so there is a translate of $e'$ by $\Stab(w)$ which is an edge
leading from $w$ to a translate $w'$ of $v$. 

We claim that for any $r,r'$ the subset $X_{w'} \intersect \left( N_r Z \intersect N_{r'} X_{e'}\right)$ does
not coarsely separate $X_{w'}$. More quantitatively, the set $X_{w'} - \left( X_{w'} \intersect N_r Z\intersect
N_{r'} X_{e'} \right)$ has at most one $R$ deep component in $X_{w'}$; the constant $R$ depends only on
$r,r'$, and the quasi-isometry constants of $f$.

To prove the claim, let $e''$ be the last edge on the edge path from $e'$ to $w'$, and let $X_{w'}$ be coarse
$\PD(k)$. We may assume that $e''$ has positive depth, or we would take a shorter path. 

In the case that $X_{e''}$ is not coarse $\PD(k-1)$ then, by Lemma~\ref{LemmaCodim2Fatness}, no subset of
$X_{w'}$ which is coarsely contained in $X_{e''}$, and in particular no subset of $X_{e'}$, can coarsely
separate $X_{w'}$. This holds quantitatively as well, proving the claim in this case.

In the case that $X_{e''}$ is coarse $\PD(k-1)$ then, by Lemma~\ref{LemmaCoarseJordan}, no subset of $X_{w'}$
which is \emph{strictly} coarsely contained in $X_{e''}$ can coarsely separate $X_{w'}$. By our assumption that
$X_{e'}$ is not coarsely contained in $Z$, for any $r,r' \ge 0$ the set $N_r Z \intersect N_{r'} X_{e'}$ is
strictly coarsely contained in $X_{e'}$, and so again the subset $X_{w'} \intersect \left( N_r Z \intersect
N_{r'} X_{e'} \right)$ does not coarsely separate $X_{w'}$. Again the claim follows by making this more
quantitative. 

Choose $r'$, depending only on $X$, so that any point in $X_{e'}$ is connected to a point in $X_{w'}$ by a
path of length $\le r'$. Fix $r$. Given $s$, consider points $x_-,x_+ \in X_{e'}$ at distance $\ge s$ from
$N_r Z$. Connect $x_-,x_+$ to $y_-,y_+\in X_{w'}$ by paths of length $r'$. By choosing $s$ sufficiently
large, these paths are in the complement of $N_r Z$, and the points $y_-,y_+$ are $R$ deep in $X_{w'} -
\left(N_r Z\intersect N_{r'} X_{e'}\right)$. The choice of $s$ depends only on $r'$, $r$, $R$, and the
uniform properness gauge for the embedding $Z \inject X_v$. It follows that $y_-,y_+$ are in the same
component of $X_{w'} - \left(N_r Z \intersect N_{r'} X_{e'}\right)$, and so $x_-,x_+$ are in the same
component of $X - N_r Z$. This proves that $X_{e'}$ does not cross $Z$ in $X$; more quantitatively, all
points in $X_{e'}$ that are $s$ deep in $X - N_r Z$ are contained in a single component of $X - N_r Z$.

\subparagraph{Step 2:} Suppose now that $X_v$ crosses $Z$. Fix a sufficiently large $r \ge 0$, and let $X - N_r
Z = \U^- \union\U^+$ be an open decomposition of $X-N_r Z$, such that each of $\U^-,\U^+$ is deep in~$X$, and
$X_v$ has deep intersection with each of $\U^-,\U^+$. By Lemma~\ref{LemmaCoarseJordan}, exactly two of the
components of $X_v - N_r Z$ are deep in $X_v$; one of these components, denoted $\U^-_v$, must
be contained in $\U^-$ while the other, denoted $\U^+_v$, is contained in $\U^+$. Any set $A \subset X_v$
that crosses $Z$ in $X_v$ also crosses $Z$ in $X$, because $A$ has deep intersection with both of
$\U^-_v$ and $\U^+_v$, implying that $A$ has deep intersection with both of $\U^-,\U^+$.

From what we have proved so far, each edge space $X_{e'}$ incident to $X_v$ does not cross $Z$ in $X$, and
so either $X_{e'} \csubset{c} Z$ or $X_{e'}$ has deep intersection with exactly one of $\U^-,\U^+$. More
quantitatively, there exists some $R \ge 0$ such that if $X_{e'}\not\csubset{c} Z$ then $X_{e'}$ has deep
intersection with exactly one of $\U^-,\U^+$ while $X_{e'}$ does not even have $R$ deep intersection with
the other. 

Now we focus on those edge spaces $X_{e'}$ incident to $X_v$ that are coarse $\PD(n-1)$ spaces; these are
the vertices of the crossing graph $\epsilon_v$, and we write $X_{e'} \in \Vertices(\epsilon_v)$. By
moving a uniformly finite Hausdorff distance, we may regard $X_{e'}$ as a subset of $X_v$.

\paragraph{Case 1:} Suppose that there exists $X_{e'} \in \Vertices(\epsilon_v)$ such that $X_{e'}
\csubset{c} Z$. It follows that $X_{e'} \ceq{c} Z$, by Lemma~\ref{LemmaBigInSmall}, together with the fact that
$X_{e'}$ and $Z$ are both coarse $\PD(n-1)$. The crossing graph $\epsilon_v$ has at least one other vertex. By
connectivity of $\epsilon_v$, there exists an edge space $X_{e''}$ incident to $X_v$ that crosses $X_{e'}
\ceq{c} Z$ in $X_v$, and so $X_{e''}$ also crosses $Z$ in $X$, contradicting Step~1.

\paragraph{Case 2:} 
Suppose that for each $X_{e'} \in \Vertices(\epsilon_v)$, $X_{e'}$ has deep
intersection with exactly one of $\U^-,\U^+$, and so $X_{e'}$ may be labelled either positive or negative,
depending on whether $X_{e'}$ has deep intersection with $\U^+$ or with $\U^-$. Since the action of
$\Stab(X_v)$ on $X_v$ is cocompact, the set $\U^+$ contains an $R$-deep point that lies on some
$X_{e'}\in\epsilon_v$. 
But this implies that $X_{e'}$ intersects $\U^+$ deeply (this is the most
crucial place where we make use of the quantitative versions of our various propositions and lemmas).
This shows that $\epsilon_v$ contains a positive vertex, and a similar argument with $\U^-$ shows that
$\epsilon_v$ contains a negative vertex. By connectivity of $\epsilon_v$, it follows that there exists a
positive $X_{e_1}\in\Vertices(\epsilon_v)$ and a negative $X_{e_2} \in\Vertices(\epsilon_v)$ such that
$X_{e_1}$ is connected to $X_{e_2}$ by an edge of $\epsilon_v$. Clearly $X_{e_1}$ and $X_{e_2}$ do not
cross each other in $X_v$, and so by definition of the edges of $\epsilon_v$ it follows that there is an
edge space $X_{e_3}$ incident to $X_v$ which crosses both $X_{e_1}$ and $X_{e_2}$ in $X_v$. But this
implies that $X_{e_3}$ crosses $Z$ in $X_v$, which implies that $X_{e_3}$ crosses $Z$ in $X$, a
contradiction.

This completes the proof of Lemma~\ref{LemmaIntrinsicCrossing}.
\end{proof}

The next corollary is a consequence of our various crossing arguments above, and will be useful in what follows.

\begin{corollary} 
\label{CorStrongNonCrossing}
Let $v$ be depth zero vertex with coarse $\PD(n)$ stabiliser, and let $e$ be a positive depth edge incident to $v$. Denote the
two components of $X-X_e$ as $\V_-$, $\V_+$, so that $X_v \subset \V_-$. It follows that $X_v$ crosses no subset of $\V_+$.
\end{corollary}

\begin{proof} Let $A \subset \V^+$, and consider $r \ge 0$. There exists $s$ so that $N_r(A) \subset N_s(\V_+) = N_s(X_e) \union
\V_+$.

We break into several cases, depending on the nature of the depth zero raft containing~$v$.

First consider the case that $v$ is a one vertex raft and $X_e$ coarsely separates $X_v$. From the crossing graph condition it
follows that $X_e$ is coarse $\PD(n-1)$ and the crossing graph of $X_v$ is nonempty.
Applying Proposition~\ref{PropThickCrossing}, let $W$ be the unique deep component of $X-N_s(X_e)$ such that $X_v$ intersects
$W$ deeply. Note that there is a uniform upper bound to the distance from a point of $X_v \intersect N_s(X_e)$ to the set $W$,
in other words, the set $X_v \intersect N_s(X_e)$ is not deep in $X_v$. The connected set $W$ is contained in a unique component
$W'$ of $X - N_r(A)$, and it follows that $W'$ is the only component of $X-N_r(A)$ that $X_v$ intersects deeply.

Next consider the case that $v$ is a one vertex raft and $X_e$ does not coarsely separate $X_v$. Applying
Lemma~\ref{LemmaCodim2Fatness}, there is a unique deep component $W$ of $X-N_s(X_e)$ such that $X_v$ intersect
$W$ deeply, and now we continue as before.

Finally consider the case that $v \in R$ for some bushy raft $R$. For any other vertex $w$ of $R$, $X_v \ceq{c}
X_w$ and so $X_v$ crosses $A$ if and only if $X_w$ crosses $A$. But $R$ has vertices $w$ that are arbitrarily
far from $v$, and so we can choose $w$ so that $x_v \intersect N_s(X_e) = \emptyset$. It follows that $X -
N_s(X_e)$ contains a unique deep component
$W$ that contains $X_w$, and so $X - N_r(X_A)$ has a unique deep component $W'$ that contains $X_w$, from which it follows that
$X_v \csubset{c} W'$, that is, $W'$ is the only component of $X - N_r(A)$ that $X_v$ intersects deeply.
\end{proof}

\paragraph{Proof of the Depth Zero Vertex Rigidity Theorem \ref{TheoremVertexRigidity}.} We now have all
the ingredients for the proof. To review the notation, $\G,\G'$ are finite type,
irreducible, admissible graphs of groups, with Bass-Serre trees of spaces $\pi \from X \to T$, $\pi'
\from X' \to T'$ each satisfying the properties that no depth zero raft is a line, each depth zero vertex
space is coarse $\PD$, every other edge and vertex space is of coarse finite type, and each depth zero vertex
which is a raft satisfies the crossing graph condition.

The coarse finite type condition implies that the coarse dimension of each edge and vertex space
is well-defined. Lemmas~\ref{LemmaCodim2Fatness} and~\ref{LemmaCoarseJordan} have the following immediate
consequence:

\begin{proposition} 
\label{PropEdgeSeparatesVertex}
If $v$ is a depth zero vertex such that $v$ is a raft and $X_v$ is coarse $\PD(n)$, then for each edge $e$
incident to $v$, $X_e$ has a subset $A$ that coarsely separates $X_v$ if and only if $X_e$ is coarse
$\PD(n-1)$ and $A \ceq{c} X_e$.
\qed\end{proposition}

The conclusion of this proposition can be used to replace the coarse finite type condition, giving
a broader vertex rigidity theorem, leading in turn to broader QI-rigidity and classification theorems, in which the
coarse finite type condition is replaced by the conclusion of Proposition~\ref{PropEdgeSeparatesVertex}. This
would extend the ``bottleneck'' of depth zero vertex rigidity a little bit, by allowing certain edge groups
that are not of coarse finite type, such as non-finitely presentable and non-coarsely separating subgroups of
aspherical manifold groups. We will not pursue this issue further here.

\bigskip 

In principle we prove Theorem~\ref{TheoremVertexRigidity} by proving that $X$ and $X'$ are mutually thick, and then applying
Theorem~\ref{AbstractRaftRigidity} and Corollary~\ref{GDZVR}. In reality the logic of our proof will flow a little differently,
requiring us to apply the proofs of these results, rather than the results themselves. The way we really proceed is by
establishing depth zero raft rigidity dimension by dimension, proceeding by downward induction on the coarse dimension of raft
spaces. 

Recall the notation $\V_0(T)$ for the depth zero vertices of $T$. Given $v \in \V_0(T)$ let $R(v) \subset T$ be the depth zero
raft containing $v$. Given $n \ge 0$, let
$$\V_0^n(T) = \{v \in \V_0(T) \suchthat \,\,\text{$X_{R(v)}$ has coarse dimension $\ge n$.}\}
$$
and similarly for $T'$. To put it another way, $\V_0^n(T)$ is the set of depth zero vertices $v$ such that \emph{either} $v$ is
a one vertex raft and the coarse dimension of $X_v$ is $\ge n$, \emph{or} $v$ is contained in a bushy raft and the coarse
dimension of $X_v$ is $\ge n-1$.

We shall prove the following statement by downward induction on $n$:
\begin{description}
\item[$\textbf{VR}_n$:] For any quasi-isometry $f \from X \to X'$ there exists a quasi-isometry $f_\# \from
\V_0^n(T) \to \V_0^n(T')$ such that if $v\in \V_0^n(T)$ then $f(X_v) \ceq{[C]} X'_{f_\#(v)}$, and if $v' \in
\V_0^n(T')$ then there exists $v \in \V_0^n(T)$ such that $f(X_v) \ceq{[C]} X'_{v'}$. The constant $C$ and the
quasi-isometry constants of $f_\#$ depend only on $\G$, $\G'$, and the quasi-isometry constants of $f$. 
\end{description}
The statement $\textbf{VR}_0$ is the desired conclusion of Theorem~\ref{TheoremVertexRigidity}. 
The statement $\textbf{VR}_n$ is vacuously true for sufficiently large $n$. Assuming by induction that
$\textbf{VR}_{n+1}$ is true, we will prove $\textbf{VR}_{n}$.

There are two major steps to the proof of $\textbf{VR}_n$:
\begin{description}
\item[Step 1: Thickness.] For any vertex $v \in \V_0^n(T) - \V_0^{n+1}(T)$, and for any positive depth edge $e$ of $T'$, $X_v$ is
$e$-thick.
\end{description}
Now fix a quasi-isometry $f \from X \to X'$. For any depth zero raft $R \subset T$ such that $X_R$ has coarse dimension $n$, we
can combine Step~1 and Lemma~\ref{LemmaLocalRaftRigidity} to obtain a depth zero raft $R' \subset T'$ such that $f(X_R)
\csubset{c} X'_{R'}$, with coarse inclusion constant depending only on $X$, $X'$, and the quasi-isometry constants of $f$. The
coarse dimension of $X'_{R'}$ cannot be $>n$, because if it were then the induction hypothesis $\textbf{VR}_{n+1}$ would
produce a depth zero raft $R_1\subset T$ with $X_{R_1}$ of coarse dimension $>n$ such that $X_R \csubset{c} X_{R_1}$, a
contradiction. It follows that the coarse dimension of $X'_{R'}$ is $\le n$.
\begin{description}
\item[Step 2: Dimension preservation.] The coarse dimension of $X'_{R'}$ equals~$n$.
\end{description}
Once these two steps are established, the proof of $\textbf{VR}_n$ is finished as follows. Steps~1 and~2 allow us to define a
map $f_{\#\#}$ from the set of dimension~$n$ depth zero rafts of $T$ to the set of dimension~$n$ depth zero rafts of $T'$, so
that if $f_{\#\#}(R) = R'$ then $f(X_R)\csubset{c} X'_{R'}$. Following the proof of Theorem~\ref{AbstractRaftRigidity}, the map
$f_{\#\#}$ is a bijection: the same argument applied to a coarse inverse $\bar f \from X' \to X$ of $f$ gives a map
$\bar f_{\#\#}$ from dimension~$n$ depth zero rafts of $T'$ to dimension~$n$ depth zero rafts of $T$; the fact that no depth
zero raft is contained in a finite neighborhood of another shows that $f_{\#\#}$ and $\bar f_{\#\#}$ are inverse bijections. 
Finally, we apply the proof of Corollary~\ref{GDZVR}: the results of \cite{MSW:QTOne} and the hypothesis of coarse $\PD$ vertex
groups in depth zero combine to show that if $f_{\#\#}(R)=R'$ then for any vertex $v \in R$, $f(X_v)\ceq{c} X'_{v'}$ for some
$v' \in R'$, with coarse inclusion constant depending only on $X$, $X'$, and the quasi-isometry constants of $f$, allowing us
to define $f_\#(v)=v'$. The same arguments as in Corollary~\ref{GDZVR} prove that the map $f_\#$ is a quasi-isometry.

It remains to establish Steps~1 and~2.

\paragraph{Proof of Step 1: Thickness.} Fix a quasi-isometry $f \from X \to X'$, a vertex $v \in \V_0^n(T) - \V_0^{n+1}(T)$,
and a positive depth edge $e$ of $T'$. We must prove:
\begin{description}
\item[(a)] $f(X_v)$ does not cross $X'_e$.
\item[(b)] $f(X_v) \not\csubset{c} X'_e$.
\end{description}
Let $R_0$ be the depth zero raft of $T$ containing $v$. We break the proof into several cases.

\subparagraph{Case 1:} Suppose that $X'_e \csubset{c} X'_{R'}$ for some depth zero raft $R' \subset T'$ such that $X'_{R'}$ has
coarse dimension $k \ge n+1$. By the induction hypothesis $\textbf{VR}_{n+1}$ there exists a depth zero raft $R \subset T$ such
that $X_R$ has coarse dimension $k$ and $f(X_R) \ceq{c} X'_{R'}$. 

To prove (a), suppose that $f(X_v)$ crosses $X'_e$. It follows that $X_v$ crosses a subset of $X_R$. On the path from $v$ to
$R$ in $T$, let $v'$ be the last vertex in $R_0$, and let $e'$ be the first edge not in $R_0$, so $X_v \ceq{c} X_{v'}$ and $e'$
has positive depth. Let $\V_+$ be the component of $X-X_{e'}$ disjoint from $X_{v'}$. Applying
Corollary~\ref{CorStrongNonCrossing} it follows that $X_{v'}$ does not cross any subset of $\V_+$, and so
$X_v \ceq{c} X_{v'}$ does not cross any subset of $X_R$, a contradiction.

To prove (b), if $f(X_v) \csubset{c} X'_e$ then $X_v \csubset{c} X_{R}$, but no depth zero vertex space is coarsely contained in
a disjoint depth zero raft space.

\subparagraph{Case 2:} Suppose that $X'_e \csubset{c} X'_{R'}$ for some depth zero raft $R' \subset T'$ such
that $X'_{R'}$ has coarse dimension $k \le n-1$. Since $e$ has positive depth, there is a unique depth zero
vertex $w \in R'$ closest to $e$, and $X'_e \cstrict{c} X'_w$, that is, $X'_e$ is uniformly properly equivalent
to a strictly proper subset of $X'_w$. Note that $X'_w$ is coarse $\PD(l)$ where $l=k$ if and only if $R' =
\{w\}$ is a one vertex raft, and $l=k-1$ if and only if $R'$ is bushy. In particular, $l \le k \le n-1$. Since
$f(X_v)$ crosses $X'_e$, it follows that some subset of $X'_e$ coarsely separates $f(X_v)$, and so some
strictly proper subset of the coarse $\PD(l)$ space $X'_w$ coarsely separates $f(X_v)$.

\subparagraph{Case 2a:} $R_0 = \{v\}$. It follows that $X_v$ is coarse $\PD(n)$, and we obtain a contradition to
Lemma~\ref{LemmaCoarseJordan} which says that since $l<n$, the coarse $\PD(n)$ space $f(X_v)$ cannot be
coarsely separated by a uniformly proper embedding of a subset of the coarse $\PD(l)$ space $X'_w$, proving
(a). Nor can $f(X_v)$ be coarsely contained in $X'_w$, proving (b).

\subparagraph{Case 2b:} $R_0$ is bushy and $R'$ is bushy. It follows that $X_v$ is coarse $\PD(n-1)$, but $X'_w$
is coarse $\PD(k-1)$, and again we prove (a) by contradicting Lemma~\ref{LemmaCoarseJordan}. The proof of (b) is
similar.

\subparagraph{Case 2c:} $R_0$ is bushy and $R'=\{w\}$. Let $e'$ be the last edge on the edge path from $e$ to $w$, so
we clearly have $X'_e \csubset{c} X'_{e'} \csubset{c} X'_w$. 

The space $X_v$ is coarse $\PD(n-1)$ and $X'_w$ is coarse $\PD(k)$. We may again apply Lemma~\ref{LemmaCoarseJordan}, although
this time instead of reaching a contradiction we only reach the conclusion that $k=n-1$. Since $f(X_v)$ crosses $X'_e$, there
is a subset $A \subset X'_e$ such that $A \csubset{c} f(X_v)$ and $A$ coarsely separates $f(X_v)$. This subset $A$ is also
coarsely contained in the coarse $\PD(n-1)$ space $X'_w$. By applying Lemma~\ref{LemmaSeparatingPair} below, it follows that
$A$ coarsely separates $X'_w$. This implies that some subset of $X'_{e'}$, coarsely equivalent to $A$, coarsely separates
$X'_w$. By the crossing graph condition and the fact that $w$ is a one vertex raft with coarse $\PD(n-1)$ stabilizer it follows
that $X'_{e'}$ is coarse $\PD(n-2)$ and $A\ceq{c} X'_{e'}$. This implies that $X'_{e'} \ceq{c} A \csubset{c} X'_e$, and so
$X'_e\ceq{c} X'_{e'} \ceq{c} A$. 

To recapitulate what we know up to this point, $X'_e$ is coarse $\PD(n-2)$, and $X'_e \csubset{c}
f(X_v)$. This implies that the coarse $\PD(n-1)$ vertex space $X_v$ crosses the coarse $\PD(n-2)$ space $\bar
f(X'_e)\csubset{c} X_v$. Now we use bushiness of $R_0$, in the same manner as was used in the proof of
Corollary~\ref{CorStrongNonCrossing}. Since $X_v$ crosses
$\bar f(X'_e)$, for any vertex $u \in R_0$ the vertex space $X_u \ceq{c} X_v$ also crosses $\bar f(X'_e)$. But we can use
bushiness of $R_0$ to choose $u$ sufficiently far from $v$ so that $X_u \intersect \bar f(X'_e) = \emptyset$, contradicting
that $X_u$ crosses $\bar f(X'_e)$. 

This completes the proof of (a) in Case 2c. Since $e$ is of positive depth in $T'$, if $f(X_v) \csubset{c} X'_{e}$ then the
coarse $\PD(n-1)$ space $f(X_v)$ is strictly coarsely contained in the coarse $\PD(n-1)$ space $X'_w$, a contradiction that
proves (b).

This completes the proof of Case 2.

\bigskip

In order to handle Case 2c we used the following:

\begin{lemma}
\label{LemmaSeparatingPair}
For any coarse $\PD(n)$ spaces $A,B$, and for any subcomplexes $C \subset A$ and $D \subset B$, if $C$ and $D$
are uniformly properly equivalent, then $C$ coarsely separates $A$ if and only if $D$ coarsely separates $B$.
\end{lemma}

\begin{proof}
Let $f \from C \to D$, $\bar f \from D \to C$ be coarsely inverse uniformly proper equivalences. By thickening
up the ranges of these maps, applying uniform contractibility of $A$ and $B$, and moving $f,f'$ a bounded
amount, we get continuous maps $f \from C \to N_s(D)$, $\bar f \from D \to N_s(C)$ for some $s \ge 0$.

Assume $C$ coarsely separates $A$. A coarse Alexander duality argument, similar to arguments we've used before, shows that for
some $R_0 < R_1 < R_2 < R_3 < \cdots$, in the sequence of restriction homomorphisms,
$$\cdots \xrightarrow{i^*_3} H^{n-1}_c(N_{R_3} C) \xrightarrow{i^*_2} H^{n-1}_c(N_{R_2} C) \xrightarrow{i^*_1} 
H^{n-1}_c(N_{R_1} C) \xrightarrow{i^*_1} H^{n-1}_c(N_{R_0} C)
$$
there exist nonzero elements $\omega_i \in H^{n-1}_c(N_{R_i} C)$ such that $i^*_{i-1}(\omega_i) = \omega_{i-1}$. 

By uniform contractibility of $A$ and $B$, we may rechoose the sequence $R_0 < R_1 < R_2 < \cdots$, and find a sequence $R'_0 <
R'_1 < R'_2 < \cdots$, with the following properties: for each $i$ there exists continuous maps
$$f_i \from N_{R_i} C \to N_{R'_{i+1}} D, \quad \bar f_i \from N_{R'_i} D \to N_{R_{i+1}} C
$$
at bounded distance from the maps $f, \bar f$, respectively, so that we have $f=f_0 \subset f_1 \subset f_2 \subset \cdots$,
and $f'=f'_0 \subset f'_1 \subset f'_2 \subset \cdots$, and so that $\bar f_{i+1} \composed f_i$ is homotopic to the inclusion
$N_{R_i} C \subset N_{R_{i+2}} C$, and $f_{i+1} \composed \bar f_i$ is homotopic to the inclusion $N_{R'_i} D \subset
N_{R'_{i+2}} D$. These maps, moreover, commute with the various restriction maps. Passing to $H^*_c$ we obtain a
commutative diagram
$$
\xymatrix{ 
H^{n-1}_c(N_{R_0} C) &
     H^{n-1}_c(N_{R_1} C) \ar[l]_{i^*_0} \ar[dl] &
     H^{n-1}_c(N_{R_2} C) \ar[l]_{i^*_1} \ar[dl] &
     H^{n-1}_c(N_{R_3} C) \ar[l]_{i^*_2} \ar[dl] &
     \cdots \ar[l] \\
H^{n-1}_c(N_{R'_0} D)  &
     H^{n-1}_c(N_{R'_1} D) \ar[l]_{j^*_0} \ar[ul] &
     H^{n-1}_c(N_{R'_2} D) \ar[l]_{j^*_1} \ar[ul] &
     H^{n-1}_c(N_{R'_3} D) \ar[l]_{j^*_2} \ar[ul] &
     \cdots \ar[l]
}
$$ 
Set $\omega'_i = \bar f_i^*(\omega_{i+1}) \in H^{n-1}_c(N_{R'_i} D)$. For $i \ge 1$ we chase through the
diagram to get
$$f_{i-1}^*(\omega'_{i}) = f_{i-1}^*(\bar f_{i}^*(\omega_{i+1})) = i^*_{i} \composed i^*_{i-1} (\omega_{i+1}) = \omega_{i-1}
$$
and so $\omega'_i \ne 0$. Commutativity of the diagram shows that 
$$j^*_{i-1}(\omega'_i) = \omega'_{i-1}
$$
A coarse Alexander duality argument now shows that $D$ coarsely separates $B$.
\end{proof}

\subparagraph{Case 3:} For every depth zero raft $R' \subset T'$ such that $X'_e \csubset{c} X'_{R'}$, the raft space $X'_{R'}$
has coarse dimension $n$. Let $w \in R'$ be the point closest to $e$, and let $e'$ be the last edge on the edge path from $e$ to
$w$ in $T'$. We have $X'_e \csubset{c} X'_{e'} \csubset {c} X'_w$. Since $f(X_v)$ crosses $X'_e$, it follows that $f(X_v)$ is
coarsely separated by a subset of $X'_{e'}$.

\subparagraph{Case 3a:} Suppose that $R_0 = \{v\}$ is a one vertex raft and $R'$ is bushy, so $X_v$ is coarse $\PD(n)$,
$X'_w$ is coarse $\PD(n-1)$, and $X_v$ is coarsely separated by a strictly proper subset of $X'_w$, a contradiction. This proves
(a) in case 3a.

To prove (b), if $f(X_v) \csubset{c} X'_e$ then the coarse $\PD(n)$ space $f(X_v)$ is coarsely contained in the coarse $\PD(n-1)$
space $X'_w$, a contradiction.

\subparagraph{Case 3b:} Suppose that $R_0 = \{v\}$ and $R' = \{w\}$ are one vertex rafts, so both $X_v$ and $X'_w$ are
coarse $\PD(n)$. Since $f(X_v)$ crosses $X'_e$, some subset $A \subset X'_e$ coarsely separates $f(X_v)$, and so some
subset $A' \subset X'_{e'}$, coarsely equivalent to $A$, coarsely separates $f(X_v)$. By
Lemma~\ref{LemmaSeparatingPair}, $A'$ coarsely separates $X'_w$. By the crossing graph condition, $X'_{e'}$ is coarse
$\PD(n-1)$ and $A' \ceq{c} X'_{e'}$. Thus, $X'_e \ceq{c} X'_{e'}$. It follows that $\bar f(X'_e) \csubset{c} X_v$, and
thus $X_v$ crosses a coarse $\PD(n-1)$ subset of itself, namely $\bar f(X'_e)$. By Lemma~\ref{LemmaIntrinsicCrossing}
it follows that the crossing graph of $v$ is empty. By the crossing graph condition, this implies that $\bar f(X'_w)$
does not cross any of the edge spaces incident to $X_v$. Local raft rigidity now implies that one of two possibilities
happens. The first possibility is that there is an edge $e_0$ incident to $X_v$ such that $\bar f(X'_w)$ is coarsely
contained in the component $\V_+$ of $X-X_{e_0}$ disjoint from $X_v$, and so $\bar f(X'_{e'}) \csubset{c} \V_+$, but
since $X_v$ crosses a subset of $\bar f(X'_{e'})$ we contradict Corollary~\ref{CorStrongNonCrossing}. The other
possibility is that $\bar f(X'_w)\csubset{c} X_v$, which implies $\bar f(X'_w) \ceq{c} X_v$, which implies that $X'_w$
crosses its incident $\PD(n-1)$ edge space $X'_{e'}$, which contradicts Lemma~\ref{LemmaIntrinsicCrossing}. This proves (a) in
case 3b.

To prove (b), if $f(X_v) \csubset{c} X'_e$, then since $X'_e \cstrict{c} X'_w$, it follows that the coarse $\PD(n)$ space $f(X_v)$
is strictly coarsely contained in the coarse $\PD(n)$ space $X'_w$, a contradiction.

\subparagraph{Case 3c:} Suppose that $R_0$ is bushy and $R'=\{w\}$ is a one vertex raft, so $X_v$ is coarse $\PD(n-1)$ and
$X'_{R'}$ is coarse $\PD(n)$. Cases~1, 2a, 3a and~3b apply to $\bar f$ and $R'$, with the conclusion that
$\bar f(X_w)$ does not cross any positive depth edge space of $X$. Local raft rigidity now implies that $\bar
f(X'_w)\csubset{c} X_{R_1}$ for some depth zero raft $R_1\subset T$ such that $X_{R_1}$ has coarse dimension
$\le n$. Thus, $X_{v}$ crosses some subset $X_{R_1}$. 

If $R_1 \ne R_0$ then applying Corollary~\ref{CorStrongNonCrossing} to the last vertex in $R_0$ on the path from $v$ to $R_1$, we
get a contradiction.

If $R_1 = R_0$ then we conclude that $\bar f(X'_w) \csubset{c} X_{R_0}$. Now we can apply \cite{FarbMosher:ABC} to conclude that
$\bar f(X'_w)\ceq{c} X_L$ for some line $L \subset R_0$. Thus, the vertex space $X_v$ crosses some subset $A \subset X_L$. For any
vertex $u\in R_0$, since $X_u \ceq{c} X_v$ it follows that $X_u$ crosses $A$. But since $R_0$ is bushy we can choose $u$
sufficiently far from $L$ that $X_u \intersect L = \emptyset$, contradicting that $X_u$ crosses $L$. This proves (a) in case 3c.

To prove (b), suppose that $f(X_v) \csubset{c} X'_e$. We also have $X'_e \csubset{c} X'_{e'} \csubset{c} X'_w$, and so
$X'_{e'}$ coarsely contains a set that coarsely separates the coarse $\PD(n)$ space $X'_w$, namely, the coarse $\PD(n-1)$ space
$f(X_v)$. By the crossing graph condition it follows that $X'_{e'}$ is coarse $\PD(n-1)$ and the crossing graph of $w$ is
connected. It also follows that $f(X_v) \ceq{c} X'_{e'}$. We therefore have $f(X_v) \ceq{c} X'_e \ceq{c} X'_{e'}$.

We still know by the argument above that $\bar f(X'_w)$ is coarsely contained in some depth zero raft space $X_{R_1}$. If $R_1 \ne
R_0$ then $X_v$ is coarsely contained in $\bar f(X'_e)$ which is coarsely contained in the raft space $X_{R_1}$, but no vertex
space in one depth zero raft space can be coarsely contained in a different depth zero raft space. If $R_1 = R_0$ then, as above,
$\bar f(X'_w) \ceq{c} X_L$ for some line $L \subset R_0$. Note that $X_L$ crosses $X_v$, because $X_L$ clearly crosses $X_e$ for
any edge $e \subset L$, but $X_e \ceq{c} X_v$. It follows that $X'_w$ crosses $f(X_v) \ceq{c} X'_{e'}$. But this contradicts
Corollary~\ref{CorStrongNonCrossing}.

\subparagraph{Case 3d:} Suppose that $R_0$ and $R'$ are both bushy, so $X_v$ and $X'_w$ are both coarse $\PD(n-1)$. Since $X_v$
crosses $\bar f(X'_e)$, and since $X_u \ceq{c} X_v$ for any other vertex $u \in R_0$, $X_u$ also crosses $\bar f(X'_e)$. In
particular, $\bar f(X'_e) \intersect X_u \ne \emptyset$ for all vertices $u \in R_0$. Letting
$L$ be a line in the raft $R_0$, it follows that $\bar f(X'_e)$ crosses the coarse $\PD(n)$ space $X_L$. Letting $A = X_L
\intersect \bar f(X'_e)$, the same coarse Alexander duality argument used in the proof of Lemma~\ref{LemmaSeparatingPair} above
shows that there is a sequence $R_0 < R_1 < R_2 < \cdots$ and nonzero elements $\omega_i \in H^{n-1}_c N_{R_i} A$ such that
$\omega_i$ restricts to $\omega_{i-1}$. Letting $B \subset X'_w$ be a subcomplex which is coarsely equivalent to $f(A)$ in $X'$,
we conclude that there is a sequence $R'_0 < R'_1 < R'_2 < \cdots$ and nonzero elements $\omega'_i \in H^{n-1}_c N_{R'_i} B$ such
that $\omega'_i$ restricts to $\omega'_{i-1}$. But $X'_w$ is coarse $\PD(n-1)$, and so a coarse Alexander duality argument shows
that $B \ceq{c} X'_w$. We therefore have $X'_w \ceq{c} B \ceq{c} B' \subset X'_e$, and also $X'_e \csubset{c} X'_w$, so $X'_e \ceq{c}
X'_w$, contradicting that $e$ has positive depth. This proves (a) in case 3d.

To prove (b), suppose $f(X_v) \csubset{c} X'_e$. We also have $X'_e \csubset{c} X'_w$. Since $f(X_v)$ and $X'_w$ are both coarse
$\PD(n-1)$ and $f(X_v) \csubset{c} X'_w$, it follows that $f(X_v) \ceq{c} X'_w$, and so $f(X_v) \ceq{c} X'_e \ceq{c} X'_w$,
contradicting that $e$ has postive depth.

This completes the proof of Step 1: Thickness.

\paragraph{Proof of Step 2: Dimension preservation.} First we review our knowledge at this stage. The induction hypothesis
$\textbf{VR}_{n+1}$ gives us that $f$ induces a bijection of depth zero raft spaces of coarse dimension $>n$. By Step 1,
for each depth zero raft $R \subset T$ such that $X_R$ has coarse dimension $n$ there exists a depth zero raft $R' \subset
T'$ such that the raft space $X'_{R'}$ has coarse dimension $n' \le n$ and $f(X_R) \csubset{c} X'_{f_\#(R)}$. 

\subparagraph{Case 1:} $R' = \{w\}$ is a one vertex raft. The raft space $X_R$ contains a uniformly properly embedded coarse
$\PD(n)$ subspace $Y$: when $R = \{v\}$ is a one vertex raft then $Y=X_v$; when $R$ is bushy then $Y=X_L$ for some bi-infinite
line $L\subset R$. The map $f$ therefore restricts to a uniformly proper embedding of a coarse $\PD(n)$ space into the coarse
$\PD(n')$ space $X'_{R'}$. Applying Lemma~\ref{LemmaBigInSmall} we conclude that $n'=n$.

\subparagraph{Case 2:} $R'$ is a bushy raft. By \cite{FarbMosher:ABC} there is a uniformly proper embedding of the raft space
$X'_{R'}$ into a coarse $\PD(n'+1)$ space $Z$, such that $X'_{R'} \cstrict{c} Z$. As in Case~1, $X_R$ contains a uniformly properly
embedded coarse $\PD(n)$ space $Y$. We thus obtain a uniformly proper embedding $F \from Y \inject Z$. Applying
Lemma~\ref{LemmaBigInSmall}, $n' = n$ or $n-1$. If $n'=n-1$, then $n'+1=n$, and from
Lemma~\ref{LemmaBigInSmall} we conclude that $F(Y) \ceq{c} Z$, but this contradicts $F(Y) \csubset{c} X'_{R'} \cstrict{c} Z$.

\bigskip

This completes the proof of the Depth Zero Vertex Rigidity Theorem~\ref{TheoremVertexRigidity}.

\section{Finite depth graphs of groups}
\label{SectionFiniteDepth}

\subsection{Definitions and examples}

Let $\G$ be a finite type graph of groups with Bass Serre tree of spaces $X \to T$. Let $G =
\pi_1\G$.

Let $\VE(T)$ be the set of vertices and edges of $T$. The metric on $T$ induces a metric on
$\VE(T)$, via a natural injection $\VE(T) \to T$ which takes each vertex to itself and each edge
to its midpoint, so any vertex-to-vertex or edge-to-edge distance is an integer while any
vertex-to-edge distance is a half integer.

There is a partial ordering of $\VE(T)$ by coarse inclusion of the associated spaces. Note that
the equivalence relation generated by this partial ordering is the same as coarse equivalence of
the associated spaces. Every maximal element of $\VE(T)$ is coarsely equivalent to a vertex,
because every edge space is coarsely included into both adjacent vertex spaces. Every nonmaximal
element of $\VE(T)$ is coarsely equivalent to an edge space: let $v$ be a vertex which is not
maximal in $\VE(T)$, and so $X_v\csubset{c} X_a$ for some vertex or edge $a \ne v$; letting $e$ be
the first edge on the path from $v$ to $a$, we have $X_v \ceq{c} X_e$. 

The \emph{depth} of $a \in \VE(T)$ is defined to be the maximum length $n$ of a
strictly increasing chain $X_a=X_{a_n} \cstrict{c} \cdots \cstrict{c} X_{a_0}$. The depth can be
an integer $\ge 0$ or $\infinity$. Each maximal element of $\VE(T)$ has depth zero. Depth is
invariant under the action of $G$. Let $\VE_n(T) = \{a \in \VE(T) \suchthat \depth(a) \le n\}$. Let
$\forest_n =\union \VE_n(T)$, and note that $\forest_n$ is a subcomplex of $T$, because the depth
of an edge is no less than the depths of its endpoints. We thus have a $G$-invariant filtration of
$T$ by subforests:
$$\forest_0 \subset \forest_1 \subset \cdots \subset \forest_\infinity = T
$$
We say that $\G$ has \emph{finite depth} if $\forest_n = T$ for some finite $n$. The least
such $n$ is the \emph{depth} of $\G$. 


If $\G$ has finite depth then the depth is at most $\#\text{edges}(\G)$. To see why,
consider a strictly increasing chain $X_{a_n} \cstrict{c} \cdots \cstrict{c} X_{a_1}  \cstrict{c}
X_{a_0}$. If $1 \le i \le n$ then $X_{a_i}$ is coarsely equivalent to an edge space, so we may
assume $a_i$ is an edge. If $1 \le i < j\le n$ then $a_i \ne a_j$. If $n >
\#\text{edges}(\G)$ then the sequence $a_1,\ldots,a_{n}$ contains two distinct edges in the
same orbit of the $G$-action, $X_e\cstrict{c} X_{g\cdot e}$ for some edge $e$ and some $g
\in G$, and we obtain an infinite sequence
$$X_e \cstrict{c} X_{g \cdot e} \cstrict{c} X_{g^2 \cdot e} \cstrict{c} \cdots
$$
contradicting that $\G$ has finite depth. The same argument shows that even if $\G$ has
infinite depth, the filtration of $T$ by subforests $\forest_n$ stabilizes for finite values $n \ge
\#\text{edges}(\G)$. 

The finite depth property can be reformulated in a purely algebraic manner as follows. Recall that
the orbit map from $G$ to $X$, given by $\O(g) = g \cdot p$ where $p \in X$ is a base point, is a
quasi-isometry. For each $a \in \VE(T)$, since $\Stab(a)$ acts coboundedly on $X_a$ we have
$\O(\Stab(a)) \ceq{c} X_a$ (but this coarse equivalence is \emph{not} uniform). It follows that
the coarse inclusion lattice of vertex and edge spaces in $X$ is isomorphic to the coarse
inclusion lattice of vertex and edge stabilizers in $G$, which by
Corollary~\ref{CorollaryCommensurable} is isomorphic to the commensurability lattice of vertex and
edge stabilizers in $G$. This proves:

\begin{lemma}
\label{LemmaCommensurabilityDepth}
$\G$ has finite depth if and only if there is a bound to the length of any totally ordered subset
of the commensurability lattice of vertex and edge stabilizers in $\pi_1\G$.
\qed\end{lemma}

\paragraph{Coarse structure of finite depth graphs of groups.} 
A depth zero graph of groups is the
same as a geometrically homogeneous graph of groups. 

In general the depth zero subforest $\forest_0$ can be described as follows. Each component $R$ of
$\forest_0$ is a subtree of $T$ with the property that for each edge $e$ of $R$ and each endpoint
$v$ of $e$ the inclusion $G^e \subset G^v$ has finite index, and for each edge $e \in T-R$ incident
to a vertex $v \in R$, the inclusion $G^e \subgroup G^v$ has infinite index. A component $R$ of
$\forest_0$ will be called a \emph{depth~0 raft} of $T$, and a component of $F_0$ consisting of a
single vertex of $T$ is called a \emph{fat vertex} of $T$. 

We now define depth~$n$ rafts inductively. Consider the $G$-tree $T / \forest_{n-1}$ obtained from
$T$ by collapsing each component of $\forest_{n-1}$ to a single vertex. Note that
$\forest_n / \forest_{n-1}$ is a subforest of $T / \forest_{n-1}$. A \emph{raft in $\forest_n /
\forest_{n-1}$}, is defined to be either of the following two objects:
\begin{itemize}
\item a vertex of $\forest_n / \forest_{n-1}$ without any incident edges,
representing a component of $\forest_{n-1}$ without incident edges of depth $n$;
\end{itemize} 
or
\begin{itemize} 
\item a coarse equivalence class of vertices and edges in $\forest_n - \forest_{n-1}$ together
with any components of $\forest_{n-1}$ incident to such an edge. 
\end{itemize}
Each raft in $\forest_n / \forest_{n-1}$ is a subcomplex: if the edge $e$ of $\forest_n /
\forest_{n-1}$ is incident to the vertex $v$ then either $v$ represents a component of
$\forest_{n-1}$, or $v \not\in \forest_{n-1}$ in which case $v$ has depth~$n$ and $X_e \ceq{c}
X_v$. Each raft in $\forest_n / \forest_{n-1}$ is a subtree, because if $e,e'$ are coarsely
equivalent edges of $\forest_n/\forest_{n-1}$, and if $e=e_0,e_1,\ldots,e_k=e'$ is the edge path in
$T$ from $e$ to $e'$, then each of the edge spaces $X_{e_i}$ coarsely contains $X_e$ and $X_{e'}$,
and so each of $e_0,\ldots,e_k$ is either in $\forest_{n-1}$ or is coarsely equivalent to $e$ and
$e'$. It immediately follows that rafts in $\forest_n / \forest_{n-1}$ decompose $\forest_n /
\forest_{n-1}$ into a collection of subtrees any two of which intersect either in the empty set or
in a fat vertex of $\forest_n / \forest_{n-1}$.

A \emph{depth $n$ raft in $T$} is the pullback of a raft in $\forest_n / \forest_{n-1}$. The
intersection of any two depth $n$ rafts in $T$ is either empty or a component of $\forest_{n-1}$.
Each component of $\forest_n$ is a maximal connected union of depth $n$ rafts, called a
\emph{depth $n$ flotilla}. Note that each depth~0 flotilla consists of a single depth~0 raft, but
the same need not be true in higher depths.

Irreducibility of $\G$ immediately implies the fact that for each $n \ge 0$, each depth
$n+1$ raft $R$ is the convex hull in $T$ of the depth $n$ flotillas contained in $R$. Since the
subgroup of $\pi_1\G$ stabilizing $R$ acts cocompactly on $R$, we get the following uniform
version of this fact:

\begin{proposition}
\label{PropReducibility}
If $\G$ is irreducible, then for each $n \ge 0$ and each depth $n+1$ raft $R$ there
exists a constant $s$ such that for each depth $n+1$ vertex or edge $a$ in $R$ there is a path in
$R$ of length at most $s$ whose endpoints have depth $\le n$ and whose interior contains
$a$. 
\qed\end{proposition}

\paragraph{Examples.} Here are some examples and nonexamples of finite depth graphs of groups.

\subparagraph{Example.} Let $\G$ have three vertices $u,v,w$ on an arc with two edges $e,f$, $\bdy e
= \{u,v\}$ and $\bdy f = \{v,w\}$. Let $\G_u = \<a_1,a_2,a_3\>$ be free abelian of rank~3, $\G_v =
\<b_1,b_2,b_3\>$ free abelian of rank~3, and $\G_w = \<c_1,c_2\>$ free abelian of rank~2. Let $\G_e$
be free abelian of rank~2, mapped to $\<a_1,a_2\>$ in $\G_u$ and to $\<b_1,b_2\> = \G_v$. Let $\G_f$
be infinite cyclic, mapped to $\<b_1\>$ in $\G_v$ and to $\<c_1\>$ in $\G_w$. Let $X \to T$ be the
Bass-Serre tree of spaces. Each vertex of $T$ is a depth~0 raft. EAch lift of the edge $e$ is of
depth~1, and any component of the lift of the subgraph $\{u,e,v\}$ is a depth~1 raft. Each lift of
the edge $f$ is of depth~2, and the entire Bass-Serre tree is a single depth~2 raft, so $\G$ has
depth~2. Note that for any lift $\tilde f$ of $f$, with endpoints $\tilde v$, $\tilde w$ lifting
$v,w$, respectively, the edge space of $\tilde f$ is the coarse intersection of the vertex spaces of
$\tilde v$ and $\tilde w$, and yet $\tilde f$ is not of depth~1.

\subparagraph{Example.} Given an amalgamated free product $A *_C B$, if $A,B$ are word hyperbolic,
and if the inclusions $C \inject A$, $C \inject B$ have quasiconvex images, then the graph of
groups $A *_C B$ has depth at most one. To see why, let $X \to T$ be the Bass-Serre tree of spaces.
The subgroup $C \subgroup A$ has finite width \cite{GMRS:widths}, which means that there is a
constant $N$ so that out of any $N$ distinct conjugates of $C$ in $A$ one can find two of them that
intersect in a finite subgroup. The subgroup $C \subgroup B$ also has finite width. It follows
that for any vertex $v \in T$ and edges $e,f\subset T$ incident to $v$, either $X_e \ceq{c} X_f$
or $X_e\cintersect{c} X_f$ is strictly coarsely contained in $X_e$ and in $X_f$: the only other
alternative is that $X_e \cstrict{c} X_f$, but this implies that some conjugate of $C$ in $A$ is
strictly coarsely contained in another, violating finite width. It follows that if $e,f$ are any
two edges of $T$ then $X_e$ is not strictly coarsely contained in $X_f$, for if $X_e\csubset{c}
X_f$ and if $e=e_0,e_1,\ldots,e_K=f$ is the path in $T$ from $e$ to $f$, then $X_e\ceq{c} X_e
\cintersect{c} X_f \csubset{c} X_{e_k}$ for each $k$, implying inductively that $X_e =
X_{e_0} \ceq{c} X_{e_1}\ceq{c} X_{e_2} \ceq{c} \cdots \ceq{c} X_{e_K}= X_f$. 

The same argument shows that if $\G$ is a tree of groups with word hyperbolic vertex groups and
quasiconvex edge-to-vertex injections, and if it is true that for each vertex $v$ and edges $e \ne
f$ at $v$ no conjugate of $\G_e$ in $\G_v$ is strictly coarsely contained in $\G_f$, then $\G$ has
depth at most~$1$.

\subparagraph{Question:} Given a tree of groups with word hyperbolic vertex groups and quasiconvex
edge-to-vertex injections, is it finite depth?

\subparagraph{Example: $\PD(n)$ graphs of groups.} If $\G$ is a finite graph of groups
whose edge groups are all \Poincare\ duality groups then $\G$ has finite depth.
More generally, if every edge group of $\G$ has the property that no finite index subgroup
is isomorphic to an infinite index subgroup then $\G$ has finite depth. This property is
true for \Poincare\ duality groups because any infinite index subgroup of a $\PD(n)$ group has
smaller virtual cohomological dimension \cite{Brown:cohomology}, and a similar fact holds for
coarse $\PD(n)$ groups.

Sela proved that if $G$ is a torsion free, one ended, word hyperbolic group then $G$ is not
isomorphic to any proper subgroup \cite{Sela:RigidityII}, but it is unknown
whether $G$ can have a finite index subgroup isomorphic to an infinite index subgroup.

\subparagraph{Nonexample.} Start with a vertex group $\G_v$ which has a free subgroup $F$ of
finite rank. Choose a finite index subgroup $F_1 \subgroup F$ and an infinite index subgroup $F_2
\subgroup F$ such that $\rank(F_1) = \rank(F_2)$. Now form an HNN amalgamation, adding an edge $e$
with edge group $\G_e$ indentified on one end with $F_1$ and on the other end with $F_2$, and
clearly the edge $e$ has infinite depth. Nevertheless, as long as $F$ has infinite index in
$\G_v$, the vertex group $\G_v$ has depth zero in $\pi_1\G$, and so one can construct such
examples to which the Vertex Rigidity Theorem~\ref{TheoremVertexRigidity} applies by taking $\G_v$
to be coarse $\PD(n)$ with $n \ge 3$.

\paragraph{An algorithm, using oracles, to decide finite depth.} It seems to be hard in general to
decide whether or not a graph of groups $\G$ is of finite depth.
Lemma~\ref{LemmaCommensurabilityDepth} reduces this to the existence of an algorithm that solves
the following: 

\begin{description}
\item[Oracle:] Given a finitely generated group $G$ and a finite collection of finitely
generated subgroups $H_1,\ldots,H_n \subgroup G$, compute the strict coarse inclusion lattice for
the collection of subgroups conjugate to one of $H_1,\ldots,H_n$: given $g,g' \in G$ and $i,j \in
1,\ldots,n$, decide whether $g H_i g^\inv$ is commensurable to a subgroup of $g' H_j g'{}^\inv$.
\end{description}

We can improve matters slightly as follows. If we are not willing to apply the oracle to the
ambient group $\pi_1\G$ and its vertex and edge groups, we may at least be willing to apply the
oracle to the vertex groups of $\G$ and their incident edge groups. If so, then we obtain the
following relative algorithm to decide whether $\G$ has finite depth.

To start the relative algorithm, we can use the oracle to tell us which edge-to-vertex injections
in $\G$ have finite index image, and from that we can compute the depth zero rafts. If the depth
zero rafts are empty then $\G$ does not have finite depth. Otherwise, crush each depth zero raft
to a point, after which we may assume that each depth zero raft is a fat vertex.

For each fat vertex group $G_v$, apply the oracle to the incident edge groups $H_1,\ldots,H_n$,
deciding which of $H_1,\ldots,H_n$ are not maximal in the commensurability lattice of $G_v$.
Remove those edges from the graph of groups, and continue removing edges inductively: at each
vertex where an incident edge has already been removed, apply the oracle to tell whether some
new incident edge group has a conjugate that is strictly coarsely contained in a removed edge
group, and if so remove the new edge. This process must stop after finitely many steps because the
graph of groups is finite. If all edges have been removed then $\G$ does not have finite depth and
the algorithm stops. Otherwise, the edges and nonfat vertices that are left have depth~1, and their
union gives the depth~1 flotillas.

Now the induction can be continued: assuming that depth $n$ flotillas have been identified,
collapse each depth $n$ flotilla to a point and continue the induction in the same manner.

\subparagraph{Question:} Is the oracle decidable for quasiconvex subgroups of word hyperbolic groups?

\bigskip

Note that if $G$ is word hyperbolic and $H_1,H_2 \subgroup G$ are quasiconvex subgroups, then
$H_1 \csubset{c} H_2$ if and only if the boundary of $H_1$ is a subset of the boundary of $H_2$ in
the boundary of $G$, and $H_1 \cstrict{c} H_2$ if and only if the inclusion of boundaries is
proper.

\subsection{Proof of the Vertex--Edge Rigidity Theorem \ref{TheoremVERigidity}}

Recall the setting of the theorem: we are given $\G,\G'$ finite type, finite depth graphs of
groups, $X\to T$, $X' \to T'$ their Bass-Serre trees of spaces, and $f \from X \to X'$ a
quasi-isometry that coarsely respects depth zero vertex spaces, meaning: 
\begin{itemize}
\item There exists $K \ge 1$, $C \ge 0$ and a $K,C$ quasi-isometry $f_\#\from\V_0(T)\to\V_0(T)$ such
that for each $v \in \V_0(T)$ we have $f(X_v)\ceq{[C]} X'_{f_\#(v)}$, and for each $v' \in \V_0(T')$
there exists $v \in \V_0(T)$ such that $X'_{f_\#(v)} \ceq{[C]} X'_{v'}$. 
\end{itemize}
We shall prove, by induction on $n$, that $f$ coarsely respects vertex and edge spaces of
depth $\le n$, meaning:
\begin{itemize}
\item There exists a quasi-isometry $f_\# \from \VE_n(T) \to \VE_n(T')$ such that for each $a \in
\VE_n(T)$ we have $f(X_a) \ceq{[C']} X'_{f_\#(a)}$, and for each $a' \in \VE_n(T')$ there exists $A
\in \VE_n(T)$ such that $X_{f_\#(a)} \ceq{[C']} X'_{a'}$. 
\end{itemize}
This map $f_\#$ will be an extension of the $f_\#$ already given. Also, the constant $C'$ and the
quasi-isometry constants of the extended $f_\#$ will depend only on the depth zero constants $K,C$
and the quasi-isometry constants of $f$. 

When $n$ is equal to the maximum of the depths of $\G$ and $\G'$, the statement above gives the
conclusion of Theorem~\ref{TheoremVERigidity}.

To start the induction, first extend the given quasi-isometry $f_\# \from \V_0(T) \to \V_0(T')$ to a
map $f_\# \from \VE_0(T) \to \VE_0(T')$ so that for any edge $e \in \VE_0(T)$, $f_\#(e) = f_\#(v)$
where $v$ is either endpoint of $e$. The extension clearly exhibits that $f$ coarsely respects
vertex \emph{and edge} spaces of depth zero.

Preparatory to the general inductive step, we consider the inductive step from $n=0$ to $n=1$. Note
that the collection of depth zero vertices and edges in $T$ coarsely determines the collection of
depth zero rafts, by an explicit formula using coarse language: up to coarse equivalence, a depth
zero raft space is the same thing as the union of a coarse equivalence class of depth zero vertices
and edges. The same formula coarsely determines depth zero rafts in $X'$, in terms of depth zero
vertex and edge spaces. Since $f$ coarsely respects depth zero vertex and edge spaces, and since $f$
coarsely respects any particular formula in coarse language, it immediately follows that $f$
coarsely respects depth zero rafts, in the sense that $f$ induces a bijection between depth zero
rafts of $T$ and depth zero rafts of $T'$, characterized by saying that rafts $R \subset T$, $R'
\subset T'$ correspond under this bijection if and only if $f(X_R) \ceq{c} X'_{R'}$. 

Moreover, the quasi-isometry $f$ \emph{uniformly} coarsely respects depth zero rafts, meaning
that there exists $C' \ge 0$, depending only on the depth zero constants and the quasi-isometry
constants of $f$, such that for any depth zero rafts $R\subset T$, $R'\subset T'$, we have $f(X_R)
\ceq{c} X'_{R'}$ if and only if $f(X_R) \ceq{[C']} X'_{R'}$. The existence of $C'$ follows from the
fact that $X'_{R'}$ is the union of $X'_{a'}$ for $a'\in R'$, and $f(X_a) \ceq{[C]} X'_{a'}$ for $a'
= f_\#(a) \in R'$, and $X'_{a'} \ceq{[kC]} f(X_a)$ for some $a\in R$ such that $f_\#(a)$ has distance
at most $k$ from~$a'$.

Since the collection of depth zero flotillas is identical to the collection of depth zero rafts,
we can also say that $f$ uniformly coarsely respects the depth zero flotillas. 

We now recast the induction as follows. Given $n \ge 0$, suppose that $f$ coarsely respects the
vertex and edge spaces, rafts, and flotillas of depth $\le n$. We shall prove:
\begin{itemize}
\item The collection
of depth~$n$ flotillas in $\G$ coarsely determines the collection of depth $n+1$ vertex and edge
spaces, in the sense that the depth $n+1$ vertex and edge spaces can be expressed in terms of the
depth $n$ flotillas using formulas in coarse language: coarse equivalence, coarse inclusion,
strict coarse inclusion, coarse intersection, as well as ordinary union. 
\end{itemize}
The proof of this
statement will use irreducibility of $\G$. The formulas produced in the course of the proof do not
depend on $\G$, and so the same formulas determine the depth $n+1$ vertex and edge spaces of $\G'$
in terms of the depth $n$ vertex and edge spaces of $\G'$, since $\G'$ is also irreducible. It
immediately follows that $f$ coarsely respects the depth $n+1$ edge and vertex spaces. We then
show, in a similar fashion using formulas in coarse language that do not depend on the graph of
groups, that the depth $n+1$ edge and vertex spaces coarsely determine the depth $n+1$ rafts, and
that the depth $n+1$ rafts coarsely determine the depth $n+1$ flotillas. It follows that $f$
coarsely respects rafts and flotillas in depth $n+1$, completing the induction. 

In all of the arguments involved in the induction step, the operations of coarse set theory come
in non-uniform versions without the constants, such as $A\ceq{c}B$, versus uniform versions with
the constants, such as $A\ceq{[r]}B$. We will write out the proof in detail using nonuniform
coarse language, indicating later how to make it uniform, that is, how to show that the depth
$n+1$ constants for $f$ depend only on the depth $n$ constants.

Recall that $\forest_n$ denotes the subforest of $T$ consisting of vertices and edges of depth
$\le n$. Let $\abs{\forest_n}$ be the component set of $\forest_n$, so the elements of
$\abs{\forest_n}$ are precisely the depth $n$ flotillas in $T$. For each $A \ne B \in
\abs{\forest_n}$ consider the coarse intersection $\Xsub{A} \cintersect{c}
\Xsub{B}$. We define a partial ordering on the set of distinct unordered pairs $\{A,B\}$ in $\abs{\forest_n}$
by coarse inclusion of coarse intersections:
$\{A,B\}\csubset{c} \{A',B'\}$ means
$$\left( \Xsub{A} \cintersect{c} \Xsub{B} \right) \csubset{c} \left( \Xsub{A'} \cintersect{c}
\Xsub{B'} \right)
$$
We similarly define $\{A,B\} \ceq{c} \{A',B'\}$, and we say that $\{A,B\}$ is \emph{maximal} with
respect to $\csubset{c}$ if $\{A,B\} \csubset{c} \{A',B'\}$ implies $\{A,B\} \ceq{c} \{A',B'\}$.

For each depth $n+1$ vertex or edge $a$, there exists a distinct unordered pair $\{A,B\}$ in
$\abs{\forest_n}$ such that $\Xsub{A} \cintersect{c} \Xsub{B} = X_a$. To see why, let $F \subset
T$ be the depth $n+1$ raft containing $a$. Since $\G$ is irreducible, we may apply
Proposition~\ref{PropReducibility} to obtain a distinct unordered pair $\{A,B\}$ in
$\abs{\forest_n}$ such that
$A,B
\subset F$, $a$ is contained in the interior of the segment $\overline{AB}$, and every vertex or
edge $b$ in the interior of $\overline{AB}$ has depth $n+1$. For each such $b$ it follows that
$\Xsub{A} \cintersect{c} \Xsub{B} \ceq{c} X_b \ceq{c} X_a$.

We claim that for each distinct unordered pair $\{A,B\}$ in $\abs{\forest_n}$ the following are
equivalent:
\begin{enumerate}
\item[(1)] \label{ItemMaximalPair}
$\{A,B\}$ is maximal with respect to $\csubset{c}$;
\item[(2)] \label{ItemEquivToEdgeSpace}
$\Xsub{A} \cintersect{c} \Xsub{B}$ is coarsely equivalent to some depth~$n+1$ edge space $X_e$;
\item[(3)] \label{ItemUniqueRaft}
There exists a unique depth~$n+1$ raft $R$ of $T$ containing both $A$ and~$B$.
\end{enumerate}
Note that the argument above shows that each depth $n+1$ edge space $X_e$ occurs as in
(\ref{ItemEquivToEdgeSpace}), for some distinct unordered pair $\{A,B\}$ in $\abs{\forest_n}$.
We also claim that if the equivalent statements (\ref{ItemMaximalPair}),
(\ref{ItemEquivToEdgeSpace}), (\ref{ItemUniqueRaft}) hold then, letting $R$ be the raft in
(\ref{ItemUniqueRaft}), the following are also true:
\begin{enumeratecontinue}
\item[(4)] 
\label{ItemRaftVertex}
For each depth~$n+1$ edge or vertex $a$ of $T$, we have $a \in R$ if and
only if $\Xsub{A} \cintersect{c} \Xsub{B} \ceq{c} X_a$. 
\item[(5)] 
\label{ItemFlotillaRafts}
For each depth~$n$ flotilla $S$ of $T$, we have $S \subset R$ if and only if $X_e \csubset
{c} \Xsub{S}$ for each depth~$n+1$ edge $e$ of $R$.
\end{enumeratecontinue}

Once the claim is established, then we can show that the depth $n$ flotillas coarsely determine
the depth $n+1$ vertex and edge spaces, which coarsely determine the depth $n+1$ rafts, which
coarsely determine the depth $n+1$ flotillas, as follows. First, the equivalence of (1), (2), (3)
shows that depth~$n+1$ vertex and edge spaces are coarsely determined by the depth~$n$ flotillas: a
depth~$n+1$ vertex or edge space is precisely what you get by taking the coarse intersection of
$\Xsub{A}$ and $\Xsub{B}$, where $\{A,B\}$ is a distinct unordered pair in $\abs{\forest_n}$ that
is maximal with respect to $\csubset{c}$. Second, each coarse equivalence class of depth~$n+1$
vertex and edge spaces coarsely determines a unique depth~$n+1$ raft $R$, the depth~$n+1$ edges and
vertices contained in $R$ are coarsely determined by property~(\ref{ItemRaftVertex}), and the
depth~$n$ flotillas contained in $R$ are coarsely determined by property~(\ref{ItemFlotillaRafts}),
so the entire raft $R$ is coarsely determined. Third, given two depth~$n+1$ rafts $R \ne R'$, the
property that $R\intersect R' \ne \emptyset$ is determined by property~(\ref{ItemFlotillaRafts}),
by saying that $\Xsub{R} \cintersect{c} \Xsub{R'} = \Xsub{S}$ for some depth~$n$ flotilla $S$;
note that if $R \intersect R' = \emptyset$ then $\Xsub{R}\cintersect{c} \Xsub{R'}$ is coarsely
contained in $X_a$ for some vertex or edge $a$ of depth $\ge n+1$ lying on the shortest path
connecting $R$ to $R'$, and $X_a$ is not coarsely equivalent to $X_S$ for any depth $n$ flotilla
$S$. Finally, the property that $R,R'$ are in the same flotilla is determined by existence of a
finite chain of depth $n+1$ rafts $R=R_0,R_1,\ldots,R_K=R'$ such that $R_{k-1}\intersect R_k \ne
\emptyset$ for $k=1,\ldots,K$. 

Now we prove (1)$\implies$(3). Let $\{A,B\}$ be maximal with respect to $\csubset{c}$. Consider
the segment $\overline{AB}$ in $T / F_n$. No interior vertex or edge $a$ of $\overline{AB}$ has
depth $\ge n+2$, because otherwise $\Xsub{A} \cintersect{c} \Xsub{B}\csubset{c} X_a \cstrict{c} X_b$
for some depth $n+1$ vertex or edge $b$, but $X_b \ceq{c} \Xsub{A'} \cintersect{c} \Xsub{B'}$ for
some distinct unordered pair $\{A',B'\}$ in $\abs{\forest_n}$, contradicting maximality of $\{A,B\}$. Any
two edges $e,f$ in the interior of $\overline{AB}$ are coarsely equivalent, because otherwise
$\Xsub{A}\cintersect{c} \Xsub{B} \csubset{c} X_e \cintersect{c} X_f \cstrict{c} X_e \ceq{c} \Xsub{A'}
\cintersect{c} \Xsub{B'}$ for some distinct unordered pair $\{A',B'\}$ in $\abs{\forest_n}$, again
contradicting maximality of $\{A,B\}$. This shows that $A,B$ are vertices of a raft in $F_{n+1} /
F_n$, and so $A,B$ are contained in some depth $n+1$ raft $R$ of $T$. The raft $R$ is unique,
because distinct depth $n+1$ rafts are either disjoint or intersect in a single depth $n$ flotilla.

To prove (3)$\implies$(2), if $A,B$ satisfy (3) then $\Xsub{A} \cintersect{c} \Xsub{B}$ is coarsely
equivalent to $X_e$ for any depth $n+1$ edge $e$ in the interior of $\overline{AB}$.

To prove (2)$\implies$(1), suppose $A,B$ satisfy (2) but not (1), and so there is a distinct
unordered pair $\{C,D\}$ in $\abs{\forest_n}$ such that $\{A,B\} \cstrict{c} \{C,D\}$. The interior of
the segment $\overline{CD}$ contains some edge $f$ of depth $\ge n+1$, and we have
$$X_e \csubset{c} \Xsub{A} \cintersect{c} \Xsub{B} \cstrict{c} \Xsub{C} \cintersect{c} \Xsub{D}
\csubset{c} X_f
$$
contradicting that $e$ has depth $n+1$.

To prove (4), it already follows from (2) and (3) that $\Xsub{A} \cintersect{c} \Xsub{B} \ceq{c}
X_a$ for each depth~$n+1$ edge or vertex $a$ in $R$. Conversely, if $\Xsub{A}
\cintersect{c} \Xsub{B} \ceq{c} X_a$ where $a$ is a depth $n+1$ vertex or edge of $T$, then $X_a
\ceq{c} X_b$ for any depth $n+1$ vertex or edge $b$ of $R$, and so $a$ is in $R$.

To prove (5), consider a depth $n$ flotilla $S$. If $S \subset R$ then there exists a depth $n+1$
edge $e$ of $R$ incident to $S$, and so $X_e \csubset{c} X_S$, implying that $X_a \csubset{c} X_S$
for any depth $n+1$ edge or vertex $a$ in $R$. Conversely, suppose that $S \not\subset R$. The
shortest path in $T$ connecting $R$ to $S$ contains some edge $b$ of depth $\ge n+1$ such that
$X_b$ is \emph{not} coarsely equivalent to $X_a$ for any depth $n+1$ vertex or edge space $a$ of
$R$. It follows that $X_a\cintersect{c} X_b \cstrict{c} X_a$. But we also have $X_a
\cintersect{c} X_S \csubset{c} X_b$, implying that $X_a \cintersect{c} X_S \cstrict{c} X_a$, which
implies that $X_a$ is \emph{not} coarsely contained in $X_S$.

To complete the proof of the Vertex--Edge Rigidity Theorem, we sketch how to make the above
arguments uniform. 

For example, we showed above that if $a$ is a depth $n+1$ vertex or edge then $X_a
\ceq{c} X_A \cintersect{c} X_B$ for some distinct pair of depth $n$ flotillas $X_A$, $X_B$, by
taking $a$ to lie in the interior of a path $\overline{AB}$ in $F_{n+1} / F_n$.
Proposition~\ref{PropReducibility} shows that we can choose $\overline{AB}$ to have uniformly
bounded length, which implies that there exists a constant $r$ such that $X_A \cintersect{[r]}
X_B$ represents the coarse intersection $X_A \cintersect{c} X_B$, and such that $X_a \ceq{[r]} X_A
\cintersect{[r]} X_B$.  This shows that the depth $n+1$ edge and vertex spaces are
\emph{uniformly} coarsely determined by the depth $n$ flotillas. Assuming by induction that $f$
\emph{uniformly} coarsely respects depth
$n$ flotillas, it follows that the extension of the map $f_\# \from \VE_n(T) \to
\VE_n(T')$ to a map $f_\# \from \VE_{n+1}(T) \to \VE_{n+1}(T')$ is a quasi-isometry with constants
depending only on $f$, and with a uniformly bounded Hausdorff distance between $f(X_a)$ and
$X'_{f_\#(a)}$ for $a \in \VE_{n+1}(T)$. Similar arguments show that $f$ uniformly coarsely
respects depth $n+1$ rafts and flotillas.

This finishes the proof of the Vertex--Edge Rigidity Theorem.

\subsection{Reduction of finite depth graphs of groups}

In this section we give an application of the Vertex--Edge Rigidity
Theorem~\ref{TheoremVERigidity}. 

We saw in Section~\ref{SectionIrreducible} how Higman's group was
used to produce an example of a graph of groups $\G$, not of finite depth, such that two different
reduction processes applied to $\G$ produced two irreducible graphs of groups $\G_1,\G_2$, where the
set of quasi-isometry types of vertex and edge groups in $\G_1$ is not equal to the set of
quasi-isometry types of vertex and edge groups in $\G_2$. 

In contrast, Section~\ref{SectionIrreducible} contains another example where the integer
Heisenberg group is used to produce an example of a finite depth graph of groups $\G$ in which two
different reduction processes applied to $\G$ produced irreducible graphs of groups $\G_1,\G_2$
whose vertex and edge groups, while not isomorphic, were abstractly commensurable. The following
proposition shows that this phenomenon is general among finite depth graphs of groups.

We show that the reduction process applied to any finite type, finite depth graph of groups
$\G$ results in a graph of groups whose vertex and edge groups are well-defined up to abstract
commensurability, in fact they are well-defined up to commensurability in $\pi_1\G$. The
statement of this application does not involve quasi-isometries, and probably there is a proof not
involving quasi-isometries, involving a careful analysis of what happens to the commensurability
types of the edge and vertex groups under the collapsing process. However, a proof using the
Classification Theorem is quite short as we shall see.

\begin{proposition}
\label{PropCollapse}
Let $\G$ be a finite type, finite depth graph of groups. Let $\G_1,\G_2$ be two irreducible
graphs of groups each obtained from $\G$ by iteratively collapsing along reducing edges. Then the
set of abstract commensurability types of vertex and edge groups of $\G_i$ is independent of $i$.
More precisely, if $T_1,T_2$ are the Bass-Serre trees of $\G_1,\G_2$, on which $\pi_1\G$ acts,
then the set of $\pi_1\G$ commensurability classes of vertex and edge stabilizers in
$T_i$ is independent of $i$.
\end{proposition}

\begin{proof}
Suppose $\G'$ is obtained from a finite type graph of groups $\G$ by collapsing an edge. We may
identify $\pi_1\G \approx \pi_1\G'$, acting on both of the Bass-Serre trees of spaces $X \to T$,
$X' \to T'$, and this identification induces a quasi-isometry $F \from X' \to X$.

We claim that $F$ coarsely respects depth zero edge and vertex spaces. To see why, let
$e$ be the edge of $\G$ which is collapsed, with incident vertices $v \ne w$ so that
$\G_e \to \G_v$ is an isomorphism. Let $f \from \G' \to \G$ be the collapse map, let
$\tilde f \from T' \to T$ be an equivariant lift, and let $F \from X' \to X$ be an
induced quasi-isometry. Note that $e$ has depth zero only if $\G_e\to \G_w$ has finite index
image; in this case it follows that for each vertex or edge $a$ of $T'$, the stabilizer of
$\tilde f(a)$ is a finite index supergroup of the stabilizer of $a$, and so $F(X'_a) \ceq{c}
X_{f(a)}$. On the other hand, suppose that $\G_e \to \G_w$ has infinite index image, and so
neither $e$ nor $v$ has depth zero. It follows that $f$ induces an isomorphism between the depth
zero subforest of $T$ and the depth zero subforest of $T'$, and for each depth zero edge or vertex
$a$ of $T'$, the stabilizer of $f(a)$ equals the stabilizer of $a$, and $F(X'_a) \ceq{c} X_{f(a)}$.
This proves the claim.

Now let $X_i \to T_i$ be the Bass-Serre tree of spaces for $\G_i$. By iterating a sequence of edge
collapses we obtain a map from $\G$ to $\G_i$ that induces an identification $\pi_1\G \approx
\pi_1\G_i$, and so we obtain an identification $\pi_1\G_1 \approx\pi_1\G_2$, which induces a
quasi-isometry $F \from X_1 \to X_2$. By applying the above claim it follows that $F$ coarsely
respects depth zero vertex and edge spaces. Since $\G_1,\G_2$ are irreducible we may apply
Theorem~\ref{TheoremVERigidity}, from which it follows that $F$ coarsely respects the entire
collection of vertex and edge spaces. Let $a$ be a vertex or edge of $T_1$ and $a'$ a vertex or
edge of $T_2$ such that $F((X_1)_a) \ceq{c} (X_2)_{a'}$. It follows that $\Stab(a)$ and
$\Stab(a')$ are coarsely equivalent in $\pi_1\G$ with respect to the word metric, and so
$\Stab(a)$ and $\Stab(a')$ are commensurable in $\pi_1\G$, by
Corollary~\ref{CorollaryCommensurable}. 
\end{proof}

The proof of Proposition~\ref{PropCollapse} seems to indicate that if $\G$ is a finite type,
finite depth graph of groups, then the vertex and edge groups of a complete reduction $\G'$ of $\G$
can be picked out from among the vertex and edge groups of $\G$ itself by some quasi-isometrically
invariant property. This is indeed true, at least up to commensurability, and we can state this
property very simply. Let $T,T'$ be the Bass-Serre trees of $\G,\G'$ respectively. First of all,
the commensurability classes of the stabilizers of depth zero vertex and edge groups are invariant
under reduction. Second, given a vertex or edge $a$ of $T$ of depth $\ge 1$, the commensurability
class of the stabilizer of $a$ survives as the commensurability class of a vertex or edge
stabilizer of $T'$ if and only if $X_a$ can be expressed as the coarse intersection of a finite
collection of depth zero vertex spaces of~$X$.

\section{Tree Rigidity}
\label{TreeRigidity}

\subsection{Examples and motivations}

Let $\G$ be a finite type graph of groups, $X \to T$ the Bass-Serre tree of spaces, $G = \pi_1\G$. 
In the Classification Theorem~\ref{TheoremClasses} we saw that, in many cases, any quasi-isometry
of $X$ coarsely respects the tree structure and so gives a quasi-isometry of $T$; this is true, for
example, when $\G$ has finite depth and satisfies the separation hypotheses. Thus, if $H$ is
quasi-isometric to $G$ we have a cobounded quasi-action of $H$ on $T$. We now aim to turn this
quasi-action into an action. This means building a tree of spaces $X' \to T'$ with an $H$ action,
and quasi-isometries $X \to X'$, $T \to T'$ which commute with projections and are $H$-coarsely
equivariant. 

For the moment let us restrict our attention to the following question: given a quasi-action of a
group $H$ on a tree $T$, is it quasiconjugate to an action? This is not always true. In this
section we present some counterexamples, which explain some of the difficulties involved, and help
to further motivate the finite depth hypothesis.

One important source of examples of exotic quasi-actions on trees is isometric actions on
$\R$-trees: 

\begin{lemma}  Every $\R$-tree $\tau$ is quasi-isometric to a simplicial tree $T$, and so every
cobounded action of a group $G$ on $\tau$ is quasi-conjugate to a cobounded quasi-action of $G$ on
$T$.
\end{lemma}

\begin{proof}
Pick a base point $p \in\tau$. Let $B_n = \{x \in \tau \suchthat d(p,x) \le n\}$, and let $V_n =
\bdy B_n = \{x \in \tau\suchthat d(p,x) = n\}$; $V_0 = \{p\}$. The vertex set of the required
simplicial tree $T$ is $V=V_0 \union V_1 \union V_2 \union \cdots$. An edge is constructed from
each point $v \in V_n$ to the unique point in $V_{n-1}$ which is connected to $v$ by a segment in
$\closure(B_n - B_{n-1})$. This produces a simplicial tree $T$ quasi-isometric to $\tau$;
the obvious map $T \to \tau$ is a $(1,2)$-quasi-isometry.
\end{proof}

In an earlier draft of this paper we had posed the question of whether
every quasi-action on a tree arises, as in the above lemma, from an
action on an $\R$-tree. But before we ever distributed this paper, Jason
Manning had anticipated the question and produced counter-examples
\cite{Manning:pseudocharacters} --- he produces quasi-actions on trees of
groups which cannot act non-trivially on $\R$-trees.

\paragraph{Example.} Even in the context of Bass-Serre complexes, there
are quasi-actions on trees not quasi-conjugate to actions.   Let $G$ be
the group $(\Z\oplus\Z)*\Z = (\Z *_\Id) * \Z$. The group $G$ acts on  its
Cayley complex $X$, which consists of pairwise disjoint planes, each 
plane being the Cayley complex for some conjugate of $\Z \oplus \Z$,
such  that each lattice point of such a plane is connected via two
segments to  two other planes. The graph of groups associated to the
splitting $G = (\Z *_\Id) *\Z$ consists of a mapping torus raft
corresponding to  the splitting $\Z \oplus \Z = \Z*_\Id$, together with
another edge labelled $0$ leading from the vertex of the raft to another
vertex labelled $\Z$. The Cayley complex $X$ also serves as the
Bass-Serre tree of spaces, with respect to a projection map $\pi \from X
\to T$ which, on each plane $P$, is a linear projection to a line in $T$
whose point pre-images are lines in $P$ parallel to one of the axes of
$P$.

Now consider a different quotient map on $X$, collapsing each plane $P$ in $X$ to a line $L_P$ by a
linear projection whose point pre-images have irrational slope. Doing this in a $G$-equivariant
manner, the result is an $\R$-tree $\tau$ on which $G$ acts. The tree consists of the pairwise
disjoint lines $L_P$, each with a countable dense set of points at which two segments are
attached, each such segment connecting $L_P$ to other lines. Now choose a lattice of points in the
line $L_P$, and move the attaching point of each segment to a closest lattice point. These bounded
moves produce a quasi-isometry from $\tau$ to a new simplicial tree $T'$. This quasi-isometry has
the effect of quasiconjugating the $G$-action on $\tau$ to a cobounded quasi-action on $T'$. But
this quasi-action is not quasi-conjugate to a simplicial action, because each $\Z\oplus\Z$
subgroup of $G$ would stabilize a line leading to a contradiction as above.

We claim, on the other hand, that there are quasi-isometries $F \from X\to X$ and $f\from T\to T'$
so that the compositions $X \to T\xrightarrow{f} T'$ and $X\xrightarrow{F} X \to T'$ differ by a
bounded amount. We therefore obtain a quasi-action of $G$ on $X$, satisfying vertex rigidity, but
which is not quasi-conjugate to any action on any Bass-Serre tree of spaces for $G$.

To justify the claim, first note that any irrational rotation $\phi_\alpha \from \R^2 \to \R^2$ is
coarsely equivalent to a bijection of the integer lattice $\Z^2$. This is an application of the
Marriage Lemma  (see \cite{HalmosVaughan} or \cite{GraverWatkins}). Divide $\R^2$ into fundamental
domains $\{D_{i,j}\}_{(i,j)\in\Z^2}$ of area $1$ for the action of $\Z^2$, $D_{i,j} = [i,i+1)
\cross [j,j+1)$. The sets $\{\phi_\alpha(D_i) \}_{i \in \Z^2}$ still partition the plane into sets
of area $1$. It follows that, within a fixed distance of $\phi_\alpha(i,j)$, there is at least $1$
and no more than a constant number of lattice points; and similarly, within a fixed distance of
each $(i,j)$ there is at least $1$ and no more than the same constant number of points of
$\phi_\alpha(\Z\oplus\Z)$. The Marriage Lemma now applies, to set up a bijection between
$\Z\oplus\Z$ and $\phi_\alpha(\Z\oplus\Z)$ which moves each point a bounded distance. The
fundamental domain trick used in the above argument is also found in \cite{BuragoKleiner:nets} and
in \cite{McMullen:lipschitz}. 

Applying this argument to each plane of $X$ produces the desired quasi-isometry $F \from X \to
X$.  Thus we see that the no line-like  rafts hypothesis is necessary.

\subsection{Outline of the Tree Rigidity Theorem.} 

We repeat here the statement of the Tree Rigidity Theorem:

\begin{theorem*}[\ref{TheoremTreeRigidity}]
Let $\G$ be a finite type, finite depth, irreducible graph of groups, with Bass-Serre tree of
spaces $X\to T$. Suppose that no depth zero raft of $T$ is a line. Let $H$ be a finitely generated
group which quasi-acts properly and coboundedly on $X$, coarsely preserving vertex and edge
spaces. Then there is a finite type, finite depth, irreducible graph of groups $\G'$, with
Bass-Serre tree of spaces $X' \to T'$, and there is an isomorphism $H\approx\pi_1\G'$, such that the
induced quasi-isometry $f \from X' \to X$ coarsely respects vertex and edge spaces.
\end{theorem*}

From the conclusion that $f$ coarsely respects vertex and edge spaces, one obtains an induced
quasi-isometry $f_\# \from T \to T'$ with the property that $X_{h \cdot a} \ceq{[r]}
X_{f_\#(a)}$, defined by choosing $a'$ as in the theorem and setting $f_\#(a)=a'$.

Over the next several subsections we develop tools for the proof of
Theorem~\ref{TheoremTreeRigidity} by considering special cases.  The general proof is given in
Section~\ref{SectionTreeRigidityProof}, by piecing together and adapting the tools introduced for
the various special cases.

The first special case occurs when the graph of groups $\G$ is homogeneous, meaning that each
edge-to-vertex injection has finite index image. Homogeneity is equivalent to saying that the
Bass-Serre tree $T$ has bounded valence. In this case, the entire tree $T$ is a single depth zero
raft. Since $T$ is not a line, by hypothesis, Proposition~\ref{PropTrichotomy} implies that either
$T$ is a point or $T$ is bushy. When $T$ is a point the conclusion of
Theorem~\ref{TheoremTreeRigidity} follows trivially. The case when $T$ is bushy is handled by the
main results of \cite{MSW:QTOne}:

\begin{theorem}
\label{homogeneous}
If $T$ is a bounded valence, bushy tree, then any quasi-action of a group $H$ on $T$ is
quasi-conjugate to an action of $H$ on another bounded valence, bushy tree $T'$. It follows that
Theorem~\ref{TheoremTreeRigidity} holds for under the special assumption that the graph of group
$\G$ is homogeneous.
\qed\end{theorem}

To motivate later developments, we recall briefly the two major steps of the proof of this theorem
from \cite{MSW:QTOne}. First, although $H$ need not act on $T$, it does have an induced action on the
space of ends of $T$, which is used to construct an action of $H$ on a simply connected
\nb{2}complex $K$ quasiconjugate to the given quasi-action on $T$. Second, one adapts Dunwoody's
theory of tracks to promote the action of $H$ on $K$ to an action on a tree $T'$ semiconjugate to
the action on $K$. 

We will consider two other special cases, designed to separate out two different techniques of
proof. The special case considered in Proposition~\ref{PropNoTracks} occurs when no two edge
spaces are coarsely equivalent; this case illustrates the technique of induction on
depth. And the special case considered in Proposition~\ref{PropDepthOne} occurs when all edge
spaces are of depth one; the proof in this case is close to the proof of Theorem~\ref{homogeneous},
and follows a similar two step outline.

\subsection{Special case: isolated edge spaces}

We first give the proof under the special assumption that the edge spaces are \emph{isolated},
meaning that no two edge spaces are coarsely equivalent. 

\begin{proposition}[The case of isolated edge spaces]
\label{PropNoTracks}
Let $\G$ be a finite type, finite depth, irreducible graph of groups, with Bass-Serre tree of
spaces $X \to T$ and with $G = \pi_1\G$. Assume that edge spaces in $X$ are isolated. If $H$ is a
finitely generated group quasi-acting properly and coboundedly on $X$, coarsely respecting the
vertex and edge spaces, then the conclusions of Tree Rigidity hold.
\end{proposition}

\begin{proof} By hypothesis, no two edge spaces of $X$ are coarsely equivalent. 

Observe now that no edge space of $X$ is coarsely equivalent to a vertex space. To see why,
suppose on the contrary that the edge space $X_{e'}$ is coarsely equivalent to the vertex space
$X_v$. Letting $e$ be the first edge on the path from $v$ to $e'$ it follows that $X_{e}$ is coarsely
equivalent to $X_v$, and so $\Stab(e)$ has finite index in $\Stab(v)$. If the index is $>1$ then we
easily contradict isolated edge spaces, so we may assume $\Stab(e) = \Stab(v)$. By irreducibility it
follows that the opposite end of $e$ is located at a vertex $v'$ in the orbit of $v$. The finite
depth property implies that $\Stab(e)$ has finite index in $v'$. A group element taking $v$ to $v'$
takes $e$ to a distinct edge whose stabilizer is commensurable to $\Stab(e)$, again
contradicting isolated edge spaces.

Observe next that no two vertex spaces of $X$ are coarsely equivalent, for suppose that $X_v$ and
$X_w$ are coarsely equivalent, with $v \ne w$. We have $X_v \ceq{c} X_v \cintersect{c} X_w
\csubset{c} X_e$ for any edge $e$ on the shortest path in $T$ connecting $v$ to $w$, in particular
$X_v \csubset{c} X_e$ where $e$ is the first edge from $v$ to $w$. But we also have $X_e
\csubset{c} X_v$ and so $X_v \ceq{c} X_e$, contradicting the previous observation.

We can summarize these observations as follows: the depth zero subforest $\forest_0$ consists of
the vertex set $\V$; and for $n \ge 1$ each depth $n$ raft contains exactly one depth $n$ edge.
A depth $n$ flotilla can, as usual, contain more than one depth $n$ raft, and so it can contain
more than one depth $n$ edge.

Notation: let $\E_n$ denote the edges of depth $n$, and so $\forest_n$ is a disjoint union
$\forest_n = \V \union \E_1 \union\cdots\union \E_n$.

An immediate consequence of the above observations is that the group $H$ \emph{acts} on the set
$\VE$, preserving each of the sets $\V$ and $\E_n$, $n \ge 1$. To see why, for
each vertex or edge $a \in\VE$ and $h\in H$ the vertex or edge $h \cdot a \in \VE$ is
uniquely characterized by the coarse equation $h \cdot X_a \ceq{c} X_{h \cdot a}$. The existence
of $a'$ such that $h \cdot X_a \ceq{c} X_{a'}$ follows from the hypothesis that $h$ coarsely
respects edge and vertex spaces; and the uniqueness of $a'$ comes from the observations above. The
vertices in $\VE$ are characterized as the maximal elements with respect to coarse inclusion. 
The hypothesis that $H$ coarsely respects vertex and edge spaces contains in it the statement that
depth is preserved, and hence $H$ preserves $\E_n$ for each $n \ge 1$.

It is also immediate that the action of $H$ on $\VE$ agrees uniformly with the quasi-action of $H$
on $X$, in the sense that there is a constant $r \ge 0$ such that $h \cdot X_a \ceq{[r]} X_{h \cdot
a}$ for all $a \in \VE$, $h\in H$.

The subtle part of this lemma is that although $H$ acts on $\VE = \V \union \E$, the
action need not preserve the incidence relation between $\E$ and $\V$. We will give an
example of this phenomenon after the proof. Our goal, therefore, is to redefine the incidence
relation between $\E$ and $\V$ in an $H$-equivariant manner, so that the result is a tree;
moreover, we can do this so that the new tree is $H$-equivariantly quasi-isometric to old tree
$T$, by making sure that the new incidence relation differs from the old by at most a bounded
amount. We define the new incidence relation inductively on $\E_1,\E_2,\ldots$.

For the basis step of the induction, we automatically have an action of $H$ on $\V = \forest_0$.
(In the general proof of Tree Rigidity, the basis step will require invoking the
homogeneous case \cite{MSW:QTOne}).

Next we verify the first step of the induction by proving that the action of $H$ on $\forest_1 =
\V \union \E_1$ automatically preserves the edge--vertex incidence relation. This case exhibits
some features of the general inductive step, though it is atypical in that it misses one important
feature as we'll see later. The key point, as was explained in the proof of the Vertex--Edge
Rigidity Theorem~\ref{TheoremVERigidity}, is that a pair of vertices
$(u,v)$ is the endpoint set of a depth~1 edge $e$ if and only if the coarse intersection $X_u
\cintersect{c} X_v \ceq{c} X_e$ is maximal, with respect to coarse inclusion, among all pairwise
coarse intersections of distinct vertex spaces. The ``only if'' direction is clear, as explained
in the proof of Theorem~\ref{TheoremVERigidity}. The ``if'' direction follows from the hypothesis
that $X$ has isolated edge spaces: if $(u,v)$ is maximal, then the coarse intersection $X_u
\cintersect{c} X_v$ is realized by each edge space along the shortest path connecting $u$ to $v$,
and isolation of edge spaces implies that there can be only one edge on this path. As maximality
of the coarse intersection $X_u \cintersect{c} X_v$ is clearly invariant under the $H$ action on
pairs of vertices, $H$ preserves the incidence relation between $\E_1$ and $V$, and so $H$ acts on
$\forest_1$.

The second step of the induction will be typical of the general inductive step: we describe how to
modify the forest $\forest_2$ relative to $\forest_1$, to obtain a forest $\forest'_2$ on which $H$
acts. First we repeat some of the arguments of the Vertex--Edge Ridigity
Theorem~\ref{TheoremVERigidity}. Consider the quotient tree $T / \forest_1$ obtained by shrinking
each component of $\forest_1$ to a point. The tree $T / \forest_1$ contains the forest $\forest_2
/ \forest_1$. The vertex--edge incidence relation in the forest $\forest_2 / \forest_1$ is
characterized in a quasi-isometrically invariant manner by the property that two components $A,B$
of $\forest_1$ are the endpoint set of an edge $e$ of $\forest_2 / \forest_1$ if and only if the
coarse intersection of $X_A \cintersect{c} X_B \ceq{c} X_e$ in $X$ is maximal, with respect to
coarse inclusion, among all pairwise coarse intersections of distinct components of $\forest_1$. 
This immediately implies that $H$ acts on the forest $\forest_2 / \forest_1$, preserving the
edge--vertex incidence relation. So far this is similar to the first step of the induction.

We can think of $\forest_2$ as obtained from $\forest_2 / \forest_1$ by blowing up each vertex,
replacing it by the corresponding component of $\forest_1$. 

Consider a directed edge $e$ of $\forest_2 / \forest_1$ from a vertex $A$ of $\forest_2 /
\forest_1$ to a vertex $B$. The edge $e$ is identified with a certain depth 2 edge of $\forest_2$
from a vertex $a \in A$ to a vertex $b \in B$. These attachments are made equivariantly with
respect to the $G$ action, but they need not be equivariant with respect to the $H$ action. By
sliding the attaching points of each such edge $e$ we will produce a new forest $\forest'_2$ on
which $H$ acts, and an $H$-equivariant quasi-isometry $\forest_2 \to \forest'_2$, without changing
the map on $\forest_1$ or on $\forest_2 / \forest_1 = \forest'_2 / \forest_1$. There are two
constraints on sliding the attaching points: one constraint, as just described, comes from the
requirement of $H$-equivariance; the other constraint is that we cannot slide the attaching points
very far if we want the resulting map $\forest_2 \to \forest'_2$ to be a quasi-isometry. 

Consider the subgroup $\Stab(e)$ of $H$ which stabilizes $e$ under the action on $\forest_2 /
\forest_1$. The subgroup $\Stab^*(e) \subgroup \Stab(e)$ which preserves orientation has index at
most~2. The action of $\Stab^*(e)$ fixes the tail end $A \in \V(\forest_2 / \forest_1)$ and the
head end $B \in \V(\forest_2 / \forest_1)$ of $e$. The group $\Stab^*(e)$ acts on the tree $A$,
but $\Stab^*(e)$ does not necessarily fix the vertex~$a$, so we must find a new point of
$A$ to which the tail end of $e$ can be attached. The vertex $a$ has bounded orbit under the
action of $\Stab^*(e)$ on $A$, because $\Stab^*(e)$ quasi-acts on $T$ preserving both of the
subtrees $A$ and $B$, and the vertex $a$ has distance~1 from $B$, and so each point in the orbit
$\Stab^*(e) \cdot a$ has uniformly bounded distance from $B$, which implies that $\Stab^*(e) \cdot
a$ has uniformly bounded distance from the point $a$. The convex hull of the orbit $\Stab^*(e)
\cdot a$ is a subtree of bounded diameter in $A$, and hence it contains a fixed point $a'$. Either
$a' \in \V(T)$ or $a'$ is the midpoint of an edge that is inverted by the action of $\Stab^*(e)$
and we subdivide this edge at $a'$. We now attach the tail end of $e$ to the vertex $a'$. 

For each $H$-orbit of oriented edges of $\forest_2/\forest_1$, we choose one oriented edge
representing $e$ this orbit, attach its tail end to the appropriate point $a'$ as described, and
then extend $H$-equivariantly. This defines the forest $\forest'_2$ with an $H$-action and an $H$
almost equivariant quasi-isometry $\forest_2 \to \forest'_2$.

The general inductive step is exactly the same as the second inductive step: assuming
that we have an action of $H$ on a forest $\forest'_i$, and an $H$ almost equivariant map from
$\forest_i$ (appropriately subdivided) to $\forest'_i$, the edges of $\forest_{i+1} / \forest_i$
are attached to $\forest'_i$ in an $H$-equivariant manner, yielding a forest $\forest'_{i+1}$ on
which $H$ acts, and an $H$ almost equivariant quasi-isometry from $\forest_{i+1}$ (appropriately
subdivided) to $\forest'_{i+1}$. The induction stops with the desired action of $H$ on
$T'=\forest'_n$ and $H$-almost equivariant quasi-isometry $T \to T'$. 
\end{proof}

\paragraph{Example.} Here is an example in the context of Proposition~\ref{PropNoTracks} where the
action of $H$ on $\VE(T)$ is not compatible with the incidence relation between edges and vertices
of $T$.

Consider the following graph $\G_1$ of virtually free groups:
$$ 
\xymatrix{
\<a,b,d\> \ar@{-}[r]^(.4){\<a,b\>} & \<a,b,c\> \semidirect \Z/3 \\
 & \<a,e\> \ar@{-}[u]^{\<a\>}
}
$$
where $\Z/3$ acts on $\<a,b,c\>$ by cyclically permuting the generators $a,b,c$.

Define a homomorphism $\pi_1\G_1 \to \Z / 3$, whose restriction to $\<a,b,c\> \semidirect \Z/3$
is the standard projection onto $\Z/3$, and whose restriction to the other vertex groups is
trivial. The kernel of this homomorphism is the fundamental group of the following graph $\G_2$ of
free groups:
$$
\xymatrix{
& \<c,e_3\> \ar@{-}[d]_{\<c\>} & \<c,a,d_3\> \ar@{-}[dl]_{\<c,a\>} \\
\<a,b,d_1\> \ar@{-}[r]^{\<a,b\>} & \<a,b,c\> & \<b,e_2\> \ar@{-}[l]_{\<b\>} \\
\<a,e_1\> \ar@{-}[ur]_{\<a\>} & \<b,c,d_2\> \ar@{-}[u]_{\<b,c\>}
}
$$
The group $\Z/3$ acts on this graph by permuting $a,b,c$, by permuting $d_1,d_2,d_3$, and by
permuting $e_1,e_2,e_3$.

Now slide the $\<a\>$ edge across the $\<a,b\>$ edge to produce the graph of groups
$\G_3$, with fundamental group isomorphic to $\pi_1\G_2$:
$$\xymatrix{
& \<c,e_3\> \ar@{-}[d]_{\<c\>} & \<c,a,d_3\> \ar@{-}[dl]_{\<c,a\>} \\
\<a,b,d_1\> \ar@{-}[r]^{\<a,b\>} & \<a,b,c\> & \<b,e_2\> \ar@{-}[l]_{\<b\>} \\
\<a,e_1\> \ar@{-}[u]_{\<a\>} & \<b,c,d_2\> \ar@{-}[u]_{\<b,c\>}
}
$$
Let $X_i \to T_i$ be the tree of spaces for $\G_i$, $i=1,2,3$. The trees of spaces $X_1 \to T_1$,
$X_2 \to T_2$ are identical. The presentation of $\G_2$ and $\G_3$ makes it clear that edge spaces
are isolated. The group $\pi_1\G_1$ acts on $X_1$ and on $X_2$, strictly
respecting vertex and edge spaces. There is a vertex and edge space respecting quasi-isometry $X_2
\to X_3$, which quasiconjugates the action of $\pi_1\G_1$ on $X_2$ to a quasi-action on $X_3$,
coarsely respecting vertex and edge spaces, but this quasi-action does \emph{not} respect the
incidence relation between edges and vertices.

\subsection{Special case: all edges have depth one}

The general proof of Tree Rigidity will follow the same sort of induction on edge depth used in the
previous special case. But the proof is complicated by the fact that in general a quasi-action that
coarsely respects vertex and edge spaces need not strictly respect the maximal edges, that is, the
quasi-action of $H$ on $\VE(T)$ need not be compatible with any true action on $\VE(T)$. This occurs,
for example, in the geometrically homogeneous case, covered by Theorem~\ref{homogeneous}, in which
the entire tree $T$ has depth zero. In order to highlight these issues without the distractions of
the induction, we shall focus on the special case where every edge has depth one; later, the general
induction step will follow the proof of this special case. 

Note that if we combine the two special cases, so that every edge space is isolated and of depth
one, Proposition~\ref{PropNoTracks} applies to show that the $H$ quasi-action on $T$ is actually an
action, and we are done. In general, the new special case is more difficult. The proof in this case
involves all of the extra ingredients that will be needed later: the induction step in
the general case will be proved by invoking Proposition~\ref{PropDepthOne}.

\begin{proposition}[Depth one edges]
\label{PropDepthOne}
Let $\G$ be finite type, irreducible graph of groups, with $G = \pi_1\G$, and with
Bass-Serre tree of spaces $X \to T$. Suppose that every edge of $T$ has depth one. Let the group
$H$ quasi-act properly and coboundedly on $X$, coarsely respecting the edge and vertex spaces.
Then there is a finite type, irreducible, graph of groups $\G'$ in which every edge has depth
one, with fundamental group $H$ and Bass-Serre tree of spaces $X' \to T'$, and there is an
$H$-coarsely equivariant quasi-isometry $X\to X'$ that coarsely respects vertex and edge spaces.
\end{proposition}

\begin{proof} Note that since every edge has depth one, every vertex is either fat or of depth
one: a nonfat, depth zero vertex $v$ of $T$ would be incident to an edge $e$ with $X_v \ceq{c}
X_e$ and so $e$ would have depth zero, a contradiction. Note that a vertex of $T$ is fat if and
only if it has infinite valence. The vertex set $\V(T)$ is a disjoint union $\V(T) =\V^\infinity(T)
\union \V^\fin(T)$ where $\V^\infinity$ denotes infinite valence and $\V^\fin$ finite valence. As
usual we let $\E(T)$ denote the edge set, and we have a disjoint union $\VE(T) =
\V(T)\union\E(T)$.

Since the quasi-action of $H$ on $X$ coarsely respects vertex and edge spaces, and since any two
fat vertex spaces have infinite Hausdorff distance, it follows that $H$ \emph{acts} on the
collection of fat vertex spaces $\V^\infinity(T)$; this action is characterized by the property
that for each $v \in \V^\infinity(T)$ and $h \in H$ we have $h \cdot X_v \ceq{[C]} X_{h \cdot v}$,
for a constant $C \ge 0$ independent of $v$ and $h$. We do not, however, have an action of $H$ on
the collection of nonfat vertex spaces, nor on the collection of edge spaces. 

Fix a base point $p \in X$, and let $\O \from H \to X$ be the induced orbit map, a quasi-isometry,
defined by $\O(h) = h \cdot p$. A judicious choice of $p$ will be made later.

We now describe the proof of Proposition~\ref{PropDepthOne} in three steps, the first two of which
parallel the two major steps in the proof of Theorem~\ref{homogeneous}.

\paragraph{Step 1: An action of $H$ on a simply connected 2-complex.} In order to initiate the
ultimate goal of making $H$ act on a tree, we start by producing an $H$ action on a certain simply
connected, \nb{2}dimensional cell complex $K$, and an $H$-equivariant quasi-isometry $\pi
\from K \to T$, so that certain finiteness properties are satisfied.

To state the result of step~1, a map $\pi \from K \to T$ from a \nb{2}dimensional cell complex $K$
to a tree $T$ is \emph{tight} if $\pi$ maps each vertex of $K$ to a vertex of $T$, and $\pi$ maps each
edge $e \subset K$ with endpoints $u,v$ to the convex hull of $f(u),f(v)$ in $T$, either collapsing
to a point or mapping homeomorphically to an arc of $T$. We make no restrictions on how the map
acts on \nb{2}cells, but see Figure~\ref{FigureSimplexMap} below for the special case of a
\nb{2}simplex.

The main work of Step~1 is:

\begin{lemma}[An action on a simply connected 2-complex]
\label{Lemma2complex} Let $\G$, $G$, $X \to T$, and $H$ be as in
Proposition~\ref{PropDepthOne}. There exists a simply connected \nb{2}dimensional cell complex
$K$, a cocompact action of $H$ on $K$, and an $H$-almost equivariant, tight quasi-isometry
$\pi\from K \to T$, such that the following properties hold:
\begin{enumerate}
\item \label{ItemEquivBij}
The map $\pi$ restricts to an $H$-equivariant bijection $K^0 \to \V^\infinity(T)$. 
\item \label{ItemEdgeRaft} For each edge $e$ of $K$, the image $\pi(e)$ is contained in a
depth one raft of $T$.
\item \label{ItemEdgeFinite}
The complex $K$ is uniformly locally finite along edge interiors.
\item \label{ItemTrackFinite}
The cardinality of the sets $\pi^\inv(x) \intersect K^1$, for $x \in T - T^0$, is uniformly
finite.
\end{enumerate}
\end{lemma}

The full proof is given in Section~\ref{Section2complex} below; here we give a sketch.

Property (\ref{ItemEquivBij}) dictates the construction of the \nb{0}skeleton: take $K^0 =
\V^\infinity(T)$ equipped with the given $H$-action. To construct the \nb{1}skeleton so as to
satisfy property (\ref{ItemEdgeRaft}), use the fact that $T$ is the union of its depth one rafts,
together with the fact that the set $\V^\infinity(T)$ is coarsely dense in~$T$, to produce a
collection of edges attached to $K^0$ which form a connected graph, so that the tightened images of
these edges in $T$ have uniformly bounded length, and so that each image is contained in some
depth one raft. The set of vertex pairs of $K^0$ to which edges are so far attached need not be
$H$-equivariant, but we can enlarge this set to make it $H$-equivariant, so as to extend the
$H$-action on $K^0$ to an $H$-action on the \nb{1}skeleton $K^1$. Faces are attached as needed, in
an $H$-equivariant fashion, to kill loops in the \nb{1}skeleton, producing the simply connected
\nb{2}complex $K$. The construction of $K$ easily yields the $H$-almost equivariant, tight
quasi-isometry $f \from K \to T$. The finiteness properties (\ref{ItemEdgeFinite}) and
(\ref{ItemTrackFinite}) are a consequence of the assumption that edges of $T$ have depth one. 

The main difference with the homogeneous case, Theorem~\ref{homogeneous}, is that in the current
situation the vertex set of $K$, with its $H$-action, is already present as $\V(T)$. In the
homogeneous case the vertices of $K$ were constructed by a completely different process couched in
terms of the action of $H$ on the space of ends of $T$. Also, in the homogeneous case
$K$ is locally finite at vertices, which eases the construction of edges and faces.

\bigskip

As a consequence of Lemma~\ref{Lemma2complex}, we can compare the cell stabilizers of the action
of $H$ on $K$ to the vertex and edge spaces of $X$:

\begin{corollary} 
\label{CorollaryYCells}
For each $v \in K^0$ we have $\Stab(v) = \Stab(\pi(v))$. Moreover:
\begin{enumerate}
\item For each $v \in K^0$ we have $\O(\Stab(v)) \ceq{c} X_{\pi v}$.
\item For each edge $e\subset K$ and each edge $e' \subset T$ contained in $\pi e$ we have
$\O(\Stab(e)) \ceq{c} X_{e'}$.
\item For each \nb{2}cell $\sigma \subset K$ and each edge $e' \subset T$ contained in
$\pi(\bdy\sigma)$ we have $\O(\Stab(\sigma)) \ceq{c} X_{e'}$.
\end{enumerate}
\end{corollary}

The coarse equivalences in this corollary, between stabilizer groups and vertex and edge spaces,
are \emph{not} uniform. We will overcome this lack of uniformity later.

\begin{proof} The first sentence follows because $\pi \from K^0 \to \V^\infinity(T)$
is an $H$-equivariant bijection; we identify these two sets under this bijection. Conclusion (1)
follows from the Coboundedness Principle~\ref{PropCoboundednessPrinciple} applied to the
quasi-action of $H$ on $X$ which coarsely respects the pattern of fat vertex spaces of $X$.

Consider an edge $e \subset K$, with $\bdy e = \{v',w'\} \subset K^0$, and so $\bdy\pi(e) = \{\pi
v',\pi w'\} \in \V^\infinity(T)$. There is a subgroup $\Stab^*(e) \subgroup \Stab(e)$ of index at
most~2 that does not invert $e$, and we have $\Stab(e) \ceq{c} \Stab^*(e) = \Stab(v')
\intersect\Stab(w') \ceq{c}\Stab(v') \cintersect{c} \Stab(w')$ where the last coarse equivalence
comes from Lemma~\ref{LemmaCoarseIntersectionSubgroup}. It follows that $\O(\Stab(e)) \ceq{c}
\O(\Stab(v'))\cintersect{c} \O(\Stab(w')) \ceq{c} X_{\pi v'} \cintersect{c} X_{\pi w'}$, by
applying (1). But since $\pi(e)$ is contained in a depth one raft of $T$, all the edges along
$\pi(e)$ have coarsely equivalent edge spaces in $X$, and so these edge spaces are all coarsely
equivalent to the coarse intersection $X_{\pi v'} \cintersect{c} X_{\pi w'}$. If $e'$ is any such
edge we therefore have $X_{e'} \ceq{c} \O(\Stab(e))$.

Finally, consider a \nb{2}cell $\sigma \subset K$, an edge $e \subset \bdy\sigma$, and an edge $e'
\subset T$ contained in $\pi(e)$. Since $K$ is locally finite along edges, the subgroups
$\Stab(\sigma)$ and $\Stab(e)$ are commensurable, and so $\Stab(\sigma) \ceq{c} \Stab(e)$. It
follows that $\O(\Stab(\sigma)) \ceq{c} \O(\Stab(e)) \ceq{c} X_{e'}$, by (2).
\end{proof}

\paragraph{Step 2: An action of $H$ on a tree.} We produce an $H$ action on a tree~$T$
from the $H$-action on $K$ by using Dunwoody's tracks \cite{Dunwoody:Accessible}. In the
homogeneous case \cite{MSW:QTOne} we were able to directly quote a result of Dunwoody, but in the
present case, the lack of local finiteness at vertices of $K$ requires a few new ideas; similar ideas
have been used in \cite{Niblo:ends}. 

First, choose an $H$-equivariant triangulation of $K$ without edge inversions, producing a
simplicial \nb{2}complex denoted $L$. The $H$-equivariant map $K \to T$ may be redefined, on the
new vertices, the subdivided edges, and the subdivided \nb{2}cells, to produce a tight,
$H$-equivariant map $\pi\from L \to T$. We also impose a tightness condition on the restriction of
$\pi$ to each \nb{2}simplex, depicted in Figure~\ref{FigureSimplexMap}.

\begin{figure}
\centeredepsfbox{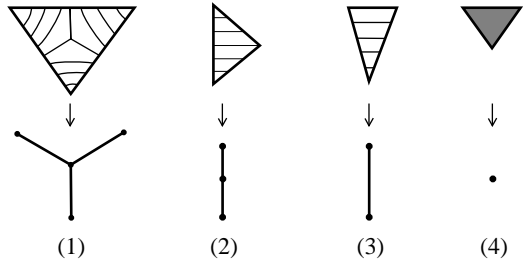}
\caption{A tight map of a simplex to a tree: (1) the vertices map to the endpoints of a triod; (2)
the vertices map to three distinct points on an arc; (3) the vertices map to two points on an arc;
(4) the vertices map to a single point.}
\label{FigureSimplexMap}
\end{figure}

A \emph{track pattern} in a simplicial \nb{2}complex $L$ is a \nb{1}dimensional complex
$P$ embedded in $L$, such that for each \nb{2}simplex $\sigma$ in $L$, $P \intersect
\sigma$ is a disjoint union of finitely many arcs, each component of which connects
points in the interiors of two distinct edges of $\sigma$. For each edge $e$ of $L$ and
each $x\in P\intersect e$, each \nb{2}simplex $\sigma$ incident to $e$ therefore contains
a component of $t \intersect \sigma$ incident to $x$.  A \emph{track} $t$ is a connected
track pattern. Given a track pattern $P$, a subset of $P$ is a track if and only if it is a
component of $P$. Two tracks $t,t'$ are \emph{isotopic} if there is an ambient isotopy of $L$
preserving each cell and taking $t$ to $t'$.

Suppose that $L$ is simply connected. Given a finite track $t$ in $L$, the set $L-t$ has two
components, and $t$ has a regular neighborhood which is a trivial $I$-bundle over~$t$; if each of
the components of $L-t$ is unbounded then we say that $t$ is an \emph{essential} track in
$L$. More generally, a finite track pattern $P$ in $L$ has finitely many components, and $P$ has an
essential component if and only if there is an \emph{essential decomposition} $L - P = U \union V$
where $U,V$ are disjoint, nonempty, and unbounded, and each of $U,V$ is a union of components of
$L-P$. This fact is proved by induction on the number of components of $P$: if $t$ is an
inessential component of $P$ then the bounded component $W$ of $L-t$ can be conglomerated with
whichever of $U$ or $V$ lies opposite $t$ from $W$, showing that $L-(P-t)$ has an essential
decomposition. 

For example, if $f \from L \to T$ is a tight map from a simplicial 2-complex to a tree,
then for each $x \in T - T^0$ the set $f^\inv(x)$ is a track pattern in
$L$; see Figure \ref{FigureSimplexMap}. Moreover, if $x$ separates $T$ into two unbounded
components, and if $f^\inv(x)$ is finite, then $f^\inv(x)$ has an essential component.

\begin{theorem}[Tracks Theorem] Let $H$ be a group, $L$ a simply connected simplicial
$2$-complex on which $H$ has a cocompact action, $T$ a simplicial tree with no valence
$1$ vertices, and $\pi\from L\to T$ a tight quasi-isometry. Suppose that $L$ is uniformly
locally finite along edges, and that the tracks $\pi^\inv(x)$, $x \in T-T^0$ are uniformly
finite. Then there is an $H$-equivariant track pattern $P$ in $L$ such that the
components of $P$ are finite and pairwise nonisotopic, and each component of $L-P$ is a
bounded \nb{2}complex (though not necessarily finite).
\label{TheoremTracks}
\end{theorem}

Note that conclusions~(\ref{ItemEdgeFinite}) and~(\ref{ItemTrackFinite}) of
Lemma~\ref{Lemma2complex} exactly match the finiteness conditions in the hypotheses of this
theorem, so the theorem applies.

From the track pattern $P\subset L$ we construct a dual tree $T'$, with one vertex $v_C$ for
each component $C$ of $L-P$, and one edge $e_t$ for each track~$t$ in~$P$. The action of $H$ on
$L$ and $P$ induces an action on $T'$. There is an $H$-equivariant map $L \to T'$ defined as
follows: for each track $t$ in $P$ the map $\pi' \from L \to T'$ takes a regular neighborhood
of $t$ to the corresponding edge $e_t$; for each component $C$ of $L-P$ the map $\pi'$ takes $C$
minus the regular neighborhoods of incident tracks to the corresponding vertex $v_C$. 

We claim that the map $\pi' \from L \to T'$ is a quasi-isometry. This immediately follows from the
fact that the components of $P$ and of $L-P$ are bounded, and that there is a uniformly
finite number of components of $P$ intersecting each edge of $L$, which is a consequence
of cocompactness of the $H$-action on~$L$.

From the quasi-isometries $\pi' \from L \to T'$, $\pi \from L \to T$, the first $H$-equivariant and
the second $H$-almost equivariant, we obtain an $H$-almost equivariant quasi-isometry $T'
\to T$. In Step~3 we shall move this quasi-isometry a bounded amount to obtain an explicit
$H$-almost equivariant quasi-isometry $\rho \from\VE(T') \to \VE(T)$ satisfying some important
properties.

The action of $H$ on $T'$ is cocompact. The action may have edge inversions; if so, subdivide each
inverted edge at its barycenter to obtain an action without inversions. This subdivision does not
change the commensurability classes of vertex and edge stabilizers. We may therefore assume that
$H$ has no edge inversions, and so the quotient $T' / H = \G'$ is a graph of groups with
$H=\pi_1\G'$. Note that we have not yet been careful to insure irreducibility of~$\G'$.

Let $X' \to T'$ be a tree of spaces for $\G'$. The group $H$ acts, properly discontinuously and
cocompactly, on $X'$. As usual, the $H$ quasi-action on $X$ and the $H$ action on $X'$ induce an
$H$-almost equivariant quasi-isometry $f \from X \to X'$. 

We now drop the simplicial triangulation $L$, which was needed only to apply the Tracks
Theorem~\ref{TheoremTracks}, and revert to the original cell structure on~$K$. The ambient
topological space and the set of locally infinite vertices is unchanged---indeed, every vertex of
$K$ is locally infinite. The construction of the dual tree $T'$ from the track pattern $P \subset
K$, which depends only on topology and not on cell structure, is also unchanged.

\paragraph{Step 3: Coarse respect of vertex and edge spaces.} It remains to make $\G'$ irreducible,
prove that the result has finite depth, and prove that $f \from X \to X'$ coarsely respects vertex
and edge spaces. The main work is to analyze the subgroups of $H$ that stabilize the
vertices and edges of~$T'$, compare them to the subgroups stabilizing cells of $K$, and use
Corollary~\ref{CorollaryYCells} to compare them with the vertex and edge spaces of $X$.

We decompose $\VE(T') = \V(T') \union \E(T') = \V^\infinity(T') \union \V^\fin(T') \union \E(T')$,
where $\V^\infinity$ denotes vertices of infinite valence and $\V^\fin$ denotes finite valence. We
shall give a more explicit construction of the $H$-almost equivariant quasi-isometry $\rho \from
\VE(T') \to \VE(T)$, as the union of a bijection $\V^\infinity(T') \to \V^\infinity(T)$ a map
$\V^\fin(T') \to \E(T)$, and a map $\E(T') \to \E(T)$. It will be evident from this construction
that $\rho$ agrees, within uniformly finite distance, with $\pi \composed \pi^\inv$. Moreover, we
shall obtain (nonuniform) coarse equivalences 
$$\O(\Stab(a)) \ceq{c} X_{\rho a} \quad\text{for all}\quad a \in \VE(T')
$$

Consider an infinite valence vertex $v_C \in \V^\infinity(T')$. The corresponding
component $C$ of $K-P$ has a frontier in $K$ consisting of infinitely many tracks. The only way
this can happen is if $C$ contains at least one vertex of $K$, because $K$ is uniformly
locally finite along edges and $C$ is a bounded 2-complex. The question arises whether $C$ can
contain more than one vertex of $K$, and the answer is no:

\begin{lemma} 
\label{LemmaOneFatVertex}
Every component $C$ of $K-P$ contains at most one vertex of $K$. 
\end{lemma}

\begin{proof} We shall identify each vertex $v \in K$ with its image $\pi v \in \V^\infinity(T)$.

Suppose $C$ contains $2$ or more vertices $v_0,v_1$ of $K$. 

We claim that $C$ contains an infinite sequence of vertices $v_0,v_1,v_2,\ldots$ whose images
under $\pi$ lie, in order, along some ray in the tree $T$, diverging to the
end of that ray. This contradicts the fact that the image of the bounded subset $C \subset K$
under the quasi-isometry $\pi \from K \to T$ is a bounded subset of~$T$.

The claim is proved by induction. Suppose we have obtained a sequence $v_0,\ldots,v_k$ of
vertices of $K$, each contained in $C$, whose images under $\pi$ lie in
order along a segment of $T$. Let $e$ be the first edge of $T$ from $v_k$ to $v_{k-1}$,
and let $T_e$ be the component of $T-e$ containing $v_{k-1}$. We shall show that $\Stab(v_k)$, a
subgroup of $\Stab(C)$, contains an element $h$ such that $h \cdot v_{k-1} \not\in T_e$, for then
$v_{k+1} = h \cdot v_{k-1}$ lies in $C$ and the sequence $v_0,\ldots,v_{k+1}$ projects to a
sequence in order along a segment of $T$.

To construct the desired $h \in \Stab(v_k)$, note that the quasi-action of $\Stab(v_k)$ on $T$
keeps
$v_{k-1}$ within a bounded neighborhood of $v_k$, and so there is a constant $A$ with the
following properties for each $h \in \Stab(v_k)$:
\begin{itemize}
\item The intersection $X_{v_k} \cintersect{[A]} X_{h\cdot v_{k-1}}$ realizes the coarse
intersection of $X_{v_k}$ and $h \cdot X_{v_{k-1}}$. 
\item Letting $e_h$ denote the first edge of $T$
from $v_k$ to $h \cdot v_{k-1}$, we have $X_{v_k} \cintersect{[A]} X_{h\cdot v_{k-1}} 
\csubset{[A]} X_{e_h}$.
\end{itemize}
But the quasi-action of $\Stab(v_k)$ on $X_{v_k}$ is cobounded, and the coarse containment of $X_e$
in the vertex space $X_{v_k}$ is a strict coarse containment, and so there exists $h \in
\Stab(v_k)$ such that $X_{v_k} \cintersect{[A]} X_{h\cdot v_{k-1}} \not\csubset{[A]} X_e$, implying
that $h \cdot v_{k-1} \not\in T_e$. 
\end{proof}

As a consequence of Lemma~\ref{LemmaOneFatVertex}, the map $\pi' \from K \to T'$ can be chosen to
restrict to an $H$-equivariant bijection between $K^0$ and $\V^\infinity(T')$. Thus, $\pi'{}^\inv
\composed \pi$ restricts to an $H$-equivariant bijection between $\V^\infinity(T')$ and
$\V^\infinity(T)$, which we take as the definition of the $\rho\restrict
\V^\infinity(T')$. For each $v \in \V^\infinity(T')$ we clearly have $\Stab(v) = \Stab(\rho v)$,
and so $\O(\Stab(v)) \ceq{c} X_{\rho v}$.

Now we turn to stabilizers of edges and of finite valence vertices of $T'$. 

For any edge $e$ of $K$, a subgroup $\Stab^*(e) \subgroup \Stab(e)$ of index at most two acts
without inversion on $e$, and so $\Stab^*(e)$ preserves each point of the finite set $P\intersect
e$, each track of $P$ that intersects $e$, each component of $e-P$, and each component of $K-P$
that intersects $e$.

We define $\rho \restrict \V^\fin(T')$, with values in $\E(T)$, as follows. Consider a vertex $v_C
\in \V^\fin(T')$ corresponding to a component $C$ of $K-P$ that contains no vertex of $K$. We have
$\Stab(v_C) = \Stab(C)$. The \nb{2}complex $C$ is locally finite, and $C$ is bounded by
Theorem~\ref{TheoremTracks}, so $C$ is finite. Choose an edge $e'$ of $K$ that intersects $C$;
such an edge exists, for otherwise $C$ is contained in the interior of a \nb{2}cell of $K$ and so
$\bdy C$ is an inessential track in $P$, a contradiction. As noted above $\Stab^*(e')$ preserves
$C$, and since there are only finitely many edges of $K$ intersecting $C$ it follows that
$\Stab^*(e')$ has finite index in $\Stab(C) = \Stab(v_C)$, and so $\Stab(e')$ is commensurable to
$\Stab(v_C)$. Choosing $e$ to be any edge of $T$ such that $f(e') \supset e$, we define
$\rho(v_C)=e$. We have $\Stab(v_C) \ceq{c}
\Stab(e')$, and so $\O(\Stab(v_C)) \ceq{c} \O(\Stab(e')) \ceq{c} X_{e} = X_{\rho(v_C)}$, where the
second coarse equation follows from Corollary~\ref{CorollaryYCells}. 

We define $\rho \restrict \E(T')$, with values in $\E(T)$, in a similar manner: for any edge $e_t$
of $T'$ corresponding to a component $t$ of $P$, choose an edge $e'$ of $K$ that intersects $t$,
choose an edge $e$ of $T$ contained in $f(e')$, and define $\rho(e_t) = e$. The subgroups
$\Stab(e_t)$ and $\Stab(e)$ are commensurable, and as before using Corollary~\ref{CorollaryYCells}
we have $\Stab(e_t)\ceq{c}\Stab(\rho e_t)$ and so $\O(\Stab(e_t)) \ceq{c} X_{\rho e_t}$. 

This completes the description of the quasi-isometry $\rho \from \VE(T') \to \VE(T)$. For each $a
\in \VE(T')$ we have proved $\O(\Stab(a)) \ceq{c} X_{\rho a}$. Also, for each $b\in \VE(T)$ there
exists $a \in \VE(T')$ such that $X_b \ceq{c} X_{\rho a}$. When $b \in \V^\infinity(T)$ this
follows from the fact that $\rho \from \V^\infinity(T')\to \V^\infinity(T)$ is a bijection. When
$b$ is an edge this follows from the fact that the map $\pi \from K \to T$ contains an edge $e$
such that $b \subset \pi(e)$, and either $e$ intersects a component $t$ of $P$ corresponding to an
edge $e_t \in \E(T')$ and we take $a=e_t$, or the endpoints of $e$ lie in some component $C$ of
$K-P$ corresponding to a vertex $v_C \in \V^\fin(T')$ and we take $a=v_C$.

We can now prove that $\G'$ has depth one, as follows. The coarse inclusion lattice of vertex and
edge spaces of $X'$ is isomorphic to the coarse inclusion lattice of subgroups of $H$ stabilizing
vertices and edges of $T'$, which as we have just shown is isomorphic to the coarse inclusion
lattice of vertex and edge spaces of $X$. Since the last lattice has depth one, so does the first
lattice, and so $\G'$ has depth one.

Also, the complete reduction process applied to $\G'$ does not eliminate any coarse equivalence
class of vertex and edge spaces, because every edge space in $X$ is the coarse intersection of two
depth zero vertex spaces, and so the same is true in $X'$. Thus, we may assume that $\G'$ is
irreducible.

The last task is to prove that the quasi-isometry $f \from X \to X'$ coarsely respects vertex and
edge spaces. We already know that the coarse inclusion lattices of vertex and edge spaces in $X$
and in $X'$ are isomorphic, indeed the map $\rho \from \VE(T') \to \VE(T)$ provides the desired
isomorphism, but we need to promote this to actual coarse respect with uniform constants, and our
arguments so far have been rather sloppy about constants. But we can get around this problem as
follows.

Knowing that $\G$ and $\G'$ are irreducible, the Vertex--Edge Rigidity
Theorem~\ref{TheoremVERigidity} applies, and it remains only to check that $f$ coarsely respects
depth zero vertex spaces. We do this in the general setting of the Coboundedness
Principle~\ref{PropCoboundednessPrinciple} as follows. 

The set of depth zero vertex spaces of $X$ is a pattern in the sense of
Section~\ref{SectionCobounded}, and similarly for $X'$. Both of these patterns satisfy the local
finiteness condition (a) and the coarse discreteness condition (b) of
Section~\ref{SectionCobounded}. The quasi-actions of $H$ on $X$ and on $X'$ satisfy the coarse
action condition (c) with respect to the patterns of depth zero vertex spaces. The
action condition (d) is therefore also satisfied. We already have a bijection between the depth
zero vertex spaces of $X$ and of $X'$, induced by the bijection $\rho \from \V^\infinity(T') \to
\V(T)$, with the property that $\Stab(\rho(v')) = \Stab(v')$ for each $v' \in \V^\infinity(T')$. 
We are therefore reduced to proving the following:

\begin{proposition}
\label{PropFromTreeToPattern}
Let $X,X'$ be metric spaces, let $\A,\A'$ be patterns in $X,X'$ satisfying the local finiteness
condition (a), and let the group $H$ quasi-act properly and coboundedly on $X$ and on $X'$,
satisfying the action condition (d) with respect to the patterns $\A$ and $\A'$. Let $f \from X
\to X'$ be a quasi-isometry that is coarsely $H$-equivariant, and suppose that there exists a
bijection $f_\# \from \A \to \A'$ such that for each $A \in \A$ we have $\Stab(A) =
\Stab(f_\#(A))$. Then $f$ coarsely respects the patterns $\A,\A'$, that is, for all $A \in \A$ we
have $f(A) \ceq{[C]} f_\#(A)$, where the constant $C$ depends only on $X$, $X'$, $\A$, $\A'$, the
quasi-isometry constants of $f$, and the quasi-action constants of $H$.
\end{proposition}

\begin{proof} Choose base points $p \in X$, $p' \in X'$. Consider $A \in \A$ and $A'=f_\#(A) \in
\A'$. The group $S = \Stab(A) = \Stab(A')$ quasi-acts properly and coboundedly on $A$ and on $A'$,
by the Coboundedness Principle~\ref{PropCoboundednessPrinciple}, and it follows that $A \ceq{c} S
\cdot p$ and $A' \ceq{c} S \cdot p'$. Clearly $f(S \cdot p) \ceq{c} S \cdot p'$, and putting it
all together we see that $f(A) \ceq{c} A'$. 

We still need to establish uniformity of this coarse equivalence, that is, we need to show that
$f(A) \ceq{[C]} A'$ where the constant $C$ is independent of $A \in \A$. By the Coboundedness
Principle~\ref{PropCoboundednessPrinciple}, there are finitely many orbits of the action of $H$ on
$\A$ and on $\A'$. Choosing a representative $A \in \A$ of each orbit and the corresponding
orbit representative $A' \in \A'$, we obtain a uniform constant $C$ over the choice of
representatives. By the action condition (d), the quasi-action of $H$ on $X$ coarsely respects the
pattern $\A$, and the quasi-action on $X'$ coarsely respects $A'$; applying these quasi-actions,
and enlarging $C$ as needed we obtain a uniform constant $C$ independent of $A \in\A$. Moreover, if
the orbit representatives of $\A$ are chosen to intersect a ball around $p$ whose size depends
only on the quasi-action constants, and similarly for $\A'$, the required dependencies for the
constant $C$ follow.
\end{proof}

This finishes the proof of Proposition~\ref{PropDepthOne}, modulo more detailed proofs of
Lemma~\ref{Lemma2complex} and the Tracks Theorem~\ref{TheoremTracks}.
\end{proof}

\subsubsection{Proof of Lemma \ref{Lemma2complex}: an action on a 2-complex.}
\label{Section2complex}

We shall use the abbreviated terminology ``raft of $T$'' to mean a depth one raft of $T$. This
cannot cause confusion, because a depth zero raft of $T$ is the same thing as a fat vertex.

We first construct $K^0$, then $K^1$, and then $K^2$, at each stage constructing the $H$-action
and the $H$-almost equivariant quasi-isometry $\pi$. Then we will establish the necessary
finiteness properties. 

The zero skeleton $K^0$ is identified with the set of fat vertices $V^\infinity(T)$. As is evident
by now, there is a well-defined action of $H$ on the set of fat vertices of $T$, characterized by
$h\cdot X_v \ceq{c} X_{h \cdot v}$, and so we obtain the $H$-action on $K^0$ and the
$H$-equivariant map $\pi \from K^0 \to V^\infinity(T)$. As usual, we will blur the distinction
between $v \in K^0$ and $\pi v \in V^\infinity(T)$.

Since $G=\pi_1\G$ acts cocompactly on the tree $T$, the set of fat vertices in $T$ is coarsely
dense, meaning that each point of $T$ is a uniformly bounded distance from a fat vertex.
As a consequence, there exists $R$ and $C$ such that, for any two fat vertices $u, v$ of
$T$, there is a sequence $u=x_{0},x_{1},\ldots,x_{n}=v$ of fat vertices for which
$d(x_{i},x_{i+1})\leq R$ and with $n$ bounded by $C\, d(u,v)$.  If the geodesic segment
from $x_{i}$ to $x_{i+1}$ contains any fat vertices, the sequence can be refined. Thus
the following is natural:

Let $K^1_\#$ be the graph with vertex set the fat vertices of $T$, and edges
connecting pairs $\{u,v\}$ for which:
\begin{itemize}
\item[(i)] $d(u,v)\leq R$
\item[(ii)] The geodesic segment in $T$ between $u$ and $v$ has no interior fat vertices
(hereafter called a fat-free path).
\end{itemize}
The graph $K^1_\#$ has a canonical tight map $\pi \from K^1_\# \to T$, which by the 
above discussion is a quasi-isometry. By properties (i), (ii) above, $\pi(e)$ is contained in a
raft of $T$ for each edge $e \subset K^1_\#$, indeed more is true: there exists a constant $C$
such that:
\begin{itemize}
\item[(iii)] For any edge $e$ of $K^1_\#$ with endpoints $u,v \in K^1_\#$, and for any
edge $e'$ in $\pi(e)$, we have $(X_u \cintersect{[C]} X_v) \ceq{[C]} X_{e'}$. 
\end{itemize}

While $H$ acts on the vertices of $K^1_\#$, the action does not extend to the edges. We fix this
as follows. Given two vertices $u,v$ of $K^1_\#$, the property that $u,v$ are contained in the
same raft of $T$ is $H$-equivariant, and this property holds whenever $u,v$ are endpoints
of an edge of $K^1_\#$. If $u,v \in K^0$ are connected by an edge in $K^1_\#$, if $h \in H$, and if
$h\cdot u$, $h\cdot v$ are not connected by an edge in $K^1_\#$, then connect them, and map this
edge homeomorphically onto the unique path in $T$ connecting $\pi(h \cdot u)$ to $\pi(h \cdot v)$;
as just explained, this path is contained in a raft of $T$. We have now completed the
construction of the \nb{1}skeleton $K^1$, of the $H$ action on $K^1$, and of the $H$-almost
equivariant quasi-isometry $\pi \from K^1\to T$ taking each edge of $K$ to a path contained in
a single raft of $T$. 


For edges of $K^1$, property (iii) above continues to hold with a larger value of $C$,
using the fact that the $H$-quasi-action on $X$ respects fat vertex spaces. For each edge $e$ of
$K^1$ let $X_e$ denote the union of the edge spaces $X_{e'}$ over all $e' \subset \pi e$, and with
$\bdy e = \{u,v\}$ we have $(X_u\cintersect{[C]} X_v) \ceq{[C]} X_e$ and we also have $X_{e'}
\ceq{[2C]} X_e$.

An embedded loop in $K^1$ is \emph{minimal} if any two points on the loop are connected by a
segment of the loop which is geodesic in $K^1$. To make $K^1$ simply connected we simply attach a
2-cell along  every minimal loop. The resulting \nb{2}complex is $K$. The minimal loops are
$H$-equivariant and so $H$ acts on $K$. Since $K^1$ is quasi-isometric to $T$, there is an upper
bound on the lengths of minimal loops.  Also, the 2-cells do not change the quasi-isometry type,
so the quasi-isometry $\pi \from K^1 \to T$ has an extension to an $H$-almost equivariant
quasi-isometry $K \to T$, also denoted $\pi$.

The next lemma plays a key role in establishing the finiteness properties of $K$. Roughly speaking
it says that the boundary of each \nb{2}cell of $\sigma$ is mapped, via $\pi$, into some raft of
$T$.

\begin{lemma} 
\label{LemmaWhateveritwas}
There is an $R \ge 0$ so that for every 2-cell $\sigma$ in $K_{0}$ and any
edges $e_1,e_2 \subset \bdy\sigma$, we have $X_{e_1} \ceq{[R]} X_{e_2}$.
\end{lemma}

\begin{proof} Since $\bdy\sigma$ has uniformly bounded length, to prove the lemma it
suffices to find $T \ge 0$ so that if $e_1,e_2$ are adjacent edges in $\bdy\sigma$
then $X_{e_1}\ceq{[T]} X_{e_2}$. Let $v$ be the vertex of $K$ at which $e_1,e_2$ meet.

Notice that if $e,e'$ are edges of $K$ and if the paths $\pi e$, $\pi e'$ contain a common
edge $e''$ of $T$ then $X_e \ceq{[2C]} X_{e''} \ceq{[2C]} X_{e'}$ and so $X_e \ceq{[4C]}
X_{e'}$.

If the map $\pi$ folds $\bdy\sigma$ at the vertex $v$ it follows that $\pi e_1$, $\pi e_2$
contain a common edge of $T$, and so $X_{e_1} \ceq{[4C]} X_{e_2}$. 

Suppose that the map $\pi$ does not fold $\bdy\sigma$ at the vertex $v$. Let $w=\pi v$ and
let $\pi^\inv(w) \intersect \bdy\sigma = \{v\} \union Y$. Since the map $\pi$
is a bijection from vertices of $K$ to fat vertices of $T$, each $y \in Y$ is an
interior point of some edge $a_y$ of $\bdy\sigma$. Each of the paths
$\pi(a_y)$ determines an edge in the link of $w$, each of the paths $\pi(e_1)$,
$\pi(e_2)$ determines a vertex in the link of $w$, and the union of these vertices and
edges is a connected subset of the link of $w$, since $T$ is a tree. We can
therefore find a sequence $e_1=a_0, a_1, \ldots, a_{M-1}, a_M = e_2$ of
edges in $\bdy\sigma$ such that for $m=1,\ldots,M$ the paths $\pi(a_{m-1})$ and
$\pi(a_{m})$ contain a common edge of $T$, implying that $X_{a_{m-1}} \ceq{[4C]}
X_{a_{m}}$. We therefore have $X_{e_1} \ceq{[4CM]} X_{e_2}$. Since $M$ is bounded by the
length of $\bdy\sigma$ we are done.
\end{proof}

Now we can prove the remaining contentions needed to establish Lemma~\ref{Lemma2complex}:

\begin{lemma} 
\label{LemmaCofiniteness}
The following properties hold:
\begin{enumerate}
\item $K/H$ is finite.
\item $K$ is uniformly locally finite along edges, that is, there is a uniform upper bound
on the number of 2-cells whose  boundary contains any fixed edge $e$ of $K$.
\item There is a uniform upper bound on the number of edges of $K$ whose image crosses
any fixed edge of $T$.
\end{enumerate}
\end{lemma}

\begin{proof} 

Finiteness of vertices, edges, and faces of $K/H$ is an application of part (1) of the
Coboundedness Principle~\ref{PropCoboundednessPrinciple} with appropriate patterns, as follows.
In the case of vertices, we simply use the pattern of fat vertex spaces in $X$. For edges, each
edge $e$ of $K$ has a corresponding edge space $X_e$, this pattern is locally finite, and for
$h \in H$ we clearly have $h \cdot X_e \ceq{[R]} X_{h \cdot e}$ for some uniform $R$. For
\nb{2}cells, Lemma~\ref{LemmaWhateveritwas} shows that if we assign to each 2-cell
$\sigma$ the pattern element $X_\sigma = \union_{e \subset \sigma} X_e$, then this
pattern is locally finite and for $h \in H$ we have $h \cdot X_\sigma \ceq{[R]} X_{h \cdot
\sigma}$ for some uniform $R$.

To prove that $K$ is uniformly locally finite along an edge $e$, from
Lemma~\ref{LemmaWhateveritwas} there is a constant $R$ such that if the boundary a
\nb{2}cell $\sigma$ contains $e$ then for any other edge $e' \subset \bdy\sigma$ we have
$X_{e'} \ceq{[R]} X_e$. Since the pattern of edge spaces is locally finite, there is a
uniformly finite number of such edges $e'$. Since $\bdy\sigma$ has uniformly bounded
length and $\sigma$ is completely determined by its boundary, there is a uniformly finite
number of such \nb{2}cells $\sigma$.

The proof of (3) is nearly identical to the proof of (2), using property (iii) for edges of $K$ in
place of Lemma~\ref{LemmaWhateveritwas}.
\end{proof}

This finishes the construction of $K$ and the proof of its properties, thus finishing
Lemma~\ref{Lemma2complex}.

\subsubsection{Proof of the Tracks Theorem \ref{TheoremTracks}}
\label{Tracks}

Stallings discovered the principle that analogies with \nb{3}manifold techniques could be applied
to solve problems about general finitely generated groups \cite{Stallings:ends}. An example of
this principle is Dunwoody's theory of tracks in \nb{2}complexes \cite{Dunwoody:Accessible}, which
are analogous to Haken's normal surfaces in a triangulated \nb{3}manifolds \cite{Haken:normal} ---
as it turns out, a normal surface intersects the \nb{2}skeleton in a track, and the complexity of
that track measures the complexity of the normal surface. Dunwoody \cite{Dunwoody:Accessible} used
tracks to prove accessibility of finitely presented groups, an analogy to Haken's
theorem on boundedness of hierarchies of compact \nb{3}manifolds.

Meanwhile, a completely different approach to \nb{3}manifold theory was developed by Meeks and Yau
\cite{MeeksYau:Equivariant}, using minimal surfaces with respect to a Riemannian metric on the
manifold to solve problems about \nb{3}manifolds. Minimal surfaces and normal surfaces were
combined by Jaco and Rubinstein \cite{JacoRubinstein:PLminimal} in their theory of PL minimal
surfaces in \nb{3}manifolds, resulting in a particularly elegant setting for \nb{3}manifold
theory. Casson's Beijing notes on \nb{3}manifolds \cite{Casson:ThreeDimensional} gave further
simplifications in PL minimal surface theory, by using hyperbolic ideal triangles to impose a
``hyperbolic structure'' on the \nb{2}skeleton minus the \nb{0}skeleton. 

The theme of our proof is that Casson's approach to PL minimal surfaces, combined with Dunwoody's
tracks, can be used to solve problems about general groups acting on simply connected simplicial
\nb{2}complexes, even in contexts where the simplicial complex might not be locally finite. This
theme was discovered independently by Niblo, who uses it to give a new proof of Stallings' ends
theorem \cite{Niblo:ends}.

\bigskip

Recall the setting of the Tracks Theorem: a simply connected, simplicial \nb{2}complex $L$, a
cocompact action of a group $H$ on $L$, and a quasi-isometry $\pi \from L \to T$ to a tree, so
that $L$ is uniformly finite along edge interiors, and so that the tracks $\pi^\inv(t)$, for $t$ in
the interior of an edge of $T$, are uniformly finite.

In the course of the proof we will make use of the ends of $T$ and of~$L$. Ends cannot be defined
in Freudenthal's sense, using complementary components of compact sets, because $T$ and $L$ are
not locally compact. Instead, we define ends in terms of complementary components of bounded
subsets: the collection of bounded subsets forms a direct system under inclusion; the
collection of unbounded complementary components of bounded subsets forms an inverse system
under inclusion; and an end is defined to be an element of the inverse limit. This definition
applies both in $T$ and in $L$.

Alternatively, the tree $T$ with its simplicial metric is a Gromov hyperbolic space, the map $\pi
\from L \to T$ is a quasi-isometry and so $L$ is also Gromov hyperbolic, and $\pi$ induces a
bijection of Gromov boundaries. For a simplicial complex quasi-isometric to a tree there is a
bijection between the ends and the Gromov boundary. It follows that $\pi \from L \to T$ induces a
bijection of ends.

Now we follow Casson by imposing a ``hyperbolic structure'' on $L-L^0$ as follows.
For each edge $e$ put a metric on $e - (\bdy e)$ making it isometric to $\R$. Choose a
\emph{midpoint} $\midpoint(e)\in e$. For each 2-simplex $\sigma$ extend the metrics on the
components of $(\bdy \sigma) -\sigma^0$ to a metric on $\sigma - \sigma^0$ which is isometric to a
hyperbolic ideal triangle, so that for each edge $e$ of $\sigma$ the point $\midpoint(e)$ matches
the base of the perpendicular from the opposite ideal vertex. This matching condition guarantees
that $L-L^0$ is a complete metric space. It also guarantees that any simplicial isomorphism of $L$
restricts to an isometry of $L-L^0$. 

For each track $\tau$ in $L$, define the \emph{weight} to be $w(\tau) = \card{\tau \intersect
L^1}$, and define the \emph{length} $\ell(\tau)$ to be the sum of the hyperbolic lengths of the
components of $\tau$ intersected with each 2-simplex. Define the PL-length of $\tau$ to be the
ordered pair $(w(\tau),\ell(\tau)) \in \Z_+ \cross (0,\infinity)$. We compare PL-lengths using the
dictionary ordering on $\Z_+ \cross (0,\infinity)$.

Consider an essential track $\tau$ in $L$. We say that $\tau$ is \emph{weight minimal} (called
simply ``minimal'' in \cite{Dunwoody:Accessible}) if $w(\tau)$ is minimal among all essential
tracks. We say that $\tau$ is \emph{PL-minimal} if $(w(\tau),\ell(\tau))$ is minimal among all
essential tracks in $L$. We will need the fact that a weight minimal track intersects any
\nb{2}simplex in at most one arc (\cite{Dunwoody:Accessible}, Proposition 3.1 and following
comments).

Given a real number $\epsilon \in (0,1)$, for any vertex $v$ there is a closed regular
neighborhood $N^\epsilon(v)$ whose boundary link $\link^\epsilon(v)$ consists of a union of
horocyclic segments of length $\epsilon$ in each 2-simplex incident to $v$; this follows from the
matching condition. Moreover $N^\epsilon(v)\intersect N^\epsilon(w) = \emptyset$ for $v \ne w$;
this follows because in an ideal hyperbolic triangle, the three horocyclic segments of length
$\epsilon$ about the three vertices are disjoint. A component of $\link^\epsilon(v)$ is called an
\emph{$\epsilon$-horocyclic track near $v$}; we take the liberty of failing to mention either the
$\epsilon$ or the $v$. Note that an $\epsilon$-horocyclic track is isotopic in $L-L^0$ to an
$\epsilon'$-horocyclic track for any $\epsilon'\le \epsilon$.

\begin{lemma} Suppose also that $L$ has an essential track.
Then there exists an essential track $\tau \subset L$ which is either
PL-minimal or horocyclic.
\label{LemmaPLOrHorocyclic}
\end{lemma}

\begin{proof} Let $w_\Min$ be the minimum of the positive integer $w(\tau)$. Let $\tau_i$ be a
sequence with $w(\tau_i)=w_\Min$ and with $\ell(\tau_i)$ approaching the infimum. Since $H$ acts
cocompactly on $L$ and isometrically on $L-L^0$, after pulling back by appropriate elements of $H$
and passing to a finite subsequence we may assume that each $\tau_i$ intersects a certain edge
$e$. There are only finitely many isotopy classes of tracks with weight $w_\Min$ that
intersect $e$, and so passing to a further subsequence we may assume that the isotopy class of
$\tau_i$ is constant.

If any $\ell(\tau_i)=0$ then we are done, for in that case $w(\tau_i)=w_\Min=1$ and $\tau_i$ is a
point on an edge $e$ which does not bound any \nb{2}simplex. So, we may assume that $\ell(\tau_i)
\ne 0$.

Suppose there exists a vertex $v$ which is an accumulation point in $L$ of the sequence of tracks
$\tau_i$. Passing to a subsequence we may assume that $v =\lim x_i$ for some $x_i \in \tau_i$.
Since $v$ is infinitely far away from any point in $L-L^0$ with respect to the hyperbolic
structure, and since $\ell(\tau_i)$ is bounded, it follows that there is an $\epsilon>0$ such that
$\tau_i \subset N^\epsilon(v)$ for sufficiently large $i$. But now it is clear that $\tau_i$ is
isotopic to a horocyclic track near~$v$.

Now suppose that no vertex is an accumulation point of the tracks $\tau_i$. Since the isotopy class
of $\tau_i$ in $L-L^0$ is fixed it follows that $\tau_i$ is contained in a bounded subset of
$L-L^0$, and so by the Ascoli-Arzela theorem $\tau_i$ converges to a track $\tau$ with
$(w(\tau),\ell(\tau))$ minimized.
\end{proof}

\begin{lemma} There exists $\epsilon>0$ depending only on $L$ with the following property. Suppose
that $\tau_1,\tau_2$ are two essential tracks in $L$, each of which is either PL-minimal or
$\epsilon$-horocyclic. Then either
$\tau_1=\tau_2$ or $\tau_1 \intersect \tau_2 = \emptyset$.
\label{LemmaDisjointTracks}
\end{lemma}

\begin{proof} If no PL-minimal tracks exist then the lemma is obvious, with any $\epsilon \in
(0,1)$.

Suppose that PL-minimal tracks exist, and let $(w,\ell)$ be the least PL-length. If $\ell=0$ then
the lemma is again obvious, so we assume $\ell>0$.

Here's how to choose $\epsilon$. First pick an arbitrary number $\epsilon' \in (0,1)$. Now choose
$\epsilon \in (0,\epsilon')$ so that the distance from $\link^\epsilon(v)$ to
$\link^{\epsilon'}(v)$ is at least $\ell$ for any vertex $v$; $\epsilon = \epsilon' e^{-\ell}$ will
do nicely. 

Given a vertex $v$, clearly there are no PL-minimal tracks entirely contained in
$N^{\epsilon'}(v)$, for any track in $N^{\epsilon'}(v)$ may be shortened by isotoping towards $v$.
It follows that every PL-minimal track is disjoint from $N^\epsilon(v)$. The lemma then follows as
long as one or both of $\tau_1,\tau_2$ are $\epsilon$-horocyclic.

For the rest of the proof we may therefore assume that both
$\tau_1,\tau_2$ are PL-minimal. A simple variation argument (see
\cite{Casson:ThreeDimensional}) shows that a PL-minimal track $\tau$
has the following ``minimal surface'' properties: 
\begin{itemize}
\item Each segment of intersection of $\tau$ with a 2-simplex is 
geodesic.
\item For oriented edge $e$ of $L$, each point $x \in\tau\intersect e$,
and each \nb{2}simplex $\sigma$ of $L$, letting $\theta(\tau,\sigma,x)$ be
the angle between $\tau\intersect\sigma$ and $e$ at the point $x$, we have
$$C(\tau,x) \equiv \Sum_{\sigma} \cos \bigl( \theta(\tau,\sigma,x) \bigr) = 0
$$
\end{itemize}
The fact that $C(\tau,x)=0$ can be seen by noting that $C(\tau,x)$ is the
infinitesmal rate of change of $\ell(\tau)$ with respect to the position of
$x$ on $e$.

\bigskip

Suppose now that $\tau_1 \intersect \tau_2 \ne \emptyset$. We argue by contradiction, first in a
special case and then in general. 

\subparagraph*{Special Case:} $\tau_1$ and $\tau_2$ are in general position: $\tau_1 \intersect
\tau_2$ consists of transverse intersection points in the interiors of 2-simplices. Following
Dunwoody, let $J = L^1 \intersect (\tau_1 \union \tau_2)$, and noting that $\card{J \intersect
\bdy\sigma}$ is even for each \nb{2}simplex $\sigma$, we obtain a finite track pattern $T_J$,
well-defined up to isotopy, satisfying $T_J \intersect L^1 = J$. 

Applying Proposition 3.2 of \cite{Dunwoody:Accessible}, $T_J$ is a union of two essential,
weight-minimal tracks, $T_J = \tau'_1 \union \tau'_2$, and so $w(\tau'_i) = w(\tau_i) = w$, for
$i=1,2$. 

\begin{remark} In \cite{Dunwoody:Accessible}, a track was defined to be essential if each of its
two complementary components has infinitely many vertices. In our present context, it is more
appropriate to require the complementary components to be unbounded, and the proof of Dunwoody's
Proposition 3.2 is easily checked with that change.
\end{remark}

Note that $T_J$ can be constructed in its isotopy class by doing a certain cut-and-paste on
$\tau_1 \union \tau_2$: for each \nb{2}simplex $\sigma$, and each $x \in \sigma \intersect \tau_1
\intersect \tau_2$, there are two choices for the cut-and-paste, but the choice is determined by
the overall structure of $\sigma \intersect \tau_1$ and $\sigma \intersect
\tau_2$ (see Figure \ref{FigureCutAndPaste}). We thus have $T_J = \tau'_1
\union \tau'_2$ where 
$$\ell(\tau'_1) + \ell(\tau'_2) = \ell(\tau_1) +
\ell(\tau_2) = 2\ell
$$
and so there exists $i=1,2$ such that $\ell(\tau'_i) \le \ell$.

\begin{figure}
\centeredepsfbox{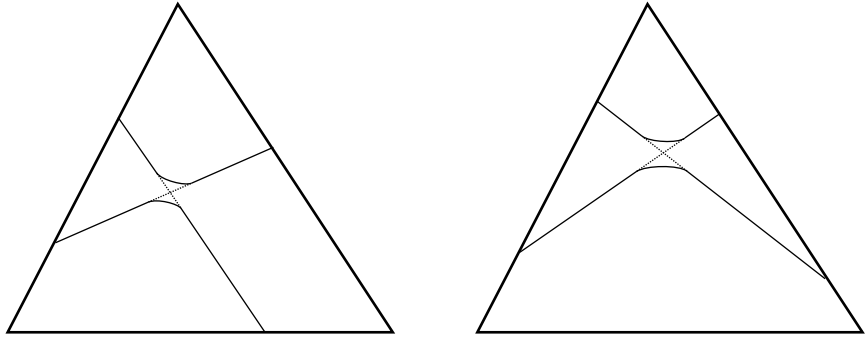}
\caption{If $\tau_1,\tau_2$ are PL-minimal then they each intersect a
2-simplex $\sigma$ in at most one arc. If $\tau_1,\tau_2$ are in general
position then these arcs intersect in at most one point $x$. Letting $J =
L^1 \intersect (\tau_1\union \tau_2)$, the track $T_J$ can be obtained
from $\tau_1, \tau_2$ by cut-and-paste. There two possibilities for the
cut-and-paste at $x$, but the choice is determined by the arcs in
$\sigma$. After cutting-and-pasting, PL length can be reduced.}
\label{FigureCutAndPaste}
\end{figure}

However, because of the assumption that $\tau_1 \intersect \tau_2$ is nonempty, say $x\in \sigma
\intersect \tau_1 \intersect \tau_2$ for some \nb{2}simplex $\sigma$, it follows that the
cut-and-paste operation produces a corner of $\tau'_i$ at $x$ with nonzero exterior angle. The
length $\ell(\tau'_i)$ can therefore be shortened below the minimum length $\ell$, a contradiction
which shows that $\tau_1$, $\tau_2$ cannot have a nonempty intersection in general position.

\subparagraph*{General Case:} $\tau_1$ and $\tau_2$ are not in general position. In particular,
there is a point $x_0 \in L^1 \intersect \tau_1\intersect \tau_2$ such that $\tau_1, \tau_2$ do
not coincide in a neighborhood of $x$. Since $C(\tau_1,x_0) = C(\tau_2,x_0) = 0$, it follows that
there is a 2-simplex $\sigma$ incident to $x_0$ such that $\theta(\tau_1,\sigma,x_0) \ne
\theta(\tau_2,\sigma,x_0)$. Let $\phi_0 = \abs{\theta(\tau_1,\sigma,x_0) -
\theta(\tau_2,\sigma,x_0)}$, so $\phi_0 > 0$.

Now we borrow Casson's version of the ``Meeks-Yau trick'': perturb $\tau_1$ and $\tau_2$, raising
their lengths by an arbitrarily tiny amount and making their cosine sums differ from zero by an
arbitrarily tiny amount, to put $\tau_1$ and $\tau_2$ in general position. Do this in such a way
as to make them intersect at a point $x$ in $\sigma$ near $x_0$ at an angle $\phi$ arbitrarily
close to $\phi_0$; the tinier the perturbation, the closer $\phi$ will be to $\phi_0$.

Now use the argument from the special case above, to produce two essential, weight minimal tracks
$\tau'_1, \tau'_2$. The tracks $\tau'_1,\tau'_2$ have corners at $x$ each with exterior angle
$\phi$, which is arbitrarily close to $\phi_0$, and their cosine sums where they intersect $L^1$
are arbitrarily small. It follows that there is a $\delta>0$, depending only on $\phi_0$, so that
by straightening $\tau'_1,\tau'_2$ near $x$ the values of $\ell(\tau'_1)$, $\ell(\tau'_2)$ can
each be reduced by at least $\delta$. However, the perturbation can be chosen so tiny that
$$\ell(\tau'_1) + \ell(\tau'_2) < \ell(\tau_1) + \ell(\tau_2) + \delta = 2 \ell+\delta
$$ 
and therefore one of $\ell(\tau'_1), \ell(\tau'_2)$ can be reduced below
$\ell$, a contradiction. 
\end{proof}

Now we put the pieces together to prove the following:

\begin{lemma}[Tracks exist]
\label{LemmaTracksExist}
Under the hypotheses of the Tracks Theorem, if $L$ is unbounded then there exists an
$H$-equivariant, essential, nonempty track pattern $P \subset L$.
\end{lemma}

\begin{proof} Since $L$ is unbounded and $L,T$ are quasi-isometric, $T$ is also unbounded. 

We use the coboundedness of the action of $H$ on $L$ to show that $L$ and $T$ have at least two
ends. For a tree with a cococompact group action even more is true: the tree is either bounded,
quasi-isometric to a line, or ``bushy'' meaning that the tree is quasi-isometric to a tree of
constant valence $\ge 3$. This was proved in \cite{MSW:QTOne} for trees of uniformly bounded
valence with the additional condition that the bushy tree of the conclusion has constant finite
valence; but if one drops the hypothesis of uniformly bounded valence one need only weaken the
conclusion to allow a bushy tree whose valence is a constant of unrestricted cardinality. The same
proof works for a tree with a cobounded quasi-action, and we can apply that to our present
situation since the cobounded action of $H$ on $L$ can be quasi-conjugated to give a cobounded
quasi-action on $T$. 

For convenience here is a detailed proof, using the coboundedness of the action of $H$ on $L$,
that if $L$ is unbounded and quasi-isometric to a tree $T$ then $L$ has at least two ends.

A subset of a metric space if \emph{$r$-deep} if that subset contains an ambient metric ball of
radius $r$. Since the tree $T$ is unbounded, there is an infinite ray leaving every compact set,
and so for each $r>0$ there exists an edge $e$ such that $T-e$ has at least two $r$-deep
components. Since $L$ and $T$ are quasi-isometric, it follows that there is a constant $k>0$ such
that for each $s>0$ there exists a point $x$ such that $L-B(x,k)$ has at least two $s$-deep
components. Since $H$ acts coboundedly on $L$, there exists $d>0$ such that for all $p,x \in L$
there exists $h \in H$ such that $d_L(p, h \cdot x) \le d$. This implies that, letting $l=k+d$,
for any $s>0$ and any point $p \in L$ the set $L-B(p,l)$ has at least two $s$-deep components. Fix
a point $p_0=q_0$, fix $s > l$, and let $C_0,D_0$ be two $s$-deep components of
$L-B(p_0,l)=L-B(q_0,l)$. We show that $C_0,D_0$ are both unbounded. For $i\ge 1$, inductively pick
$p_i$ so that $B(p_i,l) \subset B(p_i,s) \subset C_{i-1}$. There are at least two $s$-deep
components of $L-B(p_i,l)$, only one of which can contain $B(p_{i-1},l)$, so pick $C_i$ to be a
different one. Similarly, pick $q_i$ so that $B(q_i,l) \subset B(q_i,s) \subset D_{i-1}$, and pick
an $s$-deep component $D_i$ of $L-B(q_i,l)$ which does not contain $B(q_{i-1},l)$. It follows that
$p_1,p_2,p_3,\ldots$ is a sequence of points in $C_0$ which gets arbitrarily far from $p_0$, and
$q_1,q_2,q_3,\ldots$ is a sequence in $D_0$ which gets arbitrarily far from $q_0=p_0$. This shows
that $L$ has at least two ends. The tree $T$, being quasi-isometric to $L$, also has at least two
ends.

Since $T$ has at least two ends, there is a point $t$ in the interior of some edge such that $T-t$
has two unbounded components, and it follows that the finite track pattern $\pi^\inv(t)$ has an
essential component. This shows that $L$ has a finite, essential track. 

Applying Lemmas \ref{LemmaPLOrHorocyclic} and \ref{LemmaDisjointTracks}, it follows that $L$ has an
essential track $\tau$ which is either PL-minimal or $\epsilon$-horocyclic, where $\epsilon$ is
chosen as in Lemma \ref{LemmaDisjointTracks}. Since $H$ acts on $L$ isometrically, it follows that
$h(\tau)$ is PL-minimal or $\epsilon$-horocyclic for each $h \in H$. Therefore, for any $h,h' \in
H$ we apply Lemma \ref{LemmaDisjointTracks} and conclude that either $h(\tau) = h'(\tau)$ or
$h(\tau)\intersect h'(\tau) = \emptyset$. Thus, $P = H \cdot \tau$ is an $H$-equivariant track
pattern.
\end{proof}

\paragraph{Relative tracks.} We shall need a relative version of Lemma~\ref{LemmaTracksExist}. A
subcomplex $F \subset L^1$ is called \emph{peripheral} if each component is compact
and there is an open collar neighborhood $N_F$ equipped with a homeomorphism $(N_F,F) \homeo (F
\cross [0,1), F \cross 0)$. 

\begin{lemma}[Relative tracks exist]
\label{LemmaRelativeTracksExist}
Under the hypotheses of the Tracks Theorem, if $L$ is unbounded, and if $F$ is an $H$-equivariant
peripheral subcomplex of $L$, then there exists an $H$-equivariant, essential, nonempty track
pattern $P\subset L$ such that $P\intersect F = \emptyset$.
\end{lemma}

\begin{proof} Given a track $\tau \subset L$, we alter the definition of the weight $w(\tau)$ by
assigning infinite weight to any point of $\tau \intersect F$, and so a track has finite weight if
and only if it is disjoint from $F$. Assuming that an essential track of finite weight exists, the
proof of Lemma~\ref{LemmaTracksExist} goes through to show that there is a nonempty,
$H$-equivariant, essential track pattern $P$ whose components each have finite weight.

It remains to show that there exists an essential track in $L$ disjoint from~$F$. 

We generalize the concept of a track pattern as follows: a \emph{pseudopattern} is a \nb{1}complex
$P$ embedded in $L$ such that for each \nb{2}simplex $\sigma$ of $L$, each component of $\sigma
\intersect P$ is either a circle in the interior of $\sigma$ or a properly embedded arc with
endpoints disjoint from the vertices of $\sigma$; any such arc with endpoints on the same edge of
$\sigma$ is called a \emph{bight} of $P$. A \emph{pseudotrack} is connected pseudopattern.
The concept of essentiality is defined for a pseudotrack exactly as for a track. An essential
pseudotrack is not a circle in the interior of any \nb{2}simplex.

First we show that $L$ contains an essential pseudotrack disjoint from~$F$. By applying
Lemma~\ref{LemmaTracksExist} we obtain an essential track $\tau$ in $L$, but $L$ may intersect
$F$. Consider an essential decomposition $L-\tau = U \union V$. Let $F_\tau$ be the union of the
components of $F$ that intersect $\tau$, so $F_\tau$ is a finite \nb{1}complex. Let $P'$ be the
pseudopattern obtained by using product structure on the collar neighborhood $N_F = F \cross [0,1)$
to push $\tau\union (F_\tau \intersect U)$ into $L-F$. To be precise, we assume that $\tau
\intersect N_F$ has the form $A\cross [0,1)$ for some finite subset $A \subset F_\tau$, and we
take 
$$P' = \bigl(\tau - (A \cross [0,1/2))\bigr) \union \bigl((U\intersect F_\tau) \cross 1/2\bigr)
$$ 
There is an essential decomposition $L-P' = U' \union V'$, where $V'$ is the closure in $L-P'$
of $V\union \bigl( (U\intersect F_\tau) \cross [0,1/2) \bigr)$, and where $U' = (L-P') - V'$. It
follows, by the same argument for ordinary track patterns, that $P'$ contains an essential
pseudotrack $\tau'$. By construction, $\tau'$ is disjoint from $F$.

\begin{figure}
\centeredepsfbox{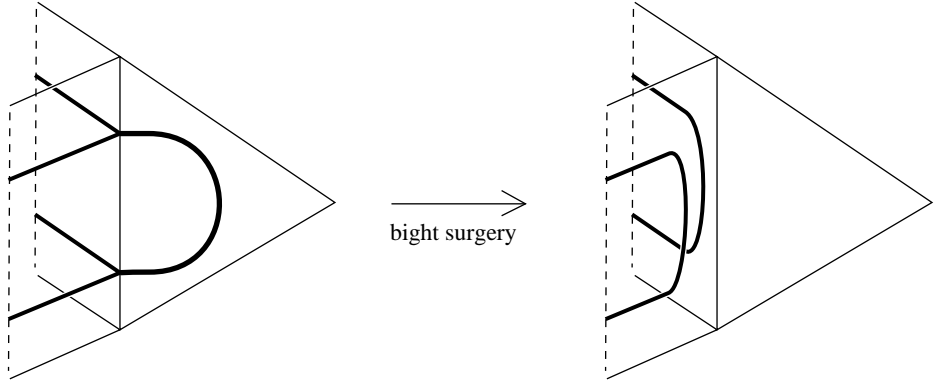}
\caption{Eliminating a bight by surgery}
\label{FigBight}
\end{figure}

To go from an essential pseudotrack $\tau'$ disjoint from $F$ to an essential track disjoint from
$F$, we inductively remove the bights, using a surgery procedure illustrated in
Figure~\ref{FigBight}. To describe this procedure, consider an essential pseudotrack $\tau'$ with
essential decomposition $L-\tau'=U' \union V'$. Let $b \subset \sigma$ be a bight of $\tau'$, with
endpoints $x_0,x_1$ on an edge $e$ of $L$. We may choose $b$ to be innermost in $\sigma$, so the
segment $[x_0,x_1] \subset e$ has interior disjoint from $\tau'$; choose the notation so that the
interior of $[x_0,x_1]$ is contained in $V'$. Let $\sigma=\sigma_0,\sigma_1,\ldots,\sigma_N$
be the \nb{2}simplices incident to $e$. Let $s_i$ be a regular neighborhood of $x_i$ in $\tau'$,
consisting of $N+1$ arcs sharing the common endpoint $x_i$, with opposite endpoint $y_{in}$ in
$\sigma_n$. For $n=1,\ldots,N$ let $a_n$ be an arc in the interior of $\sigma_n$ with $a_n
\intersect \tau' = \bdy a_n = \{y_{0n},y_{1n}\}$. The pseudopattern 
$$P'' = (\tau' - (b \union s_0 \union s_1)) \union a_1 \union \cdots \union a_N
$$
has strictly smaller weight than $\tau'$, and moreover $P'' \intersect L^1 \subset \tau' \intersect
L^1$ implying that $P'' \intersect F = \emptyset$. There is an essential decomposition $L-P'' =
U'' \union V''$ described as follows: the connected pseudopattern $b \union s_0 \union s_1 \union
a_1 \union\cdots\union a_N$ is inessential, with a bounded complementary component $W$; let
$U''$ be the closure in $L-\tau''$ of $U' \union W$, and let $V'' = (L-\tau'') - U''$. It
follows that $P''$ has a component $\tau''$ which is an essential pseudopattern disjoint from $F$
of smaller weight than $\tau'$, completing the induction.
\end{proof}

\paragraph{The Haken finiteness argument.} Consider a properly nesting sequence 
$P_1 \subset \cdots \subset P_K$
of $H$-equivariant, essential track patterns, such that no distinct components of $P_K$ are
isotopic. We claim that there is a bound $K \le K_0$ depending only on
$L$ and $H$. To see why, the quotient complex $L/H$ is finite, and so as noted in
\cite{Dunwoody:Accessible} there exists $K_0$ such that any set of pairwise disjoint tracks in
$L/H$ of cardinality $> K_0$ contains two tracks which are isotopic in $L/H$. But if $P_K / H$ has
a pair of isotopic components then $P_K$ has a pair of isotopic components. This shows that $K
\le K_0$. It follows that there exists a sequence $P_1\subset\cdots\subset P_K$ as above which is
maximal in the sense that it cannot be extended. Assuming maximality, to finish the Tracks Theorem
we need to prove that the closure of every component of $L-P_K$ is a bounded \nb{2}complex.

\paragraph{Boundedness.} Consider then a sequence $P_1 \subset \cdots\subset P_K$ as above,
and let $C$ be the closure of a component of $L-P_K$. Assuming maximality of the sequence, we must
show that $C$ is bounded.

Let $\Stab(C)$ denote the group of simplicial automorphisms of $L$ that stabilize $C$. Noting that
the quotient $C / \Stab(C)$ is the closure of a component of $L/H - P_K/H$, it follows that
$\Stab(C)$ acts cocompactly on $C$.

Note that each component of the frontier $\Fr(C)$ of $C$ in $L$ is a component of $P_K$; we call
these the \emph{frontier tracks} of $C$. The action of $\Stab(C)$ restricts to an action on the
union of the frontier tracks of $C$.

The cell structure on $C$ is described as follows. The \nb{0}skeleton is $C\intersect (L^0 \union
P_K^0)$. The \nb{1}skeleton is $C \intersect (L^1 \union \intersect P_K)$. Each \nb{2}cell $c$ of
$C$ is contained in a \nb{2}simplex $\sigma$ of $L$, $c$ is the closure of a component of
$\sigma-P_K$, and so $c$ is one of the following: $c=\sigma$; or $c$ is a quadrilateral with three
sides in $\bdy\sigma$ and one side in $\Fr(C)$; or $c$ is a pentagon with three sides in
$\bdy\sigma$ and two sides in $\Fr(C)$; or $c$ is a hexagon with three sides in
$\bdy\sigma$ and three sides in $\Fr(C)$; or $c$ is a quadrilateral with two sides in $\bdy\sigma$
and two sides in $\Fr(C)$; or $c$ is a triangle with two sides in $\bdy\sigma$ and one side in
$\Fr(C)$.

Next we prove that the inclusion $C \inject L$, which is clearly coarse Lipschitz, is actually a
quasi-isometric embedding. Consider two vertices $x,y$ of $C$ connected by a geodesic path $\gamma$
in $L^1$ consisting of at most $m$ edges or half-edges of $L$. Note that there are at most $m$
frontier tracks $\tau$ of $C$ such that $\gamma \intersect A_\tau\ne \emptyset$, because $A_\tau$
must contain a vertex of $\gamma$. For each such $\tau$, let $\gamma_\tau$ be the maximal
subsegment of $\gamma$ such that $\bdy\gamma_\tau \in A_\tau$. Since the diameter of $\tau$ is
uniformly bounded in $L$, the segment $\gamma_\tau$ has uniformly bounded length. Choosing some
order on the set of $\tau$, inductively replace each $\gamma_\tau$ by an edge path in $\tau$ of
uniformly bounded length. As each $\gamma_\tau$ is replaced, the identity of each subsequent
$\gamma_{\tau'}$ may change as one or both of its endpoints are eliminated by the previous
replacement, but this change at worst decreases the length of $\gamma_{\tau'}$. It follows that
$x,y$ are connected by a path in $C$ whose length is bounded by an affine function of the length
of $\gamma$, showing that $C \inject L$ is a quasi-isometric embedding.

We may draw several conclusions from the fact that $C \inject L$ is a quasi-isometric embedding.
First, ``boundedness'' of $C$ is unambiguously defined: $C$ is a bounded \nb{2}complex if and only
if $C$ is a bounded subset of~$L$. Second, the map $\pi \from C \to T$ induces an injection of
ends, and so by restricting the range the map $\pi \from C \to T_C$ induces a bijection of ends.
Third, if $\tau$ is a track in $L$ which is contained in $C$, is disjoint from $\Fr(C)$, and is
essential in $C$, then $\tau$ is essential in $T$ and $\tau$ is not isotopic to any component of
$P_K$. 

Assuming now that $C$ is not bounded, we prove that $P_1 \subset \cdots \subset P_K$ is not
maximal. The frontier $\Fr(C)$ is a $\Stab(C)$-equivariant, peripheral subcomplex. Applying
Lemma~\ref{LemmaRelativeTracksExist}, $C$ has a $\Stab(C)$-equivariant, essential track pattern
$P_C$ disjoint from $\Fr(C)$, and so disjoint from $P_K$. Since the inclusion $C
\inject L$ is a quasi-isometry, each component of $P_C$ is essential in $L$, and no component is
isotopic to a comonent of $P_K$. Acting on $P_C$ by the entire group $H$, we obtain an
$H$-equivariant, essential track pattern $H \cdot P_C$ disjoint from $P_K$ and with no component
isotopic to a component of $P_K$. Letting $P_{K+1} = P_K\union H \cdot P_C$, it follows that $P_1
\subset \cdot \subset P_K$ is not maximal.

There is a slight flaw in the above argument, because $C$ is not necessarily a simplicial complex.
We can triangulate $C$ in a $\Stab(C)$ equivariant manner by putting a new vertex in the
interior of any \nb{2}cell $c$ that is not already a triangle, and coning off the boundary of $c$.
But now the problem is that $P_C$, which is a track pattern with respect to the triangulation of
$C$, may only be a pseudopattern with respect to the original triangulation of $L$. In the proof
of Lemma~\ref{LemmaRelativeTracksExist}, the key step is the production of an essential track
$\tau$ in $C$ disjoint from $\Fr(C)$, but $\tau$ may only be a pseudotrack with respect to the
triangulation of $L$; however, any bights can be removed from $\tau$ by the same
inductive process described in the proof of Lemma~\ref{LemmaRelativeTracksExist}, and so we may
indeed assume that $\tau$ and $P_C$ are true track patterns in $L$.

This completes the proof of the Tracks Theorem.

\subsection{Proof of the Tree Rigidity Theorem}
\label{SectionTreeRigidityProof}
Recall the setup of the theorem. $\G$ is a finite type, finite depth, irreducible graph of
groups with Bass-Serre tree of spaces $X \to T$. The finitely generated group $H$ quasi-acts
properly and coboundedly on $X$, and this quasi-action coarsely respects vertex and edge spaces,
meaning that there is a quasi-action of $H$ on $\VE(T)$ and a constant $r \ge 0$ such
that for each $h\in H$ and $a\in\VE(T)$ we have $h\cdot X_a\ceq{[r]} X_{h \cdot a}$, and for each
$a' \in\VE(T)$ there exists $a\in\VE(T)$ such that $X_{h \cdot a} \ceq{[r]} X_{a'}$. 

The depth of a vertex or edge space $X_a$ is determined by the position of $X_a$ in the coarse
inclusion lattice of vertex or edge spaces: $X_a$ has depth zero if and only if $X_a$ is maximal in
the coarse inclusion lattice; by induction, $X_a$ has depth $n$ if and only if $X_a$ does not have
depth $<n$ and $X_a$ is maximal in the coarse inclusion lattice among vertex and edge spaces that
do not have depth $<n$. Since a quasi-isometry preserves coarse inclusions, 
it follows that $H$ preserves depth.

We must prove that there is a finite type, finite depth, irreducible graph of groups $\G'$ with
Bass-Serre tree of spaces $X' \to T'$, and with fundamental group $\pi_1\G' \approx H$, so that
the induced $H$-almost equivariant quasi-isometry $X \to X'$ coarsely respects vertex and edge
spaces.

Two special cases have been proved: Proposition~\ref{PropNoTracks} handles the case
where the edge spaces of $X$ are isolated, meaning that no two distinct edge spaces of $X$ are
coarsely equivalent; and Proposition~\ref{PropDepthOne} handles the case when every edge space has
depth one.

The induction argument for the general case is, in outline and in many details, the same as in
Proposition~\ref{PropNoTracks}. But the basis step in the general case  will need the analysis of
the homogeneous case carried out in \cite{MSW:QTOne} (see Theorem~\ref{homogeneous}), and the
induction step will require Proposition~\ref{PropDepthOne}. 

Recall some notation and terminology. The tree $T$ is filtered by the subforests $F_0 \subset F_1
\subset\cdots\subset F_N = T$, where $F_n$ consists of all vertices and edges of depth $\le n$.
The components of $F_n$ are called the depth $n$ \emph{flotillas}. For each subset $\A \subset
\VE(T)$ we denote $X_\A = \union_{a \in \A} X_a$. 

For each $n = 0,\ldots,N$ consider the pattern in $X$ defined by
$$\P_n = \{X_\A \suchthat \text{$\A$ is a depth $n$ flotilla.}\}
$$
We wish to apply the Coboundedness Principle~\ref{PropCoboundednessPrinciple} to $\P_n$, and so we
must check the sufficient conditions (a--d) given in Section~\ref{SectionCobounded}. Any bounded
subset $B$ of $X$ intersects only finitely many edge or vertex spaces of $X$, and since the depth
$n$ flotillas are pairwise disjoint unions of vertex and edge spaces, it follows that $B$
intersects only finitely many of the elements of $\P_n$; in other words, $\P_n$ satisfies the
Local Finiteness Condition (a) of Section~\ref{SectionCobounded}. If $\A \ne \A'$ are distinct depth
$n$ flotillas then, as shown in the proof of the Vertex--Edge Rigidity
Theorem~\ref{TheoremVERigidity}, the coarse intersection $X_\A\cintersect{c} X_{\A'}$ is coarsely
contained in some edge space of depth $\ge n+1$, and so $X_\A \not\ceq{c} X_{\A'}$ in $X$;
this shows that $\P_n$ satisfies the Coarse Discreteness Condition (b). The proof of the
Vertex--Edge Rigidity Theorem~\ref{TheoremVERigidity} also shows that a quasi-isometry that
coarsely respects depth zero vertex and edge spaces will also coarsely respect flotillas of depth
$n$, and so the quasi-action of $H$ coarsely respects $\P_n$, which is the Coarse Action Condition
(c). It follows that there is an action of $H$ on the set of depth $n$ flotillas so that for each
$h \in H$ and each depth $n$ flotilla $\A$ we have $h \cdot X_\A \ceq{[r]} X_{h \cdot \A}$, where
$r\ge 0$ is independent of $h$ and $\A$; this is the Action Condition (d). The Coboundedness
Principle~\ref{PropCoboundednessPrinciple} therefore applies, and so there are only finitely many
$H$-orbits of depth $n$ flotillas $\A$, and the subgroup $\Stab(\A)=\Stab(X_\A)$
quasi-acts properly and coboundedly on $X_\A$. This quasi-action coarsely respects the vertex and
edge spaces in $X_\A$, because it is the restriction of the quasi-action of $H$ on $X$---that is,
there is an induced cobounded quasi-action of $\Stab(\A)$ on the tree~$\A$. 

We now prove Tree Rigidity by induction. In the basis step, we replace $F_0$ by a forest $F'_0$ on
which $H$ acts. In the induction step, assuming we have replaced $F_n$ by a forest $F'_n$ on which
$H$ acts, we then replace $F_{n+1}$ by a forest $F'_{n+1}$ on which $H$ acts. The replacements must
always be realized by $H$-almost equivariant quasi-isometries.

Consider a depth zero flotilla $\A$ of $T$. By Proposition~\ref{PropTrichotomy}, $\A$ is either a
point, a line, or a bushy tree with a Cantor set of ends. By the hypotheses of Tree Rigidity,
$\A$ is not a line. When $\A$ is a point then we trivially have an action of $\Stab(\A)$. Consider
now the case where $\A$ is bushy, to which the results of \cite{MSW:QTOne} apply. From
Theorem~\ref{homogeneous} we may conclude that there is an action of $\Stab(\A)$ on a tree $\A'$, and
a $\Stab(\A)$-almost equivariant quasi-isometry $\A \to \A'$. But we need a little more: from
Proposition~16 of \cite{MSW:QTOne} the group $\Stab(\A)$ acts properly discontinuously and
cocompactly on a Bass-Serre tree of spaces $X_{\A'} \to \A'$ with quotient graph of groups
$\G_{\A'}$, and so we get an identification $\Stab(\A)\approx \pi_1 \G_{\A'}$; moreover, the induced
$\Stab(\A)$-almost equivariant quasi-isometry $X_\A\to X_{\A'}$ coarsely respects vertex and edge
spaces, inducing a $\Stab(\A)$-almost equivariant quasi-isometry $\A \to \A'$.

In each $H$-orbit of depth zero flotillas, pick a representative $\A_0$, and apply the above
argument to get a graph of groups $\G_{\A'_0}$ with Bass-Serre tree of spaces $X_{\A'_0} \to \A'_0$
and an identification $\Stab(\A_0) \approx \pi_1\G_{\A'_0}$ as above. We can now propogate the
construction around by the action of $H$, as follows. For each representative depth zero flotilla
$\A_0$ and any flotilla $\A \ne \A_0$ in the orbit of $\A_0$, choose $h \in H$ so that $h \cdot
\A_0 = \A$, and define $X_{\A'} \to \A'$ to be a disjoint copy of $X_{\A'_0} \to \A'_0$ with
$\Stab(\A) = h \cdot \Stab(\A_0) \cdot h^\inv$ acting properly discontinuously and cocompactly.
The union of the trees $\A'$, over each depth zero flotilla $\A$ of $X$, defines the forest
$F'_0$. The union of the trees of spaces $X_{\A'}$ defines the ``forest of spaces'' $X_{F'_0}$.
The action of each group $\Stab(\A)$, together with the maps $h\from \Stab(\A_0) \to \Stab(\A)$
where $h \cdot \A_0 = \A$, combine to give a well-defined action of $H$ on the forest of spaces
$X_{F'_0} \to F'_0$. For the moment we will not extend the metrics on the individual trees of
spaces $X_{\A'}$ to a metric on the forest of spaces $X_{F'_0}$, preferring instead to wait until
the entire tree of spaces $X_{T'} \to T'$ has been constructed. This completes the basis step.

For the induction step, we assume that the forest of spaces $X_{F_{n-1}} \to F_{n-1}$, on which $H$
quasiacts properly and coboundedly, coarsely respecting vertex and edge spaces, has been replaced
by a forest of spaces $X_{F'_{n-1}} \to F'_{n-1}$ on which $H$ acts properly discontinuously and
cocompactly, and there is an $H$-almost equivariant map $f \from X_{F_{n-1}} \to X_{F'_{n-1}}$
which preserves the decomposition into trees of spaces, restricting to a quasi-isometry between
corresponding trees of spaces that coarsely respects vertex and edge spaces within each tree of
spaces; it follows that there is an induced $H$-almost equivariant map $f_\# \from F_{n-1}\to
F'_{n-1}$ which coarsely respects the decomposition into trees, restricting to a quasi-isometry
between corresponding trees.

Consider a depth $n$ flotilla $\A$, a component of $F_n$. The group $\Stab(\A)$ quasi-acts properly
and coboundedly on $X_\A$, coarsely respecting the vertex and edge spaces in $X_\A$. Let $\F =
F_{n-1} \intersect \A$, the union of depth $n-1$ flotillas contained in $\A$. We obtain a forest
of spaces $X_\F \to \F$. From the induction hypothesis, by restricting to $\A$ we obtain a forest
of spaces $X_{\F'} \to \F'$ on which $\Stab(\A)$ acts properly discontinuously and cocompactly, and
a $\Stab(\A)$-almost equivariant map $f \from X_\F \to X_{\F'}$, that respects the decomposition
into trees of spaces, is a quasi-isometry on each tree of spaces, and that respects the vertex and
edge  spaces within each tree of spaces. There is an induced $\Stab(\A)$-almost equivariant map
$f_\# \from \F \to \F'$ that respects the decomposition into trees and is a quasi-isometry on each
tree.

Our remaining goal is to construct a tree of spaces $X_{\A'} \to \A'$ on which $\Stab(A)$ acts
properly discontinuously and cocompactly, and a $\Stab(\A)$ almost equivariant quasi-isometry
$X_\A \to X_{\A'}$ that coarsely respects vertex and edge spaces, descending to a $\Stab(\A)$
almost equivariant quasi-isometry $\A \to \A'$. Our construction will have the property that the
depth $n-1$ subforest of $\A'$ is $\Stab(\A)$ equivariantly identified with $\F'$, and that the
restriction of $X_{\A'}$ to $\F'$ is equivariantly identified with $X_{\F'}$. This will complete
the induction step, because we can carry this construction out for a representative of each orbit
of the action of $H$ on the collection of deptn $n$ flotillas, and then propogate the construction
around by $H$-equivariance.

Consider the tree $\A / \F$ obtained by collapsing each component $F$ of $\F$ to a point that
we shall denote $v_F$. The composition $X_{\A/\F} = X_\A \to \A \to \A/\F$ defines a tree of spaces
$X_{\A/\F}\to\A/\F$ with the property that every edge has depth one. Since the quasi-action of
$\Stab(A)$ on $X_\A$ coarsely respects the spaces $X_F$ for each depth $n-1$ flotilla $F \subset
\A$, as well as the depth $n$ vertex and edge spaces in $X_\A$, it follows that $\Stab(A)$
coarsely respects the depth zero vertex spaces $X_{v_F}$ of $X_{\A/\F}$, as well as the depth one
vertex and edge spaces of $X_{\A/\F}$. We have now verified the hypotheses of
Proposition~\ref{PropDepthOne}, which is the special case of Tree Rigidity in which all edges have
depth one. Applying that proposition, we obtain a tree of spaces that we shall denote $Y \to U$,
on which $\Stab(A)$ acts properly discontinuously and cocompactly, whose quotient graph of groups
is finite type, of depth one, and irreducible, and we obtain a $\Stab(A)$-almost equivariant
quasi-isometry $f \from X_{\A/\F} \to Y$ that coarsely respects vertex and edge spaces, descending
to a quasi-isometry $f_\# \from \A/\F \to U$, so that $f(X_a) \ceq{[C]} Y_{f_\#(a)}$ for each edge
or vertex $a$ of $X_{\A/\F}$, where $C$ is a constant independent of $a$; we may also choose
coarse inverses $\bar f \from Y \to X_{\A/\F}$ and $\bar f_\# \from U \to \A/\F$ with similar
properties. Since the depth zero vertex spaces of $X_{\A/\F}$ are pairwise coarsely inequivalent,
the same is true in $Y$, and so $f_\#$ restricts to a bijection of depth zero vertices whose
inverse is the restriction of $\bar f_\#$. 

Now we apply the techniques of Proposition~\ref{PropNoTracks}, which is the special case of Tree
Rigidity in which the edge spaces are isolated. Of course, in the present setting we do not have
the isolation property for edge spaces of the tree of spaces $Y \to U$. However, the
isolation property for edge spaces was designed, somewhat artificially, for a particular purpose:
so that the quasi-action \emph{strictly} respects edge spaces, and hence descends to a true
\emph{action} on the Bass-Serre tree. In the present context, the action of the group $\Stab(A)$
on $Y$ does indeed strictly respect vertex and edge spaces, thereby descending to an action on
the quotient tree $U$. 

Our goal now is the pretty much the same as in Proposition~\ref{PropNoTracks}, and can be expressed
in the present context as follows. Consider a depth zero vertex $v_F$ of $\A/\F$ corresponding to a
component $F$ of $\F$, let $v'$ be the corresponding vertex of $U$, and let $F'$ be the
corresponding component of $\F'$. We blow up the vertex $v'$, replacing it by the tree $F'$, and
for each edge $e'$ of $U$ incident to $v'$ we must assign a vertex of $F'$ to which $e'$ should be
attached. Let $\Stab^*(e')$ be the subgroup of $\Stab(A)$ which fixes $e'$ and preserves
orientation on $e'$. Let $e = \bar f_\#(e')$, an edge of $\A/\F$, and so we may also regard $e$ as
an edge of $\A$. Although we do not know whether $e$ has an endpoint in the tree $F$, we do know
that there is a vertex $w \in F$ which is uniformly close to $e$ in the tree $F$, independent of
the choice of $v_F$ and of $e'$. Since $\Stab^*(e')$ fixes $e$, it follows that, under the
quasi-action of $\Stab(A)$ on $A$, the $\Stab^*(e')$ orbit of $e$ is uniformly bounded,
independent of $v_F$ and of $e'$. This implies in turn that, under the quasi-action of
$\Stab^*(e')$ on $F$, the orbit of $w$ is uniformly bounded, which implies further that, under the
action of $\Stab^*(e')$ on $F'$, the orbit of $f_\#(w)$ is uniformly bounded. Within the convex
hull of this orbit there exists a point $w'$ of $F'$ that is fixed by the action of $\Stab^*(e')$.
Subdividing $F'$ at $w'$ if necessary, we can attach $e'$ to $w'$.

These attachments must be made in a $\Stab(A)$-equivariant manner, and we do this by choosing a
representative pair $(v',e')$ among the $\Stab(A)$-orbits of adjacent pairs of vertices and edges
of $U$, assigning the attaching point for $e'$ in $F'$, and then propogating this assignment
around by the action of $\Stab(A)$. This produces the desired tree $\A'$. The tree of spaces
$X_{\A'}$ has vertex and edge spaces of depth $\le n-1$ taken from the vertex and edge spaces of
$X_{F'}$, and $X_{\A'}$ has vertex and edge spaces of depth $n$ taken from the depth one vertex
and edge spaces of $Y$. The attachments of edge spaces to vertex spaces may clearly be made in a
$\Stab(\A)$ equivariant manner. Since $\Stab(\A)$ quasi-acts properly and coboundedly on $X_\A$,
and since it acts properly and coboundedly on $X_{\A'}$, we obtain a $\Stab(\A)$ almost
equivariant quasi-isometry $X_\A \to X_{\A'}$ which, by construction, coarsely respects vertex and
edge spaces.

This completes the proof of the Tree Rigidity Theorem~\ref{TheoremTreeRigidity}.

\vfill\break

\section{Main Theorems}
\label{MainTheorems}

With the machinery developed in the previous sections, we can now prove our main theorems. 

\paragraph{Quasi-isometric rigidity.} Given a graph of
groups $\G$ with Bass-Serre tree of spaces $X \to T$, recall the main hypotheses: 
\begin{itemize}
\item[(1)] $\G$ is finite type,
irreducible, and finite depth.
\item[(2)] No depth zero raft of the Bass-Serre tree $T$ is
a line.
\item[(3)] Each depth zero vertex group is coarse $\PD$.
\item[(4)] The crossing graph condition holds for each
depth zero vertex of $T$ which is a raft.
\item[(5)] Each vertex and edge group of $\G$ is coarse finite type. 
\end{itemize}

\begin{theorem*}[\ref{TheoremQI} Quasi-isometric rigidity theorem]
Let $\G$ be a graph of groups satisfying (1)--(5) above. If $H$ is a finitely generated group quasi-isometric
to $\pi_1\G$ then $H$ is the fundamental group of a graph of groups satisfying (1)--(5). 
\end{theorem*}

\begin{proof} Since $H$ is quasi-isometric to $\pi_1(\G)$, it follows that $H$ quasi-acts properly and
coboundedly on $X$. By Theorems \ref{TheoremVertexRigidity} and \ref{TheoremVERigidity}, this $H$ quasi-action
coarsely respects the vertex and edge spaces. Thus we can apply Theorem \ref{TheoremTreeRigidity} to get a
graph of groups decomposition of $H$ with the desired properties. 
\end{proof}

\paragraph{Quasi-isometric classification.} Given a group $G$ as above, Theorem \ref{TheoremQI} reduces the
study of the quasi-isometry class of $G$ to studying groups which have splittings similar to that of $G$. The
content of Theorem \ref{TheoremClasses} is that all such splittings of $H$ look coarsely the same as the
splitting of $G$: 

\begin{theorem*}[\ref{TheoremClasses} Quasi-isometric classification theorem] \quad
Let $\G,\G'$ be graphs of groups satisfying (1)--(5) above. Let $X \to T$, $X' \to T'$ be Bass-Serre trees of
spaces for $\G,\G'$, respectively. If $f \from X \to X'$ is a quasi-isometry then $f$ coarsely respects vertex
and edge spaces. To be precise, for any $K \ge 1$, $C \ge 0$ there exists $K'\ge 1$, $C' \ge 0$ such that if $f
\from X \to X'$ is a $K,C$ quasi-isometry then there exists a $K',C'$ quasi-isometry $f_\# \from \VE(T) \to
\VE(T')$ such that the following hold:
\begin{itemize}
\item If $a \in\VE(T)$ then $d_\Haus(f(X_a),X'_{f_\#(a)}) \le C'$.
\item If $a'\in \VE(T')$ then there exists $a \in \VE(T)$ such that $d_\Haus(f(X_a),X'_{a'})
\le C'$.
\end{itemize}
\end{theorem*}

\begin{proof}
This is an immediate consequence of the Depth Zero Vertex Rigidity
Theorem~\ref{TheoremVertexRigidity} and the Vertex--Edge Rigidity
Theorem~\ref{TheoremVERigidity}: the hypotheses of Theorem~\ref{TheoremVertexRigidity}
hold, and so we conclude that $f$ coarsely respects depth zero vertex spaces;
theorem~\ref{TheoremVERigidity} then applies to show that $f$ respects vertex and edge
spaces of all depths as claimed.
\end{proof}

\paragraph{Graphs of coarse $\PD$ groups.} In many applications, Theorems \ref{TheoremClasses} and
\ref{TheoremQI} are combined: any group $H$ quasi-isometric to $\pi_1 \G$ has a graph of groups
decomposition $\G'$ with properties similar to $\G$; and for any graph of groups decomposition $\G'$ of $H$
with properties similar to $\G$, any quasi-isometry $\pi_1 \G \approx \pi_1 \G' \approx H$ coarsely respects
vertex and edge spaces.

In most of our applications, the vertex and edge spaces are coarse $\PD$ spaces. Such graphs of groups are
automatically finite depth, because an infinite index subgroup of a $\PD(n)$ groups has coarse dimension at
most $n-1$; see Lemma~\ref{LemmaBigInSmall}.

\begin{theorem*}[\ref{PDgraphs} QI rigidity and classification for graphs of coarse $\PD$ groups] Let $\G$ be
a finite, irreducible graph of groups. Suppose that all vertex and edges groups of $\G$ are coarse \Poincare\
duality, no depth zero raft of the Bass-Serre tree
$T$ is a line, and that the crossing graph of any one vertex depth zero raft of $T$ is either connected or
empty. If
$H$ is any finitely generated group quasi-isometric to $\pi_{1}\G$ then $H$ is the fundamental group of a
finite type, irreducible graph of groups $\G'$, each of whose vertex groups and edge groups is coarse
\Poincare\ duality, no depth zero raft of the Bass-Serre tree $T'$ is a line, and the crossing graph of any
one vertex depth zero raft is either connected or empty. Furthermore, any quasi-isometry $\pi_1 \G \to
\pi_1\G'$ coarsely respects vertex and edge spaces.
\end{theorem*}

\begin{proof} This follows immediately from Theorems \ref{TheoremClasses} and \ref{TheoremQI} and the
observation above that a graph of coarse $\PD$ groups is finite depth. 
\end{proof}

It is also worth noting that for a graph of $\PD$ groups, the assumption that no depth zero raft is a line
is essentially just a normalization condition for the splitting: any such raft can be ``rolled-up'' to a single
vertex space to give a different splitting, and the new vertex group created in this way is a coarse $PD$
group of one dimension higher. Since the edges incident to this new vertex are the edges incident to the
line raft in the old graph of groups, these edge spaces have codimension $\ge 2$ in the new vertex space, so
the new vertex has empty crossing graph. The irreducibility condition is also a kind of normalization, as
explained in Section~\ref{SectionIrreducible}: for any group $G$, any finite depth graph of groups
decomposition of $G$ can be converted into an irreducible one, and the resulting graph of groups is canonical
in the sense that the commensurability lattice of vertex and edge stabilizers in the Bass-Serre tree is
independent of the conversion process. Theorem~\ref{PDgraphs} therefore applies to \emph{any} finite graph of
$\PD$ groups, subject to a translation of the crossing graph condition into a condition that is applied to
each bounded depth zero raft of the Bass-Serre tree; making this translation is somewhat tedious, and so it is
easier to state Theorem~\ref{PDgraphs} in the normalized version above.
\vfill\break


\section{Applications and Examples}
\label{Applications}

\newcommand\g{g}

The general theorems of the previous section can be applied to a wide variety of graphs of groups.  We give
several specific examples, chosen with an eye towards applying additional rigidity tools to get even stronger
conclusions. These stronger conclusions arise from considerations of ``pattern rigidity'' which we now discuss.

The idea behind pattern rigidity is that the conclusions of Theorems~\ref{TheoremClasses} and~\ref{TheoremQI}
contain more information than is immediately obvious.  Given $G$ and $H$ as in the theorems, any quasi-isometry
coarsely respects the $\VE$ lattice.  Consider a depth zero vertex space $X_v$ for $G$ and the corresponding
vertex space $X'_{v'}$ for $H$. Not only does the quasi-isometry $G \to H$ restrict to give a quasi-isometry
$X_v \to X'_{v'}$, but this quasi-isometry respects the coarse intersections of $X_v$ and $X'_{v'}$ with other
vertex spaces. This is often enough to prove that the quasi-isometry $X_v \to X'_{v'}$ induces an isomorphism
between the pattern of edge spaces incident to $X_v$ and the pattern of edge spaces incident to $X'_{v'}$. One
then invokes one of several ``pattern rigidity'' theorems to reach even stronger conclusions about the
quasi-isometry $X_v \to X'_{v'}$.

First we discuss patterns in a very general context, before bringing in various pattern rigidity results.

\subsection{Patterns of edge spaces in a vertex space.}
\label{patterns} 
Given a space $X$, a \emph{pattern} in $X$ is a collection $P$ of subsets of $X$. If $P$ is a pattern in
$X$ and $Q$ is a pattern in $Y$, a quasi-isometry $\phi:X \to Y$ is said to \emph{respect the patterns} $P$ and
$Q$ if there is  $C>0$ so that: for all $A \in P$ there exists
$B \in Q$ for which 
$$d_\H(\phi(A),B) < C
$$ 
and likewise, for all $B \in Q$ there exists $A \in P$ for which 
$$d_\H(\phi(A),B) < C
$$

There are two important and closely related examples of concern to us: patterns of subgroup cosets in a group;
and patterns of edge spaces in a vertex space.

Let $G$ be a finitely generated group and $H_1, \cdots, H_n \subgroup G$ a finite collection of subgroups.
The corresponding \emph{coset pattern} $P$ is simply the collection of all left cosets of the $H_i$. Let $X$ be
a proper, geodesic metric space on which $G$ acts properly and coboundedly, let $Y_1, \cdots, Y_n \subset
X$ be subspaces such that $Y_i$ is invariant and cobounded under the action of the subgroup $H_i$, and let $Q
= \{g \cdot Y_i \suchthat g \in G, i=1,\ldots,n\}$ be the pattern of left translates. Any $G$-equivariant
quasi-isometry $G \mapsto X$ respects the patterns $P$ and $Q$. 

The situation above applies to graphs of groups as follows. Consider a finite type graph of groups $\G$ with
Bass-Serre complex $X \to T$, and let $p \from T \to \G$ be the quotient map. Given a vertex $v \in T$, let
$Q$ be the \emph{pattern of incident edge spaces} in $X_v$, consisting of the images $Y_e$ of the attaching
maps $X_e\to X_v$, over all edges $e$ incident to $v$. Alternatively, choosing edges $e_1,\ldots,e_n$
representing the $\Stab(v)$ orbits of edges of $T$ incident to~$v$, $Q$ is just the left
translate pattern of the subsets $Y_{e_1},\ldots,Y_{e_n}$. If $e_1,\ldots,e_n$ are chosen
appropriately, the injections $\Stab(e_1),\ldots,\Stab(e_n) \inject \Stab(v)$ correspond via the
isomorphism $\Stab(v) \to \G_v$ to the edge-to-vertex injections
$\G_{p(e_1)},\ldots,\G_{p(e_n)} \inject \G_{p(v)}$, whose images define a coset pattern $Q$ in
$\G_{p(v)}$. It follows that any $\G_{p(v)}$ equivariant quasi-isometry $\G_{p(v)}\to X_v$
respects the patterns $P$ and $Q$.

In many cases the group of quasi-isometries of a group preserving a coset pattern is quite small. As we will
see in the next several sections, this phenomenon, called \emph{pattern rigidity}, leads to very strong
rigidity results for many classes of graphs of groups.

\subsection{$\hyp^n$ vertex groups and $\Z$ edge groups.} 
\label{SectionHnZ}
The first example of a pattern rigidity 
theorem is the result of Schwartz (\cite{Schwartz:Symmetric}) on patterns of geodesics in
$\hyp^n$ for $n \geq 3$.  Schwartz proves that if $M$ is a compact hyperbolic $n$-manifold, $n \ge 3$,
and $c_1, \cdots ,c_k$ is any finite collection of closed geodesics in $M$, and $\C$ is the pattern of
geodesics in $\hyp^n$ covering the $c_i$, then any quasi-isometry of $\hyp^n$ that coarsely preserves
$\C$ is a bounded distance from an isometry that preserves $\C$. It follows that the subgroup of
$\QI(\hyp^n)$ that coarsely preserves $\C$ contains $\pi_1(M)$ as a subgroup of finite index.

\begin{theorem}\label{GenHyp} Let $\G$ be a finite graph of groups whose vertex groups are fundamental groups
of compact hyperbolic $n$-manifolds, $n \ge 3$, and whose edge groups are infinite cyclic. If $H$ is a
finitely generated group quasi-isometric to $G=\pi_1(\G)$ then $H$ splits as a graph of groups whose depth
zero vertex groups are commensurable to those of $\G$ and whose edge groups and positive depth vertex groups
are virtually infinite cyclic.   
\end{theorem}

\begin{proof}
The vertex and edge spaces are coarse $\PD$ spaces, all depth zero rafts are one vertex rafts,
and there are no codimension one edge-to-vertex inclusions.  Thus theorem \ref{PDgraphs}
applies to show that $H$ splits as a graph of groups $\G'$ with depth zero vertex spaces
quasi-isometric to $\hyp^n$ and edge groups quasi-isometric to $\Z$.   

By construction the quasi-isometry respects the vertex and edge spaces of this splitting, and thus the
quasi-actions of the vertex groups on the vertex spaces of $\G$ preserve the patterns of edge spaces. 
Schwartz' theorem \cite{Schwartz:Symmetric} thus shows that the depth zero vertex groups in the splitting
$\G'$ are commensurable to the corresponding groups in the splitting $\G$.
\end{proof}

This means that to understand all groups quasi-isometric to $G$, one only needs to study graphs
of groups of the given type.  One can read off many quasi-isometry invariants from the
splittings.  One simple example is:

\begin{corollary}\label{HypManifolds} Let $G_i= A_i *_{\Z} B_i$ for $i=1,2$ be amalgamated free
products with $A_i,B_i$ fundamental groups of closed hyperbolic manifolds of dimension $\ge 3$,
with $\Z$ including into $A_i$ and $B_i$ as maximal cyclic subgroups. If $G_1$ and $G_2$ are
quasi-isometric then up to reordering, $A_1$ is commensurable to $A_2$ and $B_1$ is
commensurable to $B_2$.  Further, these commensurators must respect, up to finite index, the
cyclic subgroups.  Thus, in particular, for fixed $A$ and $B$ there are infinitely many
quasi-isometry classes of groups of the form $A*_{\Z} B$. 
\end{corollary}

\begin{proof}
The Classification Theorem shows that any quasi-isometry between $G_1$ and $G_2$ coarsely
respects the given splittings.  The commensurability statements then follow from the same
argument as in the proof of Theorem \ref{GenHyp}.   

To see the infinitude of quasi-isometry types, one needs to check that the fundamental group of
any compact hyperbolic $n$-manifold contains infinitely many cyclic subgroups not commensurable
to each other. This follows as any such manifold has infinitely many distinct primitive
closed geodesics.
\end{proof}


\paragraph{A variation: $\hyp^n$ vertex groups and free edge groups.} The proof of Theorem~\ref{GenHyp}
suggests other classes of groups to investigate. For example, consider an amalgamated free product $G_1 *_F
G_2$ where $G_1,G_2$ are fundamental groups of closed hyperbolic manifolds of dimension $\ge 3$, and $F$ is a
finite rank free group including into each of $G_1$ and $G_2$ as a malnormal subgroup. The finite depth
property follows from the malnormality hypothesis, and so our main QI-rigidity and classification theorems
apply to this group. It follows that if the finitely generated group $H$ is quasi-isometric to $G_1 *_F G_2$
then $H = \pi_1\G$ where $\G$ is a finite graph of groups whose vertex groups are quasi-isometric to either
$G_1$ or $G_2$ and whose edge groups are quasi-isometric to $F$. The edge groups are therefore virtually free,
and each vertex group has a finite normal subgroup so that the quotient is a cocompact lattice in the
isometries of hyperbolic space of some dimension. However, in this situation we do not know of any pattern
rigidity result which would imply that the vertex groups of $\G$ are commensurable to
$G_1$ or $G_2$. 

This invites some questions: Is there a pattern rigidity theorem for a free subgroup $F$ of a cocompact
lattice $G$ acting on hyperbolic space? What if $F$ is quasiconvex? What if $F$ is quasiconvex and malnormal?
In the latter situation, the limit sets of all conjugates of $F$ in $G$ would form a pairwise disjoint pattern
of Cantor sets in the boundary of hyperbolic space, which might be an interesting starting point for pondering
a generalization of Schwartz' pattern rigidity theorem \cite{Schwartz:Symmetric}.

\subsection{$\hyp^3$ vertex groups and surface fiber edge groups.} 
In any context in which one has a
pattern rigidity result analogous to the theorem of Schwartz used above, one can substantially
strengthen the general results of Theorems~\ref{TheoremQI} and~\ref{TheoremClasses}. Here is a pattern
rigidity result of Farb and Mosher for fibered hyperbolic \nb{3}manifolds; the proof of this result
ultimately depends on Schwartz' result \cite{Schwartz:Symmetric}.

\begin{theorem}[\cite{FarbMosher:sbf}] 
\label{TheoremSBF}
Let $M$ be a hyperbolic $3$-manifold which fibers over
the circle with fiber $F$. The subgroup of $\QI(\hyp^3)$ which coarsely preserves the pattern
of cosets of $\pi_1(F)$ is a subgroup $\G$ of $\Isom(\hyp^3)$ which contains $\pi_1(M)$ as a
finite index subgroup.  
\end{theorem}

We can use this to prove a rigidity theorem for graphs of groups. Given a closed hyperbolic
\nb{3}manifold $M$, a \emph{fiber subgroup} of $\pi_1(M)$ is the image of an injection $\pi_1(F) \inject
\pi_1(M)$ where $M$ fibers over $S^1$ with fiber $F$. 

\begin{theorem}\label{FiberedManifolds} Let $\G$ be a finite graph of groups in which every
vertex group is the fundamental group of a closed hyperbolic $3$-manifold, every edge group is
a surface group, and the image of every edge-to-vertex injection is a fiber subgroup. If $H$ is
any finitely generated group quasi-isometric to $\pi_1(\G)$ then $H$ splits as a graph of
groups whose depth zero vertex groups are commensurable to those of $\G$, and whose edge groups
are positive depth vertex groups are commensurable to surface groups.
\end{theorem}

\begin{proof} The case where $\G$ has no edges was proved by Gromov and Sullivan; see
\cite{CannonCooper}. We may therefore assume $\G$ has at least one edge. 

We modify the graph of groups $\G$ as follows. Suppose that $v$ is a vertex such that all the
incident edge groups are contained in the same fiber subgroup $\pi_1(F) \inject \pi_1(M)
\approx \G_v$. Blow up the vertex $v$, replacing it by a vertex $v'$ with group $\G_{v'} \approx
\pi_1(F)$ and an edge group $\pi_1(F)$ glued at both ends to $\G_{v'}$, thereby describing
$\pi_1(M)$ as an HNN amalgamation with $\pi_1(F)$ as fiber.  Call the resulting graph of groups
$\G'$. This is a graph of $\PD$ groups, and as $\G$ has at least one edge, $\G'$ has no
line-like rafts. To apply Theorem \ref{PDgraphs} we need to check the crossing graph condition.
If all vertex groups of $\G'$ are surface groups then all crossing graphs are empty. Those
vertices which are not surfaces are hyperbolic $3$-manifold groups with at least two distinct
fiber subgroups containing edge groups with finite index. 

\begin{claim} If $M$ is a closed hyperbolic $3$-manifold which fibers
over the circle in two homotopically distinct ways, with fibers $F$ and $F'$, then
every coset of $\pi_1(F)$ crosses every coset of $\pi_1(F')$.
\end{claim} 

Since $\pi_1(F)$ and $\pi_1(F')$ are normal subgroups, all the cosets of
each are parallel, so we only need to check that $\pi_1(F)$ crosses
$\pi_1(F')$.  This is immediate as a subset of $\pi_1(M)$ which stays to
one side of $\pi_1(F)$ is contained in a half ray in $\pi_1(M)/\pi_1(F)
\equiv \Z$.  For a subgroup, this means the image must be trivial, in
other words, that the subgroups of $\pi_1(M)$ which do not cross
$\pi_1(F)$ are the subgroups of $\pi_1(F)$. Since $F$ and $F'$ are
non parallel fibers, they cross, proving the claim.

We can now apply Theorem \ref{PDgraphs} to produce the desired splitting of $H$. The statement
about commensurability types of depth zero vertex groups follows from
Theorem~\ref{TheoremSBF}, just as the similar statement in Theorem~\ref{GenHyp} followed
from Schwartz' pattern rigidity. This finishes the proof of Theorem~\ref{FiberedManifolds}.
\end{proof}

The Classification Theorem together with pattern rigidity and Theorem~\ref{TheoremSBF} can
be used to give many quasi-isometry invariants of such groups, as in corollary
\ref{HypManifolds}.  Here is one example:

\begin{corollary}  Let $M_i$, $i=1,2$, be closed hyperbolic $3$-manifolds each of which which
fibers over the circle in two distinct ways, with fibers $F_i$ and $F'_i$.  Let $T_i$ be an
isomophism between a finite index subgroup of $\pi_1(F_i)$ and $\pi_1(F'_i)$, and let $G_i$ be
the $HNN$ extension of $\pi_1(M_i)$ over $T_i$.    If $G_1$ and $G_2$ are quasi-isometric then
there is an isomorphism between finite index subgroups of $\pi_1(M_1)$ and $\pi_1(M_2)$ which
(up to re-ordering) also commensurates $\pi_1(F_1)$ with $\pi_1(F_2)$ and  $\pi_1(F'_1)$ with
$\pi_1(F'_2)$.
\end{corollary}

\subsection{Surface vertex groups and cyclic edge groups.}

For surfaces, connected crossing graph has a geometric interpretation. Consider first a closed
hyperbolic surface $\Sigma$ and a collection of nontrivial conjugacy classes $c_1,\ldots,c_k$ in
$\pi_1 \Sigma$. These classes are represented by simple closed geodesics in $\Sigma$, also
denoted $c_1,\ldots,c_k$. We say that $c_1,\ldots,c_k$ \emph{fill} $\Sigma$ if the components of
$\Sigma-(c_1\union\cdots\union c_k)$ are all open discs. More generally, suppose that $\Sigma$
is a closed hyperbolic \nb{2}orbifold, and $c_1,\ldots,c_k$ are infinite order conjugacy
classes in $\pi_1\Sigma$. Each conjugacy class $c_i$ is still represented by a closed geodesic
in $\Sigma$, properly interpreted: $c_i$ can reflect across a mirror edge of the orbifold; and
$c_i$ can pass through a cone point or a corner reflector. We say that $c_1,\ldots,c_k$ fill
$\Sigma$ if each component of $\Sigma - (c_1\union\cdots\union c_k)$ is a 2-orbifold with
finite fundamental group. In both cases, $c_1,\ldots,c_k$ fill $\Sigma$ if and only if the
union of the total lifts of $c_1,\ldots,c_k$ to $\hyp^2$ has all of its complementary
components being open discs with compact closure.

\begin{lemma}  Let $\Sigma$ be a closed hyperbolic $2$-orbifold, $G=\pi_1(\Sigma)$.   Let
$c_1,\cdots, c_k$ be a collection of non-trivial isotopy classes of closed geodesics on $\Sigma$,
and let $H_i$ be an infinite cyclic subgroup generated by an element of $G$ in the conjugacy
class determined by  $c_i$.  The crossing graph of the collection of conjugates of the $H_i$ in
$G$ is connected if and only if the curves $c_i$ fill $\Sigma$.
\end{lemma}

\begin{proof}
Let $C$ be the collection of geodesics in $\hyp^2$ overing the $c_i$.  Note that  $C$ is in
one-to-one correspondence with the collection of subgroups of $G$ conjugate to the $H_i$. 
Further, two such subgroups cross in the sense of section \ref{SectionAbstractVertexRigidity}
if and only if the corresponding geodesics intersect in $\hyp^2$. Thus the crossing graph is
connected if and only if, $Z$, the union of geodesics in $C$, is a connected subset of
$\hyp^2$.  The lemma therefore says that $Z$ is connected if and only if the complementary
components of $Z$ are bounded.  

Suppose $Z$ is connected. Let $U$ be a component of $\hyp^2 - Z$. Since the collection $C$
is locally finite, the stabilizer of $U$ in $G$ acts coboundedly  on $U$. If $U$ is unbounded
this implies there is a non-trivial element of $G$ preserving $U$, and hence, as $U$ is convex,
$U$ conatains the axis of an element of $G$.  This axis is disjoint from all elements of $C$ and
disconnects $\hyp^2$. Since $Z$ is connected it must lie on one side of this axis. This is a
contradiction as $Z$ is $G$ invariant and $G$ is cocompact on $\hyp^2$. Thus all complementary
components are bounded.

Next, suppose all components of $\hyp^2 - Z$ are bounded.  If $Z$ is not connected then there
must be a component $U$ of this complement which borders two components of $Z$.  However, since
$U$ is disk, its boundary is connected, which is a contradiction. 
\end{proof}

Schwartz's pattern rigidity result does not hold in general for patterns of geodesics on
hyperbolic surfaces.  However, Kapovich and Kleiner show that pattern rigidity does hold when the geodesics
$c_1,\cdots, c_k$ fill the surface; see \cite{KapovichKleiner:LowDBoundaries} Lemma~22.  This is equivalent to
our assumption on crossing graphs, and so we have:

\begin{theorem} 
\label{TheoremStrongSurface}
Let $\G$ be a graph of groups where every vertex group is the fundamental group of a closed
hyperbolic $2$-orbifold, and whose edge groups are virtual infinite cyclic groups.  Assume that
the edge groups at each vertex correspond to curves which fill.  If $H$ is any finitely
generated group quasi-isometric to $G = \pi_1(\G)$ then $H$ has a splitting with virtual
surface groups as depth zero vertex groups and virtually cyclic positive depth vertex groups
and edge groups.  Further, the quasi-isometry $G\to H$ restricts to a commensurator on each
vertex group.
\end{theorem}

Since all cocompact lattices in $SL_2(\R)$ are (abstractly) commensurable, one must look more closely at the
consequences of the Classification Theorem to see that the graphs of groups in
Theorem~\ref{TheoremStrongSurface} are not all quasi-isometric.   The key point is that it is certainly not
true that any two finite collections of filling curves is commensurable to any other.

For explicit examples, let $T_{r}$ denote the reflection group of the hyperbolic triangle with
angles $\pi/2$, $\pi/4$, $\pi/r$. This has an index~$2r$ subgroup which is the reflection group
in a regular, right angled, hyperbolic polygon with $r$ sides. Consider the tiling of $\hyp^2$
by copies of this polygon, moved around by the group $T_r$. Each side of each tile extends to a
bi-infinite geodesic which is contained in the union of sides of the tiles. This pattern of
geodesics, denoted $\S_r$, is invariant under $T_r$. Picking one geodesic in $\S_r$, the
subgroup of $T_r$ that stabilizes this geodesic is an infinite dihedral group denoted $D_r$. 

Consider the free product with amalgamations $T_r *_{D_\infinity} T_s$ where $D_\infinity$
includes to the subgroups $D_r \subset T_r$ and $D_s\subset T_s$. We claim that the set of
integers $\{r,s\}$ is a quasi-isometry invariant, thus giving infinitely many distinct
quasi-isometry classes. To prove this it suffices to show that if $r\ne r'$ then there is no
quasi-isometry of $\hyp^2$ taking the pattern $\S_r$ to the pattern $\S_{r'}$. However, any such
quasi-isometry takes two crossing geodesics in $\S_r$ to two crossing geodesics in $\S_{r'}$,
and when $r \ne r'$ there is clearly no bijection between $\S_r$ and $\S_{r'}$ that preserves
the relation of crossing.

\subsection{Graphs of abelian groups.}

For $\Z^n$ vertex groups, the crossing graph condition is simple: a subgroup $A
\subgroup\Z^n$ crosses a subgroup $B \approx \Z^{n-1}$ if and only if $A$ is not virtually
contained in $B$, if and only if $A,B$ generate a finite index subgroup of $\Z^n$. 

The following generalized version of Theorem~\ref{TheoremBabyAbelian} is an immediate corollary
of Theorem~\ref{TheoremQI}:

\begin{corollary} \label{AbelianGraphs} Let $\G$ be a finite type, irreducible graph of groups
in which all vertex and edges groups are finitely generated abelian groups.  Assume that there
are no line-like depth zero rafts, and that at any vertex $v$ which is a depth zero raft, if
there is an incident codimension one edge then the incident edge groups generate a subgroup of
finite index in $\G_v$. If $H$ is any finitely generated group quasi-isometric to $\pi_{1}\G$
then $H$ is the fundamental group of a graph of groups, each of whose vertex groups is
quasi-isometric to a vertex or edge group of $\G$ and each of whose edge groups is
quasi-isomtric to an edge group of $\G$. In particular all vertex and edge groups are finitely
generated and virtually abelian.
\end{corollary}

For stronger results, we need to understand pattern rigidity phenomena for abelian groups.  

Given a linear subspace $V \subset \reals^n$, let $P_V$ be the set of affine subspaces of
$\reals^n$ parallel to $V$. Let $F$ be a finite collection of linear subspaces of $\R^n$.  The
affine pattern induced by $F$, denoted $P_F$, is the union of $P_V$, $V \in F$. We say that the
collection $F$, or its affine pattern $P_F$, is \emph{rigid} if for every $(K,C)$ and $R$ there
is an $R'$ such that if $f \from \reals^2 \to \reals^n$ is a $K,C$ quasi-isometry, and if $f$
coarsely respects each pattern $P_V$ with coarseness constant $R'$ for each $V \in F$, then $f$
is within $R'$ of an affine homothety.

\begin{lemma} 
\label{LemmaRigidLinear}
Let $F$ be a rigid collection of subspaces of $\R^n$.  Let $f$ be a quasi-isometry of $\R^n$ such that $f$
sends $P_F$ to some affine pattern $P_{F'}$.  Then $f$ is at bounded distance from an affine map. Further, up
to homothety there are only finitely many such affine maps. 
\end{lemma}

\begin{proof}
For each $V \in F$ there is a $V' \in F'$ such that subspaces parallel to $V$ are sent by $f$ to subspaces
parallel to $V'$.  Let $f_\#$ be this induced map $F \to F'$.  We claim there is a linear  map $T \in GL_n(\R)$
such that $TV = F_\# V$ for all $V \in F$. Let $f_\lambda(v)=\frac{1}{\lambda}f(\lambda v)$.  There is a
sequence of $\lambda_i$ going to infinity such that the $f_{\lambda_i}$ converge to a bilipschitz homeomorphism
$\phi$ of $\R^n$.  For each $V \in F$, $\phi$ sends the affine foliation parallel to $V$ to the affine
foliation parallel to $f_\#(V)$. By Rademacher's theorem, $\phi$ is differentiable almost everywhere, and the
derivatives give the desired linear map.

Thus, after changing $f$ by composition with a linear map, we may assume $F=F'$ and $f_\#$ is the identity. 
Since $F$ is rigid, this $f$ is at bounded distance from a homothety. Thus the original $f$ is at bounded
distance from an affine map.  Further, rigidity implies that, up to homothety, the affine map is determined by
the map $f_\#$. Since the collections $F$ and $F'$ are finite, there are only finitely many homothety classes
of maps.  
\end{proof}

The following is an immediate consequence of Theorems~\ref{TheoremQI} and~\ref{TheoremClasses},
together with Lemma~\ref{LemmaRigidLinear} and the observation that the assumption on the
patterns of edge attachments implies the crossing graph condition: 

\begin{corollary}\label{abrigid} Let $\G$ be a graph of groups with all vertex and edge groups
finitely generated abelian groups. Assume that for each depth zero, one vertex raft
$v$ in the Bass-Serre tree, the collection of edges spaces at the vertex space of $v$ is a rigid
affine pattern. Assume also that there are no line-like rafts of depth zero.  If $H$ is any
finitely generated group quasi-isometric to $G=\pi_1(\G)$ then $H$ splits as a graph of
virtually abelian groups and the quasi-isometry $G \to H$ is affine along each depth zero, one
vertex raft.  Moreover, the set of affine equivalence classes of edge patterns in depth zero one
vertex rafts is a quasi-isometry invariant of $G$, and so this pattern is the same for $H$ as
it is for $G$.
\qed\end{corollary}

The precise determination of which collections $F$ are rigid is complicated. One class~is:

\begin{lemma} If $F$ is a finite collection of linear subspaces of $\R^n$ which contains
$n+1$ hyperplanes in general position then $F$ is rigid.
\end{lemma}

\begin{proof}
It suffices to prove this when $F$ is a collection of $n+1$ hyperplanes in general position. 
We can choose coordinates on $\R^n$ so that the first $n$ of these hyperplanes are the
coordinate hyperplanes $x_i=0$.  By scaling these coordinates, we may assume that the final
hyperplane, $H$, is given by an equation $x_1 + \cdots + x_n =0$.  Let $f$ be a quasi-isometry
which coarsely preserves the foliations parallel to hyperplanes in $F$ to within $R$.

Coarsely preserving the foliation parallel to $x_i=0$ means that there is a constant $R \ge 0$
such that for each $c \in \reals$ there is an $k \in \reals$ such that the image of the
hyperplane $x_i = c$  is contained in the set $k-R \leq x_i\leq k+R$.  Thus, by moving $f$ a
bounded distance, we can write $f$ as a product of quasi-isometries of $\R$, $f(x_1, \cdots,
x_n) = (f_1(x_1), \cdots, f_n(x_n))$.   Composing $f$ with a translation, we can assume each
$f_i(0)=0$.    The lemma now follows by induction.  Assume $k \geq 3$  and the lemma is known
for $n < k$.  Let $H_i$ be the hypersurface $x_i=0$.  We know that $f$  coarsely preserves $H_i$
and respects the pattern of foliations of $H_i$ obtained  by intersection.   By induction, this
restriction is a homothety, so there is some $C_i$  so that for $j \neq i$, $f_j (x) = C_i x$. 
By considering different $i$, this proves the lemma.

It remains to prove the case $n=2$.  We have $f=(f_1,f_2)$ a quasi-isometry of $\R^2$ which 
preserves the foliation by lines $x + y = C$.  We may assume $f(0,0) = (0,0)$.  Since $(a,0)$ and
$(0,a)$  are on a line of the foliation, the same is (coarsely) true of the images $(f_1(a),0)$
and $(0,f_2(a))$, so by moving $f$ a bounded distance, we may assume $f_1 = f_2$.  Likewise,
the point $(a,b)$ is on the line containing $(a+b,0)$ so $(f_1(a),f_1(b))$ is within bounded
distance of the line of the foliation containing $(f_1(a+b),0)$.  This means that $f_1$ is a a
quasi-morphism: there is a $C$ so that $|f_1(a+b) - f_1(a) - f_1(b)| \leq C$.  Since every
quasi-morphism $\R \to \R$ is a bounded distance from a homomorphism \cite{Gromov:Volume}, the
lemma follows.
\end{proof}

A simple example where corollary \ref{abrigid} applies is when every vertex group of $\G$ is $\Z^2$ and every
edge group is $\Z$, with at least three distinct lines of attachment at every vertex.   Such groups have been
studied and have intricate geometric and algebraic properties (\cite{BradyBridson:Isoperimetric} and
\cite{CrokeKleiner:IdealBoundaries}).  Any two patterns of three lines in the plane are affinely equivalent. 
For four lines the cross-ratio of the slopes is a complete invariant.  Our result says that these cross-ratios 
are quasi-isometry invariants of the groups.

\subsection{Quasi-isometry groups and classification}

As these example applications show, Theorems \ref{TheoremQI} and~\ref{TheoremClasses} have
strong implications for many classes of graphs of groups.  Given a class of graphs of groups
$\cal{C}$, if  Theorem \ref{TheoremQI} applies, one can conclude that any group $G$
quasi-isometric to a group in $\cal{C}$ is in $\cal{C}$.  Theorem \ref{TheoremClasses} says
that quasi-isometries between groups in $\cal{C}$ respect the splittings.  We have seen that
considering the patterns of edge spaces at a vertex often shows that this implies the
quasi-isometries between the vertex groups are, in fact,  virtual isomorphisms.  This falls
short of a complete classification of the quasi-isometry classes of groups in $\cal{C}$.   In
general this is a difficult problem which requires a method for constructing quasi-isometries
between groups in $\cal{C}$.   One case where this difficulty does not arise is for groups $G$ which are
finite index in their quasi-isometry group $\QI(G)$.  This implies that any group
quasi-isometric to $G$ is commensurable to $G$.  

Theorem \ref{TheoremClasses} can, in many cases, be used to calculate the group $\QI(G)$. 
Consider, for example, an amalgamated product $G=A *_{\Z} B$ of hyperbolic $n$-manifold groups,
as in corollary \ref{HypManifolds}. Let $X \to T$ be the Bass-Serre tree of spaces for this
splitting, each vertex space of which is isometric to $\hyp^3$. Let $f$ be a self
quasi-isometry of $G$. According to Theorem \ref{TheoremClasses}, $f$ coarsely respects vertex
spaces. For this graph of groups, no two vertex spaces are at finite Hausdorff distance, so
this means $f$ induces a bijection of the $VT$, the vertices of $T$.  Further, since the cyclic
subgroups are maximal, any two distinct edge spaces at a vertex have bounded coarse
intersection.  This implies that the edges in $T$ can be characterized by the fact that their
endpoints are pairs of vertex spaces with unbounded coarse intersection.  Thus $f$ induces an automorphism of $T$. From the discussion in corollary \ref{HypManifolds} we know that the map $f$ induces along vertex
spaces is, within bounded distance, an isometry. However, this is not enough to imply that
$f$ is itself within bounded distance of an isometry of the Bass-Serre complex. Indeed, this is
not true.  The flexibility comes from Dehn twists: given two adjacent vertex spaces $V$ and
$V'$, the isometries induced by $f$ do not necessarily agree on the edge space connecting them,
but need only be at bounded distance. Thus one can differ from the other by
an isometry which translates along this geodesic. This is essentially all that can happen. 

\begin{definition} The group of \emph{generalized Dehn twists} of $G$, $\cal{D}(G)$ is the
group of maps $\phi: VT \to G$ such that:

\begin{itemize}

\item The map $v \mapsto \phi(v) v $ is an automorphism of $T$.

\item For every edge $e$ in $ET$, the element $\delta(e) = \phi(\iota e)^{-1} \phi(\tau e)$
preserves $e$.

\item There is an $R$ such that for all $e$ translation distance of $\delta(e)$ along $X_e$ is at
most $R$. 

\end{itemize} 

\end{definition}

A generalized Dehn twist gives a self map of $G$ defined by $x \mapsto \phi(v) x$, where $v$
is the vertex for which $x \in X_v$.   This is clearly an isometry along vertex spaces.  The
third condition guarantees the map is Lipschitz, and the first condition guarantees the map is
onto.  Since the inverse is also a generalized Dehn twist, these are quasi-isometries of
$G$.   Thus there is a natural map $\cal{D}(G) \to \QI(G)$.  This map is injective since no
non-trivial isometry of $\hyp^n$ is a bounded distance from the identity. 

There is an automorphism $\tau$ of $G$ which is the identity on $A$ and conjugation by $t$ on
$B$, where $t$ is the generator of the amalgamating subgroup.  This is a Dehn twist, and it is
easy to see that the group of automorphisms generated by $G$ and $\tau$ is precisely the group
of ``periodic'' generalized Dehn twists.  Note that by \cite{RipsSela:structure} these are, up to finite
index, all the automorphisms of $G$.  

\begin{proposition} 
\label{PropDehn}
Let $G = A *_{\Z} B$ be as above, and assume that $A$ and $B$ are the full groups
of isometries of $\hyp^n$ preserving the patterns of geodesics corresponding to the inclusions
of $\Z$. The subgroup $\cal{D}(G)$ in $\QI(G)$ is of index at most $2$.
\end{proposition}

\begin{proof}
Let $f \from X \to X$ be a quasi-isometry, descending to an automorphism of $T$ also denoted
$f$. Every edge of $T$ connects a vertex covering $A$ to one covering $B$.  Since $f$ induces an
automorphism of $T$ it either preserves or reverses these two classes of vertices.  By passing
to a subgroup of index two in $\QI(G)$ we may assume $f$ preserves these classes.   

Since $f$ preserves the $G$ orbits of vertices in $T$ there is a $\g_0$ such that $f(v)=\g_0
v$.  The map ${\g_0}^{-1} f$ thus fixes $v$.  The map this induces on the vertex space $X_v$
preserves the patterns of edge spaces, and so there is an element $\g_1$ of $G^v$ which is at
finite distance from ${\g_0}^{-1}f$ along  $X_v$.    Thus $f$ is at finite distance from
$\g=\g_0 \g_1$.   If $\g'$ is any other such element then $\g^{-1} \g'$ is an element of
$G^v$, and along $X_v$ is isometry of $\hyp^n$ at bounded distance from the identity, and
therefore is the identity.  Since the map from $G^v$ to $Iso(\hyp^n)$ is injective, $\g = \g'$.

Let $\phi$ be the map $VT \to G$ sending $v$ to the $\g$ which agrees with $f$ along $X_v$. 
Since $f$ induces an automorphism of $T$, $\phi$ satisfies the first two conditions in the
definition of a generalized Dehn twist.  The third condition holds as $f$ is a
quasi-isometry.   Thus $\phi$ is a generalized Dehn twist which agrees with $f$ up to bounded
distance.
\end{proof}

A similar conclusion holds in the setting of Theorem \ref{TheoremStrongSurface}.   While this gives
a complete description of the quasi-isometry groups in these cases, none of the groups are
rigid in the sense of being finite index in their quasi-isometry groups.  Thus these
calculations do not immediately lead to a classification.  One should be able to analyze the
uniform subgroups of the quasi-isometry group, as in \cite{MSW:MaximallySymmetric} to give a
complete classification.   As is clear from the discussion, the lack of rigidity comes from the
fact that the edge groups have non-trivial centralizer.  Thus similar groups formed by gluing
$\hyp^n$ manifolds along totally geodesic immersed $m$-manifolds for $1<m<n$ should be rigid. 
However, the analogue of Schwartz' pattern rigidity result in this context does not seem to be
known.    Along the same lines, the situation for graphs of Abelian groups is much more
flexible.  Some partial classification results are known only in low dimensions
\cite{Whyte:Tubular}.

The groups of Theorem \ref{FiberedManifolds} show a different phenomenon which
causes problems in calculating the groups of quasi-isometries.  There the edge groups are
virtual surface groups, and it is certainly true that any isometry of a vertex group which is
bounded distance from the identity along an edge group is the identity (or, at least, finite
order).   Thus there are no Dehn twist type quasi-isometries.  However, since the edge groups
are normal in the vertex groups there is no reason to expect a quasi-isometry to induce an
automorphism of the Bass-Serre tree.   The difficulty is that these graphs of groups have
non-trivial depth one rafts.   These depth one rafts are surface-by-free groups, and their self quasi-isometries are
not well understood.  The results of \cite{FarbMosher:sbf} give rigidity in some cases, which
we can use to prove rigidity here.


\begin{theorem}\label{StrongFibered} Let $\Sigma$ be a closed surface of genus at least $2$.  
Suppose $\phi$ and $\psi$ are elements of $MCG(F)$ which generate a Schottky group in the sense
of \cite{FarbMosher:quasiconvex} and that the mapping tori of $\phi$ and $\psi$ are
homeomorphic to the same hyperbolic $3$-manifold $M$.  Let $G$ be the $HNN$-extension of
$\pi_1(M)$ along the two embeddings of $\pi_1(F)$.   The group $\QI(G)$ contains $G$ as a
subgroup of finite index.
\end{theorem}

\begin{proof}

Let $T$ be the Bass-Serre tree and $X$ the Bass-Serre complex of the loop of groups defining
the $HNN$-extension.  For two adjacent vertices $u$ and $v$, the vertex spaces $X_u$ and $X_v$ 
have $2$-dimensional coarse intersection.   If $u$ and $v$ are any two vertices, then their
coarse intersection is the coarse intersection of the edge spaces along the path from $u$ to
$v$.    At any vertex $w$ there are two classes of edges corresponding to  $F$ and $F'$.    
The edge spaces from two edges in the same class are at finite Hausdorff distance, while for
edges from different classes the corresponding edge spaces have coarse intersection which is at
most one dimensional.   Thus,  $X_u$ and $X_v$ have two dimensional coarse intersection if and
only if the path $u$ to $v$ enters and leaves every vertex along edges of the same class.   If
we orient the edge in the quotient graph,  and thereby the edges of the Bass-Serre tree, this
corresponds to a ``zig-zag'' path which must reverse orientation at each vertex.   

Fix an edge $e$, and let $X'$ be depth one raft corresponding to $e$ (recall that this is the collection
of vertex and edge spaces which coarsely contain $X_e$), let $T'$ be the corresponding subtree
of $T$, and let $G'$ be the subgroup of $G$ stabilizing $X'$.  From the description in
previous paragraph one sees that $T'/G'$ is a single edge.  This gives a splitting of $G'$ as
the amalgamated product $\pi_1(M) *_{\pi_1(\Sigma)} \pi_1(M)$.   The space $X'$ is the
Bass-Serre complex of this splitting.   We modify this graph of groups by replacing the vertex
groups by $\pi_1(\Sigma)$, and adding the loops at each vertex exhibiting $\pi_1(M)$ as an $HNN$
extension.  We can now collapse the edge in the middle to get a rose of surface groups whose
fundamental group is $G'$.   Thus we see that $G'$ is a surface-by-free group.   By
construction, the map from the free group to $MCG(\Sigma)$ is the Schottky subgroup generated
by $\psi$ and $\phi$.   The results of \cite{FarbMosher:sbf} give that the group of
quasi-isometries of
$G'$ is the group of isometries of the depth one raft and contains $G'$ as a finite index subgroup.  
Let $\Delta$ be a collection of coset representatives of $G'$ in $\QI(G')$.

Choose a base point $u$ in $T$.  Let $f$ be a quasi-isometry of $G$.  By Theorem
\ref{TheoremClasses}, $f$ coarsely preserves the vertex and edge spaces.  For $f$ in a
subgroup of index at most two in $\QI(G)$ we may assume $f$ preserves the two classes of edge
spaces.  After composing with an element of $G$, we may assume $f$ fixes $u$.   Thus $f$
induces a quasi-isometry of the two depth one rafts passing through $u$.  By assumption, these are
rigid.  Thus, after further composing with an element of $G$, we may assume $f$ agrees with an
element of $\Delta$ on the depth one raft. 

We claim that for each $\delta \in \Delta$ there is at most one quasi-isometry of $G$ which
preserves $X'$ and agrees with $\delta$ there.   Given any two, $f_1$ and $f_2$, the map
${f_1}^{-1} f_2$ preserves $X'$ and is at bounded distance from the identity on $X'$.  So we
need to see that any quasi-isometry of $G$ which is bounded distance from the identity on $X'$
is at bounded distance from the identity on all of $X$.   This follows by a sort of analytic
continuation argument.   Let $v$ be adjacent to $u$ along an edge $e$ (not necessarily in the
depth one raft $X'$).  Since $f$ is at bounded distance from the identity on $X_u$, it is at bounded
distance from the identity on $X_e$.   Let $X''$ be the depth one raft of $e$.   Since $f$ coarsely
preserves $e$, it coarsely preserves $X''$.  Again by \cite{FarbMosher:sbf}, $f$ must be at
bounded distance from an isometry of this depth one raft.  From the description of these isometries in
\cite{FarbMosher:sbf},  the only isometry which is at bounded distance from the identity along a
surface subgroup is the identity (this is why there is no Dehn twist type phenomenon here).
Thus $f$ is the identity along the depth one raft $X''$ as well.  Continuing in this way, we see that $f$
is the identity on $X$.  

\end{proof}

\begin{remark}  This result points out some of the ways in which the mapping class group is not well
understood.  Consider any hyperbolic $3$-manifold $M$ which fibers over the circle with fiber $\Sigma$ in two
different ways. If we think of the $M$ as the mapping torus of two different pseudo-Anosov elements $\phi$ and
$\psi$ in $MCG(\Sigma)$ then for every $T$ in $MCG(\Sigma)$ we get an $HNN$-extension by twisting the
identification of one of the fibers with $\Sigma$ by $T$. The depth one rafts are then the surface-by-free
groups corresponding to the homomorphism $F_2 \to MCG(\Sigma)$ sending the generators to $\phi$ and $T^{-1}
\psi T$.   By analogy with classical Schottky groups, one expects that this is a Schottky subgroup of
$MCG(\Sigma)$ in the sense of \cite{FarbMosher:quasiconvex}, for any $T$ which moves the axis of $\psi$
sufficiently far from the axis of $\phi$.  This would give infinitely many examples to which the theorem
applies for every such $M$.  However, this result does not follow directly from the constructions of Schottky
groups given in \cite{FarbMosher:quasiconvex}, nor from the paper \cite{Mosher:hypbyhyp} on which the
construction depends; nevertheless we expect it will follow using similar ideas as in
\cite{Mosher:hypbyhyp}.


\end{remark}

\vfill\break

 \bibliographystyle{amsalpha} %
 \bibliography{mosher}        %

\providecommand{\bysame}{\leavevmode\hbox to3em{\hrulefill}\thinspace}
\providecommand{\MR}{\relax\ifhmode\unskip\space\fi MR }
\providecommand{\MRhref}[2]{%
  \href{http://www.ams.org/mathscinet-getitem?mr=#1}{#2}
}
\providecommand{\href}[2]{#2}
\begin{thebibliography}{GMRS98}

\bibitem[BB00]{BradyBridson:Isoperimetric}
N.~Brady and M.~R. Bridson, \emph{There is only one gap in the isoperimetric
  spectrum}, Geom. Funct. Anal. \textbf{10} (2000), no.~5, 1053--1070.

\bibitem[BK90]{BassKulkarni}
H.~Bass and R.~Kulkarni, \emph{Uniform tree lattices}, Jour. AMS \textbf{3}
  (1990), no.~4, 843--902.

\bibitem[BK98]{BuragoKleiner:nets}
D.~Burago and B.~Kleiner, \emph{Separated nets in {E}uclidean space and
  {J}acobians of bi-{L}ipschitz maps}, Geom. Funct. Anal. \textbf{8} (1998),
  no.~2, 273--282.

\bibitem[Bro82]{Brown:cohomology}
K.~Brown, \emph{Cohomology of groups}, Graduate Texts in Math., vol.~87,
  Springer, 1982.

\bibitem[BW97]{BlockWeinberger}
J.~Block and S.~Weinberger, \emph{Large scale homology theories and geometry},
  Geometric Topology (Athens, GA, 1993), AMS/IP Stud. in Adv. Math., vol. 2.1,
  Amer. Math. Soc., 1997, pp.~522--569.

\bibitem[Cas]{Casson:ThreeDimensional}
A.~Casson, \emph{Three-dimensional topology}, lecture notes.

\bibitem[CC92]{CannonCooper}
J.~Cannon and D.~Cooper, \emph{A characterization of cocompact hyperbolic and
  finite-volume hyperbolic groups in dimension three}, Trans. AMS \textbf{330}
  (1992), 419--431.

\bibitem[CK00]{CrokeKleiner:IdealBoundaries}
C.~Croke and B.~Kleiner, \emph{Spaces with nonpositive curvature and their
  ideal boundaries}, Topology \textbf{39} (2000), no.~3, 549--556.

\bibitem[Dav98]{Davis:CoxeterCohomology}
M.~W. Davis, \emph{{The cohomology of a Coxeter group with group ring
  coefficients}}, Duke Math. J. \textbf{91} (1998), 297--314.

\bibitem[DS00]{DunwoodySwenson:torus}
M.~J. Dunwoody and E.~L. Swenson, \emph{The algebraic torus theorem}, Invent.
  Math. \textbf{140} (2000), no.~3, 605--637.

\bibitem[Dun85]{Dunwoody:Accessible}
M.~J. Dunwoody, \emph{The accessibility of finitely presented groups}, Invent.
  Math. \textbf{81} (1985), 449--457.

\bibitem[FM98]{FarbMosher:BSOne}
B.~Farb and L.~Mosher, \emph{A rigidity theorem for the solvable
  {Baumslag-Solitar} groups}, Invent. Math. \textbf{131} (1998), no.~2,
  419--451, Appendix by D. Cooper.

\bibitem[FM99]{FarbMosher:BSTwo}
\bysame, \emph{Quasi-isometric rigidity for the solvable {Baumslag-Solitar}
  groups, {II}}, Invent. Math. \textbf{137} (1999), no.~3, 613--649.

\bibitem[FM00]{FarbMosher:ABC}
\bysame, \emph{On the asymptotic geometry of abelian-by-cyclic groups}, Acta
  Math. \textbf{184} (2000), no.~2, 145--202.

\bibitem[FM02a]{FarbMosher:quasiconvex}
\bysame, \emph{Convex cocompact subgroups of mapping class groups}, Geometry
  and Topology \textbf{6} (2002), 91--152.

\bibitem[FM02b]{FarbMosher:sbf}
\bysame, \emph{The geometry of surface-by-free groups}, Geom. Funct. Anal.
  \textbf{12} (2002), 915--963.

\bibitem[FS96]{FarbSchwartz}
B.~Farb and R.~Schwartz, \emph{The large-scale geometry of {Hilbert} modular
  groups}, J. Diff. Geom. \textbf{44} (1996), no.~3, 435--478.

\bibitem[Ger93]{Gersten:dimension}
S.~M. Gersten, \emph{Quasi-isometry invariance of cohomological dimension},
  Comptes Rendues Acad. Sci. Paris S\'{e}rie 1 Math. \textbf{316} (1993),
  411--416.

\bibitem[GMRS98]{GMRS:widths}
R.~Gitik, M.~Mitra, E.~Rips, and M.~Sageev, \emph{Widths of subgroups}, Trans.
  AMS \textbf{350} (1998), no.~1, 321--329.

\bibitem[Gro82]{Gromov:Volume}
M.~Gromov, \emph{Volume and bounded cohomology}, IHES Sci. Publ. Math. (1982),
  no.~56, 5--99 (1983).

\bibitem[GW77]{GraverWatkins}
J.~E. Graver and M.~E. Watkins, \emph{Combinatorics with emphasis on the theory
  of graphs}, Springer-Verlag, New York, 1977, Graduate Texts in Mathematics,
  Vol. 54.

\bibitem[Hak61]{Haken:normal}
W.~Haken, \emph{Theorie der {N}ormalflaechen}, Acta Math. \textbf{105} (1961),
  245--375.

\bibitem[Hig61]{Higman:subgroups}
G.~Higman, \emph{Subgroups of finitely presented groups}, Proc. Roy. Soc. Ser.
  A \textbf{262} (1961), 455--475.

\bibitem[HV50]{HalmosVaughan}
P.~R. Halmos and H.~E. Vaughan, \emph{The marriage problem}, Amer. J. Math.
  \textbf{72} (1950), 214--215.

\bibitem[JR88]{JacoRubinstein:PLminimal}
W.~Jaco and J.~H. Rubinstein, \emph{{PL} minimal surfaces in 3-manifolds}, J.
  Diff. Geom. \textbf{27} (1988), no.~3, 493--524.

\bibitem[KK99]{KapovichKleiner:duality}
M.~Kapovich and B.~Kleiner, \emph{Coarse {A}lexander duality and duality
  groups}, preprint, arXiv:math.GT/9911003, 1999.

\bibitem[KK00]{KapovichKleiner:LowDBoundaries}
\bysame, \emph{Hyperbolic groups with low-dimensional boundary}, Ann. Sci.
  \'Ecole Norm. Sup. (4) \textbf{33} (2000), no.~5, 647--669.

\bibitem[KL97]{KapovichLeeb:haken}
M.~Kapovich and B.~Leeb, \emph{Quasi-isometries preserve the geometric
  decomposition of {H}aken manifolds}, Invent. Math. \textbf{128} (1997),
  no.~2, 393--416.

\bibitem[Man05]{Manning:pseudocharacters}
J.~Manning, \emph{Geometry of pseudocharacters}, Geometry and Topology
  \textbf{9} (2005), 1147--1185.

\bibitem[Mar81]{Margulis:amalgams}
G.~A. Margulis, \emph{On the decomposition of discrete subgroups into
  amalgams}, Selecta Math. Soviet. \textbf{1} (1981), no.~2, 197--213.

\bibitem[McM98]{McMullen:lipschitz}
C.~T. McMullen, \emph{Lipschitz maps and nets in {E}uclidean space}, Geom.
  Funct. Anal. \textbf{8} (1998), no.~2, 304--314.

\bibitem[Mos97]{Mosher:hypbyhyp}
L.~Mosher, \emph{A hyperbolic-by-hyperbolic hyperbolic group}, Proc. AMS
  \textbf{125} (1997), no.~12, 3447--3455.

\bibitem[MS04]{MonodShalom:superrigidity}
N.~Monod and Y.~Shalom, \emph{Cocycle superrigidity and bounded cohomology for
  negatively curved spaces}, J. Diff. Geom. \textbf{67} (2004), no.~3,
  395--455.

\bibitem[MSW00]{MSW:QTannc}
L.~Mosher, M.~Sageev, and K.~Whyte, \emph{Quasi-actions on trees, research
  announcement}, preprint, \textsc{arXiv:math.GR/0005210}, 2000.

\bibitem[MSW02]{MSW:MaximallySymmetric}
\bysame, \emph{Maximally symmetric trees}, Geometriae Dedicata \textbf{92}
  (2002), 195--233.

\bibitem[MSW03]{MSW:QTOne}
\bysame, \emph{{Quasi-actions on trees I: Bounded valence}}, Ann. of Math.
  (2003), 116--154.

\bibitem[MY82]{MeeksYau:Equivariant}
W.~H. Meeks and S.~T. Yau, \emph{The equivariant {D}ehn's lemma and loop
  theorem}, Comment. Math. Helv. \textbf{56} (1982), 225--239.

\bibitem[Nib04]{Niblo:ends}
G.~Niblo, \emph{A geometric proof of {S}tallings' theorem on groups with more
  than one end}, Geometriae Dedicata \textbf{105} (2004), 61--76.

\bibitem[Pap02]{Papasoglu:GroupSplittings}
P.~Papasoglu, \emph{{Group splittings and asymptotic topology}}, preprint,
  \textsc{arXiv:math.GR/0201312}, 2002.

\bibitem[Pap05]{Papasoglu:Zsplittings}
\bysame, \emph{Quasi-isometry invariance of group splittings}, Ann. of Math.
  \textbf{161} (2005), no.~2, 759--830.

\bibitem[PW02]{PapasogluWhyte:ends}
P.~Papasoglu and K.~Whyte, \emph{Quasi-isometries between groups with
  infinitely many ends}, Comment. Math. Helv. \textbf{77} (2002), no.~1,
  133--144.

\bibitem[RS94]{RipsSela:structure}
E.~Rips and Z.~Sela, \emph{Structure and rigidity in hyperbolic groups. {I}},
  Geom. Funct. Anal. \textbf{4} (1994), no.~3, 337--371.

\bibitem[Sch97]{Schwartz:Symmetric}
R.~Schwartz, \emph{Symmetric patterns of geodesics and automorphisms of surface
  groups}, Invent. Math. \textbf{128} (1997), 177--199.

\bibitem[Sel97]{Sela:RigidityII}
Z.~Sela, \emph{{Structure and rigidity in (Gromov) hyperbolic groups and
  discrete groups in rank $1$ Lie groups. II}}, Geom. Funct. Anal. \textbf{7}
  (1997), no.~3, 561--593.

\bibitem[Ser80]{Serre:trees}
J.~P. Serre, \emph{Trees}, Springer, New York, 1980.

\bibitem[Sta68a]{Stallings:Dimension1}
J.~Stallings, \emph{Groups of dimension 1 are locally free}, Bull. Amer. Math.
  Soc. \textbf{74} (1968), 361--364.

\bibitem[Sta68b]{Stallings:ends}
\bysame, \emph{On torsion free groups with infinitely many ends}, Ann. of Math.
  \textbf{88} (1968), 312--334.

\bibitem[SW79]{ScottWall}
P.~Scott and C.~T.~C. Wall, \emph{Topological methods in group theory},
  Homological group theory, Proceedings of Durham symposium, Sept. 1977, London
  Math. Soc. Lecture Notes, vol.~36, 1979, pp.~137--203.

\bibitem[Why02]{Whyte:bs}
K.~Whyte, \emph{{The large scale geometry of the higher Baumslag-Solitar
  groups}}, Geom. Funct. Anal. \textbf{11} (2002), 1327--1353.

\bibitem[Why04a]{Whyte:FiberedGeometries}
\bysame, \emph{Geometries which fiber over trees}, in preparation, 2004.

\bibitem[Why04b]{Whyte:Tubular}
\bysame, \emph{Tubular groups}, in preparation, 2004.

\end{thebibliography}

\bigskip

\bigskip\noindent
\textsc{\scriptsize
Lee Mosher:\\
Department of Mathematics,
Rutgers University,
Newark, NJ 07102\\
mosher@andromeda.rutgers.edu
}

\bigskip\noindent
\textsc{\scriptsize
Michah Sageev:\\
Technion, Israel University of Technology,
Dept.\ of Mathematics,
Haifa 32000, Israel\\
sageevm@techunix.technion.ac.il
}

\bigskip\noindent
\textsc{\scriptsize
Kevin Whyte:\\
Dept.\ of Mathematics,
University of Illinois at Chicago
Chicago, IL 60607\\
kwhyte@math.uic.edu
}

\end{document}